%% file: hippylib.tex
\documentclass[format=acmsmall, screen=true, review=false]{acmart}

\input{defjcp14.txt}

\newcommand{\figureref}{Fig. }

\pgfplotsset{every axis/.append style={
    tick label style={font=\Large},
    label style={font=\Large},
    legend style={font=\large}},
    every axis plot/.append style={ultra thick}
}
  
\begin{document}

\acmJournal{TOMS}

\title[\bhip: a Software Framework for Large-Scale Inverse
Problems Governed by PDEs]{
\bhip: An Extensible Software Framework for Large-Scale Inverse
Problems Governed by PDEs}
\subtitle{Part I: Deterministic Inversion and
Linearized Bayesian Inference}

\author{Umberto Villa}
\orcid{0000-0002-5142-2559}
\affiliation{
\institution{Washington University in St. Louis}
\department{Electrical \& Systems Engineering}
\city{St. Louis}
\state{MO}
\country{USA}
}
\email{uvilla@wustl.edu}

\author{Noemi Petra}
\orcid{0000-0002-9491-0034}
\affiliation{
\institution{University of California, Merced}
\department{Applied Mathematics, School of Natural Sciences}
\city{Merced}
\state{CA}
\country{USA}
}
\email{npetra@ucmerced.edu}

\author{Omar Ghattas}
\orcid{}
\affiliation{
\institution{The University of Texas at Austin}
\department{Oden Institute for Computational Engineering \& Sciences}
\department{Department of Mechanical Engineering}
\department{Department of Geological Sciences}
\city{Austin}
\state{TX}
\country{USA}
}

\ccsdesc[500]{Mathematics of computing~Bayesian computation}
\ccsdesc[500]{Computing methodologies~Uncertainty quantification}
\ccsdesc[500]{Mathematics of computing~Mathematical optimization}
\ccsdesc[500]{Mathematics of computing~Partial differential equations}
\ccsdesc[300]{Applied computing~Physical sciences and engineering}
\ccsdesc[300]{Mathematics of computing~Computations on matrices}
\ccsdesc[300]{Mathematics of computing~Discretization}
\ccsdesc[300]{Mathematics of computing~Solvers}

\keywords{
Infinite-dimensional inverse problems, adjoint-based methods, inexact
Newton-CG method, low-rank approximation, Bayesian inference,
uncertainty quantification, sampling, generic PDE toolkit
}
  
\begin{abstract}
We present an extensible software framework, \hip, for solution of large-scale deterministic and Bayesian inverse problems governed by partial differential equations (PDEs) with (possibly) infinite-dimensional parameter fields (which are high-dimensional after discretization). \hip overcomes the prohibitively expensive nature of Bayesian inversion for this class of problems by implementing state-of-the-art scalable algorithms for PDE-based inverse problems that exploit the structure of the underlying operators, notably the Hessian of the log-posterior. The key property of the algorithms implemented in \hip is that the solution of the inverse problem is computed at a cost, measured in linearized forward PDE solves, that is independent of the parameter dimension. The mean of the posterior is approximated by the MAP point, which is found by minimizing the negative log-posterior with an inexact matrix-free Newton-CG method. The posterior covariance is approximated by the inverse of the Hessian of the negative log posterior evaluated at the MAP point. The construction of the posterior covariance is made tractable by invoking a low-rank approximation of the Hessian of the log-likelihood. Scalable tools for sample generation are also discussed. \hip makes all of these advanced algorithms easily accessible to domain scientists and provides an environment that expedites the development of new algorithms. 
\end{abstract}

\maketitle

\section{Introduction}

Recent years have seen tremendous growth in the volumes of
observational and experimental data that are being collected, stored,
processed, and analyzed. The central question that has emerged is:
How do we extract knowledge and insight from all of this data?
When the data correspond
to observations of (natural or engineered) systems, and these systems
can be represented by mathematical models, this knowledge-from-data
problem is fundamentally a mathematical inverse problem.  That is,
given (possibly noisy) data and (a possibly uncertain) model, the goal becomes to infer
parameters that characterize the model.  Inverse problems abound in
all areas of science, engineering, technology, and medicine. As 
just a few examples of model-based inverse problems, we may
infer: the initial condition in a time-dependent partial differential equation (PDE) model, a
coefficient field in a subsurface flow model, the ice sheet basal
friction field from satellite observations of surface flow, the earth
structure from reflected seismic waves, subsurface contaminant plume
spread from crosswell electromagnetic measurements, internal
structural defects from measurements of structural vibrations, ocean
state from surface temperature observations, and so on.

Typically, inverse problems are {\em ill-posed} and suffer from
non-unique solutions; simply put, the data---even when they are
large-scale---do not provide sufficient information
to fully determine the model
parameters. This is the usual case with PDE models that have parameters representing fields such as boundary
conditions, initial conditions, source terms, or heterogeneous
coefficients. Non-uniqueness can stem from noise in the data or model,
from sparsity of the data, from smoothing properties of the map from
input model parameters to output observables or from its nonlinearity,
or from intrinsic redundancy in the
data.  In such
cases, {\em uncertainty is a fundamental feature of the inverse
  problem}. Therefore, not only do we wish to infer the parameters, but we must
also quantify the uncertainty associated with this inference,
reflecting the degree of ``confidence'' we have in the solution.

Methods that facilitate the solution of Bayesian inverse problems
governed by complex PDE models require a diverse and advanced background
in applied mathematics, scientific computing, and statistics to
understand and implement, e.g., Bayesian inverse theory, computational
statistics, inverse problems in function space, adjoint-based first-
and second-order sensitivity analysis, and variational discretization
methods.  In addition, to be efficient these methods generally require first and
second derivative (of output observables with respect to input parameters)
information from the underlying forward PDE model,
which can be cumbersome to derive. In this paper, we present \hip, an
{\bf I}nverse {\bf P}roblems {\bf Py}thon {\bf lib}rary (\hip),
an extensible software framework aimed at overcoming these
challenges and providing capabilities for additional algorithmic
developments for large-scale deterministic and Bayesian inversion.

\hip builds on \Fe (a parallel finite element element
library)~\cite{LoggMardalWells12} for the discretization of the PDEs,
and on~\pet\cite{BalayAbhyankarAdamsEtAl14} for scalable and efficient
linear algebra operations and solvers. Hence, it is easily applicable
to medium to large-scale problems. One of the main features of this
library is that it clearly displays and utilizes specific aspects from the model setup to
the inverse solution, which can be useful not only for research
purposes but also for learning and teaching.\footnote{\hip is currently used to teach several 
graduate level classes on inverse problems at various universities, including The University of Texas at Austin,
University of California, Merced, Washington University in St. Louis, New York University, and North Carolina State.
\hip has also been demonstrated with hands-on interactive sessions  at workshops and summer schools, such as the 2015 ICERM IdeaLab,
the 2016 SAMSI Optimization Program Summer School, and the 2018 Gene Golub SIAM Summer School.}
In the \hip examples, we
show how to handle various PDE models and boundary conditions, and
illustrate how to implement prior and log-likelihood terms for the
Bayesian inference. \hip is implemented in a mixture of C++ and Python
and has been released under the GNU General Public License version 2
(GPL). The source codes can be downloaded from
\url{https://hippylib.github.io}. Below we summarize the main
algorithmic and software contributions of \hip.

\paragraph{Algorithmic contributions}
\begin{enumerate}
  \item A single pass randomized eigensolver for generalized symmetric
    eigenproblems that is more accurate than the one proposed
    in~\cite{SaibabaLeeKitanidis16}.
  \item A new scalable sampling algorithm for Gaussian random fields
    that exploits the structure of the given covariance operator. This
    extends the approach proposed in~\cite{CrociGilesRognesEtAl18} to
    covariance operators defined as the inverse of second order
    differential operators as opposed to the identity operator.
  \item A scalable algorithm to estimate the pointwise variance of
    Gaussian random fields using randomized eigensolvers. For the same
    computational cost this algorithm allows for more accurate
    estimates than the stochastic estimator proposed
    in~\cite{BekasKokiopoulouSaad07}. Our method drastically reduces
    the variance of the estimator at a cost of introducing a small
    bias.
\end{enumerate}

\paragraph{Software contributions}
\begin{enumerate}
\item A modular approach to define complex inverse problems governed
  by (possibly nonlinear or time-dependent) PDEs. \hip automates the
  computation of higher order derivatives of the
  parameter-to-observable map for forward models and observation
  processes defined by the user through \Fe.
\item Implementation of adjoints and Hessian actions needed to solve
  the deterministic inverse problem and to compute the maximum a
  posteriori (MAP) point of the Bayesian inverse problem. In addition,
  to test gradients and the Hessian action, \hip incorporates finite
  difference tests, which is an essential component of the
  verification process.
\item A robust implementation of the inexact Newton-conjugate gradient
  (Newton-CG) algorithm together with line search algorithms to
  guarantee global convergence of the optimizer.
\item Implementation of randomized algorithms to compute the low-rank
  factorization of the misfit part of the Hessian.
\item Scalable algorithms to construct and evaluate the Laplace
  approximation of the posterior.
\item Sampling capabilities to generate realizations of Gaussian random
  fields with a prescribed covariance operator.
\item An estimation of the pointwise variance of the prior
  distribution and Laplace approximation to the posterior.
\end{enumerate}

Numerous toolkits and libraries for finite element computations based
on variational forms are available, for instance \software{COMSOL
  Multiphysics}~\cite{COMSOL}, \software{deal.II}
\cite{BangerthHartmannKanschat07}, \software{dune}
\cite{BastianBlattDednerEtAl08}, \Fe
\cite{LoggMardalWells12,LangtangenLogg17},
and \software{
  Sundance}, a package from \software{Trilinos}
\cite{HerouxBartlettHowleEtAl05}. While these toolkits are usually
tailored towards the solution of PDEs and systems of PDEs, they cannot
be used straightforwardly for the solution of inverse problems with
PDEs. However, several of them are
sufficiently flexible to be extended for the solution of inverse problems
governed by PDEs.
Nevertheless, some knowledge of the structure
underlying these packages is required since the optimality systems
arising in inverse problems with PDEs often cannot be solved using
generic PDE solvers, which do not exploit the optimization structure
of the inverse problems. In~\cite{FarrellHamFunkeEtAl13} the
authors present \software{dolfin-adjoint}, a project that also builds on
\Fe and derives discrete adjoints from a forward model written in the
Python interface to \software{dolfin} using a combination of symbolic and
automatic differentiation.  While \software{dolfin-adjoint} could be used
to solve deterministic inverse problems, it lacks the
framework for Bayesian inversion. In addition, we avoid using the
adjoint capabilities of \software{dolfin-adjoint} since this does not
allow the user to have full control over the construction of
derivatives. In~\cite{RuthottoTreisterHaber2017} the
authors present \software{jInv}, a flexible parallel
software for parameter estimation with PDE forward models. The main
limitations of this software are that it is restricted to
deterministic inversion and that the user needs to provide the discretization
for both the forward and adjoint problems. Finally, the \software{Rapid Optimization Library},
ROL \cite{KouriRidzalWinckel18}, is a flexible and robust optimization package in
\software{Trilinos} for the solution of optimal design, optimal control and
deterministic inverse problems in large-scale engineering applications.
ROL implements state-of-the-art algorithms for
unconstrained optimization, constrained optimization and optimization under uncertainty, and exposes
an interface specific for optimization problems with PDE constraints. The main limitation is that the user
has to interface with other software packages for the definition and implementation of the forward and adjoint problems.
There also exist several general purpose libraries addressing uncertainty 
quantification (UQ) and Bayesian inverse
problems.  Among the most prominent we mention
\software{QUESO}~\cite{McDougallMalayaMoser17,PrudencioSchulz12},
\software{DAKOTA}~\cite{EldredGiuntaEtAl02,AdamsBohnhoffDalbeyEtAL09}, \software{PSUADE} \cite{Tong17}, 
UQTk \cite{DebusschereSargsyanSaftaEtAl17}.
All of these libraries provide Bayesian inversion capabilities,
but the underlying methods do not fully exploit the structure
of the problem or make use of derivatives and as such are not intended for
high-dimensional problems. Finally, MUQ \cite{MUQ}
provides powerful Bayesian inversion models and algorithms, but
expects forward models to come equipped with gradients/Hessians to
permit large-scale solution.

In summary, to the best of our knowledge, there is no available software
(open-source or otherwise) that provides all the discretization, optimization and
statistical tools to enable scalable and efficient solution of deterministic and Bayesian inverse problems
governed by complex PDE forward models. \hip is the first software framework that allows to tackle
this specific class of inverse problems by facilitating the construction of forward PDE models equipped with adjoint/derivative
information, providing state-of-the-art scalable optimization algorithms for the solution of the deterministic inverse problem and/or MAP point computation, and integrating tools for characterizing the posterior distribution.

The paper is structured as follows. Section \ref{sec:InfDimIP} gives a brief overview of the deterministic and Bayesian formulation
of inverse problems in an infinite-dimensional Hilbert space setting, and addresses the discretization of the underlying PDEs using the finite element method. 
Section \ref{sec:design} contains an overview of the design of the \hip software and of its components.
Section \ref{sec:algorithms} provides a detailed description of the algorithms implemented in \hip to
solve the deterministic and linearized Bayesian inverse problem, namely the inexact Newton-CG algorithm,
the single and double pass randomized algorithms for the solution of generalized hermitian eigenproblems, scalable sampling techniques
for Gaussian random fields, and stochastic algorithms to approximate the pointwise variance of the prior and posterior distributions.
Section \ref{sec:models} demonstrates \hip's capabilities for deterministic and linearized Bayesian inversion by solving two representative inverse problems: inversion for the coefficient field in an elliptic PDE model and for the initial condition in an advection-diffusion PDE model.
Last, Section \ref{sec:conclusions} contains our concluding remarks.

\section{Infinite-dimensional deterministic and Bayesian inverse problems in hIPPYlib}
\label{sec:InfDimIP}

In what follows, we provide a brief account of the deterministic
\cite{EnglHankeNeubauer96,Vogel02} and Bayesian formulation
\cite{Tarantola05, KaipioSomersalo05} of inverse problems.
Specifically, we adopt infinite-dimensional Bayesian inference
framework \cite{Stuart10}, and we refer to
\cite{PetraMartinStadlerEtAl14,
  AlexanderianPetraStadlerEtAl14,
  AlexanderianPetraStadlerEtAl16,
  Bui-ThanhGhattasMartinEtAl13}
for elaborations associated with discretization issues.

\subsection{Deterministic inverse problems governed by PDEs}
\label{sec:det}

The inverse problem consists of using available observations $\obs$ to
infer the values of the unknown parameter field\footnote{
\hip also supports deterministic and Bayesian inversion for a finite-dimensional set of parameters,
however, for ease of notation, in the present work we only present the infinite-dimensional case.}
$\ipar$ that characterize a
physical process modeled by PDEs. Mathematically this inverse
relationship is expressed as
\begin{align}\label{eq:noisemodel}
  \obs = \ff(\ipar) + \vec{\eta},
\end{align}
where the map $\ff: \iparspace \to \mathbb{R}^q$ is the so-called {\it
  parameter-to-observable} map. This mapping
can be linear or nonlinear. In the applications targeted in \hip
$\iparspace \subseteq L^2(\D)$, where $\D \subset \R^d$ is a bounded domain,
and evaluations of $\ff$ involve the solution of a PDE given $m$,
followed by the application of an observation operator to extract the
observations from the state. 
That is, introducing the state variable $\istate \in \mathcal{V}$ for
a suitable Hilbert space $\mathcal{V}$ of functions defined on $\D$, the map $\ff$ is defined as
\begin{equation}\label{eq:fwdmap}
\ff(\ipar) = \mathcal{B}(\istate), \text{ s.t. } r(\istate, \ipar) = 0,
\end{equation}
where $\mathcal{B}: \mathcal{V} \to \mathbb{R}^q$ is a (possibly nonlinear) observation operator,
and $r: \mathcal{V}\times\iparspace \to \mathcal{V}^*$---referred as the forward problem from now on---represents the PDE problem.
The observations $\obs$ contain noise due
to measurement uncertainties and model
errors~\cite{Tarantola05}. In~\eqref{eq:noisemodel}, this is captured
by the additive noise $\vec{\eta}$, which in \hip is modeled as $\vec{\eta} \sim
\GM{\vec{0}}{\ncov}$, i.e., a centered Gaussian at $\vec{0}$ with covariance
$\ncov$.
A significant difficulty when solving infinite-dimensional inverse problems
is that typically these are not well-posed (in the sense of
Hadamard~\cite{TikhonovArsenin77}). To overcome the difficulties due
to ill-posedness, we regularize the problem, i.e., we include
additional assumptions on the solution, such as smoothness. The
deterministic inverse problems in \hip are regularized via Tikhonov
regularization, which penalizes oscillatory components of the
parameter $m$, thus restricting the solution to smoothly varying
fields~\cite{EnglHankeNeubauer96, Vogel02}.

A deterministic inverse problem is therefore formulated as follows:
given finite-dimensional noisy observations $\obs \in \R^q$, one seeks
to find the unknown parameter field $\ipar$ that best reproduces the
observations. Mathematically this translates into the following
nonlinear least-squares minimization problem
\begin{equation}\label{eq:optpb}
  \min_{\ipar \in \iparspace} \! \mc{J}\LRp{\ipar} \equaldef
  \half \nor{\ff(\ipar) - \obs
  }^2_{\ncov^{-1}} + \mathcal{R}(\ipar),
\end{equation}
where the first term in the cost functional, $\mc{J}\LRp{\ipar}$,
represents the misfit between the observations, $\obs$, and that
predicted by the parameter-to-observable map $\ff(\ipar)$, weighted by
the inverse noise covariance $\ncov^{-1}$. The regularization term,
$\mathcal{R}(\ipar)$, imposes regularity on the inversion field
$\ipar$, such as smoothness. As explained above, in the absence of
such a term, the inverse problem is ill-posed, i.e., its solution is
not unique and is highly sensitive to errors in the
observations~\cite{EnglHankeNeubauer96, Vogel02}.

As we will explain in Section \ref{subsec:NCG}, to efficiently solve the nonlinear least-squares problem \eqref{eq:optpb} with parameter-to-observable-map $\ff$ implicitly defined as in \eqref{eq:fwdmap}, first and second derivative information are needed.
Using the  Lagrangian formalism \cite{Troltzsch10}, abstract expressions for the gradient and Hessian action are obtained below, and we refer to Section \ref{sec:models} for concrete examples. To this aim, we introduce an auxiliary variable $\iadj \in \mathcal{V}$, from here on referred as the adjoint, and write the Lagrangian functional
\begin{equation*}
\mathcal{L}^\G(\istate, \ipar, \iadj) := \frac{1}{2}\| \mathcal{B}(u) - \obs \|^2_{\ncov^{-1}} + \mathcal{R}(\ipar) + \,_{\mathcal{V}}\langle\iadj,  r(\istate, \ipar) \rangle_{\mathcal{V}^*},
\end{equation*}
where $_{\mathcal{V}}\langle \cdot, \cdot \rangle_{\mathcal{V}^*}$ denotes the duality pair between $\mathcal{V}$ and its adjoint. The gradient for the cost functional \eqref{eq:optpb} in an arbitrary direction $\tilde{m}\in\iparspace$ evaluated at $\ipar = \ipar_0 \in \iparspace$ is the  G\^ateaux derivative of $\mathcal{L}$ with respect to $\ipar$,
and reads
\begin{equation}\label{eq:abstract_gradient}
\left( \G(\ipar_0), \tilde{\ipar}  \right)=
\left( \mathcal{R}_{\ipar}(\ipar_0), \tilde{\ipar} \right) 
+ \,_{\mathcal{V}}\langle \iadj_0, r_{\ipar}(\istate_0, \ipar_0)[\tilde{\ipar}]\rangle_{\mathcal{V}^*}, \quad \forall \tilde{\ipar} \in \iparspace,
\end{equation}
where $\left( \mathcal{R}_{\ipar}(\ipar_0), \tilde{\ipar} \right) \in \mathbb{R}$ denotes the G\^ateaux derivative of $\mathcal{R}$ with respect to $\ipar$ in the direction $\tilde{\ipar}$ evaluated at $\ipar = \ipar_0$, and $
r_{\ipar}(\istate_0, \ipar_0)[\tilde{\ipar}] \in \mathcal{V}^*$ the  G\^ateaux derivative of $r$ with respect to $\ipar$ in the direction $\tilde{\ipar}$ evaluated at $\istate=\istate_0, \ipar = \ipar_0$.
Here $u_0$, $p_0$ are obtained by setting to zero the derivatives of $\mathcal{L}$ with respect to $\iadj$ and $\istate$; specifically, $u_0$ solves the forward problem
\begin{equation}\label{eq:abstract_fwd}
_{\mathcal{V}}\langle \tilde{\iadj},  r(\istate_0, \ipar_0) \rangle_{\mathcal{V}^*} = 0 \quad \forall \tilde{\iadj} \in \mathcal{V},
\end{equation}
and $\iadj_0$ solves the adjoint problem
\begin{equation}\label{eq:abstract_adj}
_{\mathcal{V}}\langle \iadj_0,  r_{\istate}(\istate_0, \ipar_0)[\tilde{\istate}] \rangle_{\mathcal{V}^*}
+ \langle \mathcal{B}_{\istate}(\istate_0)[\tilde{\istate}],  \mathcal{B}(\istate) - \obs \rangle_{\mathbb{R}^q}= 0, \quad \forall \tilde{\istate}  \in \mathcal{V}.
\end{equation}

In a similar way, to derive the expression for the Hessian action in an arbitrary direction $\hat{\ipar} \in \iparspace$ we introduce the second order Lagrangian functional
\begin{equation}\label{eq:abstract_metaL}
\begin{aligned}
\mathcal{L}^\H(\istate, \ipar, \iadj; \hat{\istate}, \hat{\ipar}, \hat{\iadj}) & := \left( \G(\ipar), \hat{\ipar}\right) \\
{} & +
_{\mathcal{V}}\langle \hat{\iadj},  r(\istate, \ipar) \rangle_{\mathcal{V}^*} \\
{} & +
_{\mathcal{V}}\langle \iadj,  r_{\istate}(\istate, \ipar)[\hat{\istate}] \rangle_{\mathcal{V}^*}
+ \langle  \mathcal{B}_{\istate}(\istate)[\hat{\istate} ],  \mathcal{B}(\istate) - \obs \rangle_{\mathbb{R}^q},
\end{aligned}
\end{equation}
where the first term is the gradient expression, the second term stems from the forward problem, and the last two terms represent the adjoint problem. Then, the action of the Hessian in a direction $\hat{m} \in \mathcal{M}$ evaluated at $\ipar = \ipar_0$ is the variation of $\mathcal{L}^\H$ with respect to $\ipar$ and reads
\begin{equation}\label{eq:abstract_hessianAction}
\begin{aligned}
\left(\tilde{\ipar}, \H(\ipar_0)\hat{\ipar} \right) &= \left(\tilde{\ipar},  \mathcal{R}_{\ipar\ipar}(\ipar_0)[\hat{\ipar}]\right) + 
\left(\iadj_0,  r_{\ipar\ipar}(\istate_0, \ipar_0)[\tilde{\ipar}, \hat{\ipar}]\right) \\
{} &+ _{\mathcal{V}}\langle \hat{\iadj},  r_\ipar(\istate_0, \ipar_0)[\tilde{\ipar}] \rangle_{\mathcal{V}^*}
 + _{\mathcal{V}}\langle \iadj_0,  r_{\istate\ipar}(\istate_0, \ipar_0)[\hat{\istate},\tilde{\ipar}] \rangle_{\mathcal{V}^*}, \quad \forall \ipar \in \iparspace.
\end{aligned}
\end{equation}
Here $\istate_0$, $\iadj_0$ are the solution of the forward and adjoint problems \eqref{eq:abstract_fwd} and \eqref{eq:abstract_adj}, respectively. The incremental state $\hat{\istate}$ and incremental adjoint $\hat{\iadj}$ solve the so-called \emph{incremental forward} and \emph{incremental adjoint problems}, which are obtained by setting to zero variations of $\mathcal{L}^\H$ with respect to $\iadj$ and $\istate$, respectively.
In Appendix~\ref{app-sec:MAP}, we present a Newton-type algorithm to minimize \eqref{eq:optpb} that uses the expression for the gradient \eqref{eq:abstract_gradient} and Hessian action  \eqref{eq:abstract_hessianAction} derived here.

Finally, we note that the solution of a deterministic inverse problem based on
regularization is a \emph{point estimate} of $\ipar$, which
solves~\eqref{eq:noisemodel} in a least-squares sense. A systematic integration of the prior
information on the model parameters and uncertainties associated with
the observations can be achieved using a probabilistic point of view,
where the prior information and noise model are represented by
probability distributions. In the following section, we describe the
probabilistic formulation of the inverse problem via a Bayesian
framework, whose solution is a posterior probability distribution for $m$.

\subsection{Bayesian inversion in infinite dimensions}
\label{sec:bayes}
In the Bayesian formulation in infinite dimensions, we state the
inverse problem as a problem of statistical inference over the space
of uncertain parameters, which are to be inferred from data and a
physical model.  In this setup, in contrast to the finite-dimensional
case, there is no Lebesgue measure on $\iparspace$, the
infinite-dimensional Bayes formula is given by
\begin{equation} \label{equ:bayes-abstract}
   \frac{d\mupost}{d\muprior} \propto \like(\obs | \ipar).
\end{equation}
Here, $d\mupost/d\muprior$ denotes the Radon-Nikodym
derivative~\cite{Williams1991} of the posterior measure $\mupost$ with
respect to $\muprior$, and $\like(\obs | \ipar)$ denotes the data
likelihood.  Conditions under which the posterior measure is well
defined and~\eqref{equ:bayes-abstract} holds are given in detail
in~\cite{Stuart10}.

{\it The noise model and the likelihood.} 
In our \hip framework, we assume
an additive noise model, $\obs = \ff(\ipar) + \vec{\eta}$, where
$\vec{\eta} \sim \GM{\vec{0}}{\ncov}$ is a centered Gaussian on
$\R^q$. This implies
\begin{equation}\label{equ:likelihood}
\like(\obs | \ipar) \propto \exp\Big\{-\Phi(m)\Big\},
\end{equation}
where $\Phi(m) = \half \| \ff(\ipar) -
\vec{d} \|^2_{\ncov^{-1}}$ denotes the negative log-likelihood.\\

{\it The prior.} For many problems, it is
reasonable to choose the prior to be Gaussian, i.e., $\ipar \sim
\GM{\iparpr}{\Cprior}$. 
This implies
\begin{equation}\label{equ:prior}
d\muprior(\ipar) \propto \exp\Big\{-\half\| m - \iparpr\|^2_{\Cprior^{-1}}\Big\}.
\end{equation}
If the parameter represents a spatially
correlated field defined on $\D \in \mathbb{R}^d$, the prior
covariance operator $\Cprior$ usually imposes smoothness on the
parameter. This is because rough components of the parameter field are
typically cannot be inferred from the data, and must be determined by the
prior to result in a well-posed Bayesian inverse problem.

In \hip we use elliptic PDE operators to construct
the prior covariance, which allows us to capitalize on fast, optimal
complexity solvers.  More precisely, the prior covariance operator is
the inverse of the $\nu$-th power ($\nu > \frac{d}{2}$) of a
Laplacian-like operator, namely $\Cprior := \Acal^{-\nu} = (-\gamma \,
\gbf \Delta + \delta I)^{-\nu}$, where $\gamma$, and $\delta > 0$
control the correlation length $\rho$ and the pointwise variance $\sigma^2$ of the prior operator. 
Specifically, $\rho$---empirically defined as the distance $\rho$ for which the two-points correlation coefficient is 0.1---is
proportional to $\sqrt{\gamma/\delta}$, and $\sigma^2$ is proportional to $\delta^{-\nu}\rho^{-d}$
(see e.g. \cite{LindgrenRueLindstroem11} where exact expressions for $\rho$ and $\sigma^2$ as functions of $\gamma$ and $\delta$
are derived under the assumption of unbounded domain $\D$ and constant coefficients $\gamma$ and $\delta$).
The coefficients $\gamma$ and~$\delta$ can be constant
(in which case the prior is stationary) or spatially varying.  In addition, one can consider
an anisotropic diffusion operator $\Acal = -\gamma\,\nabla \cdot
(\boldsymbol{\Theta}\, \nabla) + \delta I$, with $\boldsymbol{\Theta}$
a symmetric positive definite (s.p.d.) tensor that models, for instance, stronger correlations in a
specific direction.  These choices of prior ensure that $\Cprior$ is a
trace-class operator, guaranteeing bounded pointwise variance and a
well-posed infinite-dimensional Bayesian inverse problem
\cite{Stuart10, Bui-ThanhGhattasMartinEtAl13}.

{\it The posterior.} Using the expression for the likelihood function \eqref{equ:likelihood}
and prior distribution \eqref{equ:prior}, the posterior distribution in \eqref{equ:bayes-abstract} reads
\begin{equation}\label{equ:post-infdim}
d\mupost \propto \exp\Big\{-\half \| \ff(\ipar) -
\vec{d} \|^2_{\ncov^{-1}} -\half\| m - \iparpr\|^2_{\Cprior^{-1}} \Big\}.
\end{equation}
The \emph{maximum a posteriori} (MAP) point $\iparmap$ is defined as the parameter field that
maximizes the posterior distribution. It can be obtained
by solving the following deterministic optimization problem
\begin{equation}\label{eq:map_continuum}
\iparmap := \argmin_{\ipar \in \iparspace} (- \log d\mupost(\ipar) ) = \argmin_{\ipar \in \iparspace}
\half \| \ff(\ipar) -
\vec{d} \|^2_{\ncov^{-1}}  + \half\| m - \iparpr\|^2_{\Cprior^{-1}}.
\end{equation}
We note that, the prior information plays the role of Tikhonov regularization in \eqref{eq:optpb}; in fact the deterministic optimization problem \eqref{eq:optpb} is the same as \eqref{eq:map_continuum} for the choice $\R(\ipar) = \half\| m - \iparpr\|^2_{\Cprior^{-1}}$.
The Hessian $\mathcal{H}(\iparmap)$ of the negative log-posterior evaluated at $\iparmap$
plays a fundamental role in quantifying the uncertainty in the inferred parameter. In particular, this indicates which directions in the parameter space are most informed by the data \cite{Bui-ThanhGhattasMartinEtAl13}.
We note that when $\iFF$ is linear, due to the particular choice of
prior and noise model, the
posterior measure is Gaussian, $\GM{\iparmap}{\Cpost}$
with~\cite[Section 6.4]{Stuart10},
\begin{equation}\label{equ:mean-cov}
\Cpost = \mathcal{H}^{-1} = (\iFFadj \ncov^{-1} \iFF + \Cprior^{-1})^{-1}, 
\qquad 
\iparmap = \Cpost(\iFFadj\ncov^{-1}\obs + \Cprior^{-1}\iparpr),
\end{equation}
where $\iFFadj:\R^q \to \iparspace$ is the adjoint of $\iFF$.

In the general case of nonlinear parameter-to-observable map $\iFF$ the
posterior distribution is not Gaussian. However, under certain assumptions
on the noise covariance $\ncov$, the number $q$ of observations, and the
regularity of the parameter-to-observable map $\iFF$, the \emph{Laplace approximation}
\cite{Stigler86,Wong01,EvansSwartz00,Press03,TierneyKadane86} can be invoked to estimate
posterior expectations of functionals of the parameter $\ipar$.
Specifically, assuming that the negative log-likelihood $\Phi(\ipar)$ is strictly convex in a neighborhood of 
$\iparmap$\footnote{To guarantee a positive definite posterior
covariance operator also in the case of non-locally convex negative log-likelihood $\Phi(\ipar)$,
the inverse of the Gauss-Newton Hessian of the negative log-posterior can be used instead.
This corresponds to linearizing the parameter-to-observable map $\iFF$ around $\iparmap$.},
the Laplace approximation to the posterior constructs a Gaussian
distribution $\mulaplace$, 
\begin{equation}\label{eq:laplace_approx}
\mulaplace \sim \GM{\iparmap}{\mathcal{H}(\iparmap)^{-1}},
\end{equation}
centered at $\iparmap$ and with covariance operator
\begin{equation}\label{eq:laplace_approx_cov}
\mathcal{H}(\iparmap)^{-1} = (\mathcal{H}_{\mbox{\tiny misfit}}(\iparmap) + \Cprior^{-1})^{-1}.
\end{equation}
Here $\mathcal{H}_{\mbox{\tiny misfit}}$ denotes the Hessian of the negative log-likelihood evaluated at
$\iparmap$ (see Section \ref{sec:models} for examples of the derivation of the action of $\mathcal{H}_{\mbox{\tiny misfit}}$ using
variational calculus and Lagrangian formalism).

The \emph{Laplace approximation} above is an important tool in designing
scalable and efficient methods for Bayesian inference and UQ implemented in \hip.
It has been studied in the context of PDE-based inverse problems
to draw approximate samples and compute approximate statistics (such as the pointwise variance)
in \cite{Bui-ThanhGhattasMartinEtAl13}. Likewise, it has been exploited in \cite{PetraMartinStadlerEtAl14}
 to efficiently explore the true posterior distribution by generating high quality proposals
for Markov chain Monte Carlo algorithms, in \cite{CuiMartinMarzoukEtAl14a} to construct
likelihood informed subspaces that allows for optimal dimension reduction in Bayesian inference problems,
and in \cite{SchillingsSchwab16,ChenVillaGhattas17} to construct a dimension independent
sparse grid to evaluate posterior expectations. It has also been invoked in \cite{IsaacPetraStadlerEtAl15} for
scalable approximation of the predictive posterior distribution of a scalar quantity of interest. Finally,
its use was advocated
in \cite{LongScavinoTemponeEtAl13,LongScavinoTemponeEtAl15,LongMotamedTempone15,AlexanderianPetraStadlerEtAl16}
to approximate the solution of Bayesian optimal experimental design problems.

\subsection{Discretization of the Bayesian inverse problem}
\label{subsec:discretebayesianform}
We present a brief discussion of the finite-dimensional approximations
of the prior and the posterior distributions; a lengthier discussion
can be found in \cite{Bui-ThanhGhattasMartinEtAl13}.  We start with a
finite-dimensional subspace $\iparspace_h$ of $\iparspace \subseteq L^2(\D)$
originating from a
finite element discretization with continuous Lagrange basis functions
$\LRc{\phi_j}_{j=1}^n$~\cite{BeckerCareyOden81,StrangFix88}. The
approximation of the inversion parameter function $\ipar \in
\iparspace$ is then $\ipar_h = \sum_{j=1}^nm_j\phi_j \in \iparspace_h$, 
and, in what follows, $\bs m =\LRp{m_1,\hdots,m_n}^T\in \R^n$ denotes
the vector of the coefficients in the finite element expansion of $\ipar_h$.

The
finite-dimensional space $\iparspace_h$ inherits the $L^2$-inner product. Thus,
inner products between nodal coefficient vectors must be weighted by a
mass matrix $\M\in \R^{n\times n}$ to approximate the
infinite-dimensional $L^2$-inner product.  This $\M$-weighted inner
product is denoted by $\mip{\cdot\,}{\cdot}$, where
$\mip{\vec{y}}{\vec{z}} = \vec{y}^T \M \vec{z}$ and $\M$ is the
(symmetric positive definite) mass matrix
\[
M_{ij} = \int_\D \phi_i(\x) \phi_j(\x) \, d\x ~, \quad i,j =
1,\hdots,n.
\]
To distinguish $\R^n$ equipped with the $\M$-weighted inner product
with the usual Euclidean space $\R^n$, we denote it by $\R^n_\M$.  For
an operator $\BB: \R^n_{\M} \rightarrow \R^n_{\M}$, we denote the
matrix transpose by $\BB^T$ with entries $(B^T)_{ij} = B_{ji}$.  In
contrast, the $\M$-weighted inner product adjoint $\BB^*$ satisfies,
for $\vec y, \vec z\in \R^n$,
\[
\mip{\BB\vec y}{\vec z} = \mip{\vec y}{\BB^* \vec z},
\]
which implies that $\BB^*$ is given by
\begin{align}
  \label{eq:MM_adj}
  \Badj &= \M^{-1} \BB^T \M.
\end{align}
With these definitions, the matrix representation of the bilinear form
involving the elliptic PDE operator $\Acal^n$ defined in
Section~\ref{sec:bayes} is given by $\RR$ whose components are
\begin{align}
R_{ij} = \int_\D \phi_i(\x) \Acal^\nu \phi_j(\x) \, d\x, \quad i,j
\in \LRc{1,\hdots,n}.
\end{align}
Finally, restating Bayes' theorem with Gaussian noise and prior in
finite dimensions, we obtain:
\begin{equation}
\post(\vec{m})
\propto\exp\Big(
-\frac{1}{2}\|\FF(\vec{m}) - \vec{d}_{\text{obs}} \|^2_{\ncov^{-1}}
-\frac{1}{2}\|\vec{m}-\mpr\|^2_{\prcov^{-1}}
\Big),\label{posterior}
\end{equation}
where $\mpr$ is the mean of the prior distribution,
$\prcov := \RR^{-1}\M\in \mathbb{R}^{n\times n}$ is the
covariance matrix for the prior that arises upon discretization of
$\Cprior$, and $\ncov\in \mathbb{R}^{q \times
  q}$ is the covariance matrix for the noise.
The method of choice to
explore the full posterior is Markov chain Monte Carlo (MCMC), which 
samples the posterior so that sample statistics can be
computed.  MCMC for large-scale inverse problems is still prohibitive
for expensive forward problems and high-dimensional parameter spaces;
here we make a quadratic approximation of the negative log
of the posterior \eqref{posterior}, which results---as discussed in Section
\ref{sec:bayes} in the continuous setting---in the Laplace
approximation of the posterior given by
\begin{equation}\label{eq:lapl_discrete}
\post(\vec{\ipar}) \propto \mathcal{N}(\vec{\map},\postcov)
\end{equation}
The mean of this approximate posterior distribution, $\vec{\map}$, is the
parameter vector maximizing the posterior \eqref{posterior}, and is
known as the {\em maximum a posteriori} (MAP) point.  It can be found
by minimizing the negative log-posterior\footnote{For simplicity,
we assume that the negative log-posterior has
a unique minimum. In general, the negative log-posterior is not guaranteed to
be convex and may admit multiple minima; in this case domain specific techniques should be exploited
to locate the global minimum.},
which amounts to solving the
following optimization problem:
\begin{equation}\label{eq:objfunction-bayesian}
\vec{\map} := \argmin_{\vec{m}} (- \log \post(\vec{m})) = \argmin_{\vec{m}} \frac{1}{2}\|\FF(\vec{m}) - \vec{d}_{\text{obs}} \|^2_{\ncov^{-1}}
+ \frac{1}{2}\|\vec{m}-\mpr\|^2_{\prcov^{-1}},
\end{equation}
which is the discrete counterpart of problem \eqref{eq:map_continuum}.
Denoting with $\tH(\vec{\map})$ and $\tHmisfit(\vec{\map})$ the matrix representations of, respectively, the second derivative  of negative log-posterior $\mathcal{H}$ and log-likelihood $\mathcal{H}_{\mbox{\tiny misfit}}$ (i.e., the data misfit
component of the Hessian) in the $\M$-weighted inner product, and assuming that
$\tHmisfit(\vec{\map})$ is positive definite,
the covariance matrix $\postcov$ in the Laplace approximation is given by 
\begin{equation}
\postcov = \tH^{-1} (\vec{\map})= \left(\tHmisfit(\vec{\map}) + \prcov^{-1} \right)^{-1}.
\label{eq:posterior}
\end{equation}
For
simplicity of the presentation, in the following we will let $\mat{H} = \M \tH$ and $\Hmisfit
= \M \tHmisfit$ be the matrix representation of the Hessian with
respect to the standard Euclidean inner product. Using this notation, and recalling that
$\prcov = \mat{R}^{-1}\M$, we rewrite \eqref{eq:posterior} as
\begin{equation}
\postcov = \mat{H}^{-1} (\vec{\map})\,\M = \left(\Hmisfit(\vec{\map}) + \mat{R} \right)^{-1}\, \M.
\label{eq:posteriorEuclidean}
\end{equation}

Equation \eqref{eq:objfunction-bayesian} and \eqref{eq:posteriorEuclidean} define the
mean and covariance matrix of the Laplace approximation to the posterior in
the discrete setting. In Section \ref{sec:algorithms}, we present scalable
(with respect to the parameter dimension) algorithms to compute the discrete
MAP point $\vec{\map}$ and to efficiently manipulate the covariance matrix $\postcov$.

\section{Design and software components of \hip}\label{sec:design}
\begin{figure}
\includegraphics[width=\textwidth]{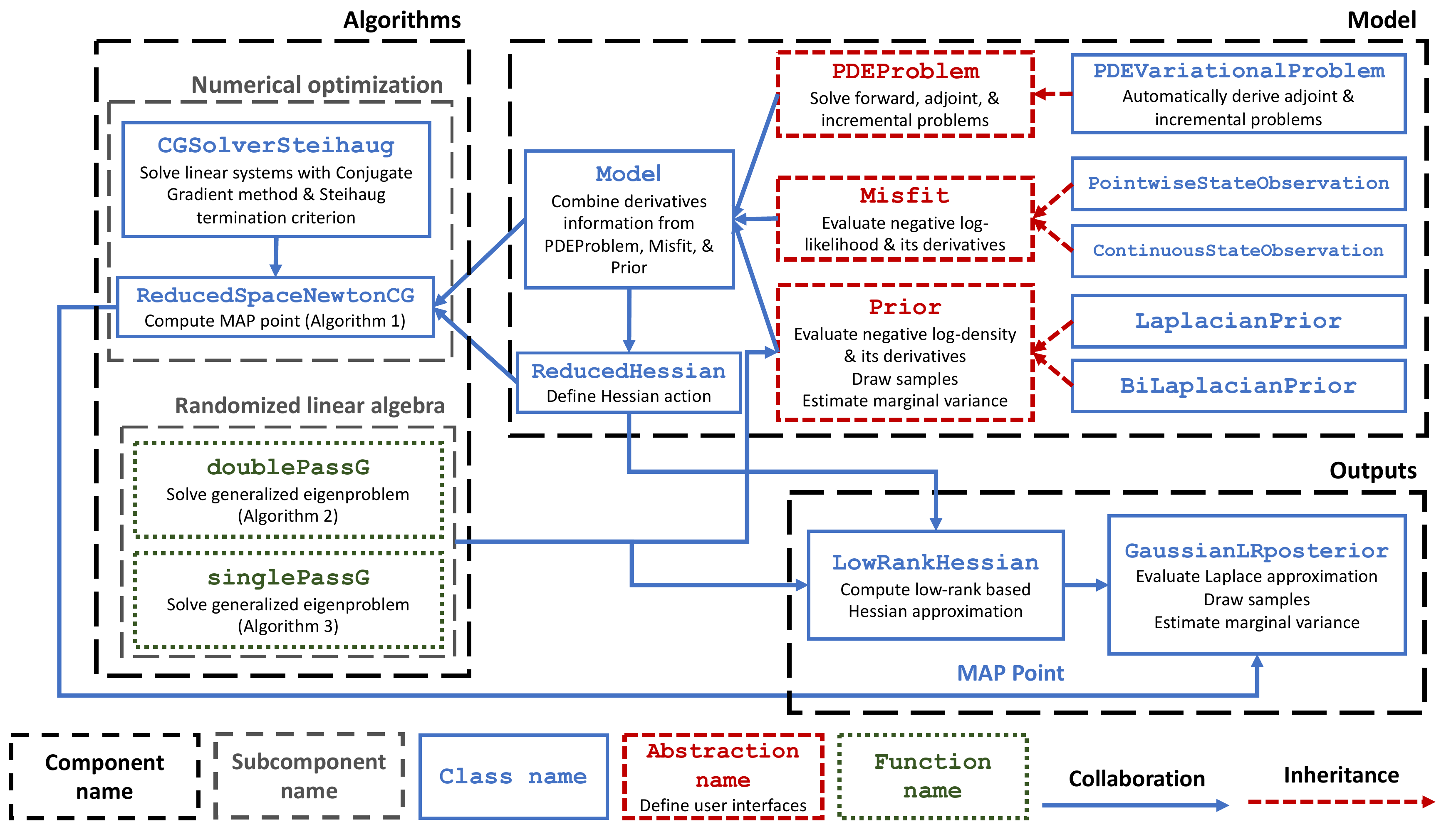}
\caption{Design of the \hip framework. Red dotted connectors represent inheritance and blue solid connectors represent collaborations among the main components (black dashed boxes), classes (blue solid boxes and red dashed boxes), and functionalities (green dotted boxes) of \hip. The legend is shown on the bottom.}
\label{fig:hippylib_design}
\end{figure}

\hip implements state-of-the-art scalable algorithms for PDE-based
deterministic and Bayesian inverse problems. It builds on \Fe (a
parallel finite element library)~\cite{LoggMardalWells12,LangtangenLogg17} for
discretization of PDEs and on \pet
\cite{BalayAbhyankarAdamsEtAl14} for scalable and efficient linear
algebra operations and solvers.  In \hip the user can express the
forward PDE and the likelihood in variational form 
using the friendly,
compact, near-mathematical notation of \Fe, which will then
automatically generate efficient code for the discretization.  Linear
and nonlinear, stationary and time-dependent PDEs are supported in
\hip.  For stationary problems, gradient and Hessian information can
be automatically generated by \hip using \Fe symbolic differentiation
of the relevant variational forms. For time-dependent problems, instead,
symbolic differentiation can only be used for the spatial terms, and the
contribution to gradients and Hessians arising from the time dynamics
 needs to be provided by the user.
Noise and prior covariance operators are modeled as inverses of
elliptic differential operators allowing us to build on fast
multigrid solvers for elliptic operators without explicitly
constructing the dense covariance operator.  
The main components, classes, and functionalities of \hip are summarized
in \figureref~\ref{fig:hippylib_design}.  These include:
\begin{enumerate}[leftmargin=*]
\item The {\bf \bhip model} component describes the inverse problem,
  i.e., the data misfit functional (negative log-likelihood), the
  prior information, and the forward problem.  More specifically, the
  user can select from among a library of data misfit
  functionals---such as pointwise observations or continuous
  observations in the domain or on the boundary---or implement new
  ones using the prescribed interface. \hip offers a library of priors
  the user can choose from and allows for user-provided priors as
  well. Finally, the user needs to provide the forward problem either
  in the form of a \Fe variational form or (for more complicated or
  time dependent problem) as a user-defined object.  When using \Fe
  variational forms, \hip is able to derive expressions for the
  gradient and Hessian-apply automatically using \Fe' symbolic
  differentiation. This means that for stationary problems the user
  will only have to provide the variational form of the forward
  problem. For more complex problems (e.g., time-dependent problems),
  \hip{} allows the user to implement their own derivatives.
\item The {\bf \bhip algorithms} component contains the numerical methods needed for
  solving the deterministic and linearized Bayesian inverse problems,
  i.e., the globalized inexact Newton-CG algorithm, randomized
  generalized eigensolvers, scalable sampling of Gaussian fields,
  and trace/diagonal estimators for large-scale
  not-explicitly-available covariance matrices. These algorithms
  are described in detail in Section~\ref{sec:algorithms}.
\item The {\bf \bhip outputs} component includes the parameter-to-observable map
  (and its linear approximation), gradient evaluation and Hessian
  action, and Laplace approximation of the posterior distribution (MAP
  point and low-rank based representation of the posterior covariance
  operator). The \hip outputs can be utilized as inputs to other UQ
  software, e.g., the MIT Uncertainty Quantification Library (MUQ), to
  perform a full characterization of the posterior distribution using
  advanced dimension-independent Markov chain Monte
  Carlo simulation, requiring derivative information.
\end{enumerate}

We refer to \cite{VillaPetraGhattas20a} for a detailed description of the modules, classes and functions implemented in \hip.

\section{\bhip algorithms}\label{sec:algorithms}
In this section we describe the main algorithms implemented in \hip
for solution of deterministic and linearized Bayesian inverse problems.
Specifically, we focus on computation of the MAP point and various operations on prior and posterior
covariance matrices.  In the linear case (and under the assumption of
Gaussian noise and Gaussian prior) the posterior distribution is also
 Gaussian, and therefore is fully characterized once the MAP point
and posterior covariance matrix are computed.  In the nonlinear case,
efficient exploration of the posterior distribution for large-scale PDE problems will require the use of a
Markov chain Monte Carlo (MCMC) sampling method enchanted by Hessian information
(see e.g. \cite{MartinWilcoxBursteddeEtAl12,CuiLawMarzouk14,
CuiMartinMarzoukEtAl14,Bui-ThanhGirolami14,BeskosGirolamiLanEtAl17}). In this case \hip
provides the tools to generate proposals for MCMC.

\subsection{Deterministic inversion and MAP point computation via inexact Newton-CG}
\label{subsec:NCG}

\hip provides a robust implementation of the inexact Newton-conjugate
gradient (Newton-CG) algorithm (e.g.,\cite{AkcelikBirosGhattasEtAl06a,
  BorziSchulz12}) to solve the deterministic inverse problem and, 
  in the Bayesian framework, to compute the maximum a posterior (MAP) point (see
Algorithm \ref{alg:InexactNewtonCG}).  The gradient and
Hessian actions---whose expressions are given in \eqref{eq:abstract_gradient}  and \eqref{eq:abstract_hessianAction}, respectively---are automatically computed via their variational form
specification in \Fe by constraining the state and adjoint variables
to satisfy the forward and adjoint problem in \eqref{eq:abstract_fwd} and 
\eqref{eq:abstract_adj} respectively. 
The Newton system is
solved inexactly using early termination of CG iterations using
Eisenstat--Walker~\cite{EisenstatWalker96} (to
prevent oversolving) and Steihaug \cite{Steihaug83} (to avoid negative
curvature) criteria. Specifically, the choice of the tolerance $\eta_i$ in Algorithm \ref{alg:InexactNewtonCG} leads to superlinear convergence of Newton's method,
and represents a good compromise between the number of Newton iterations and the computational effort to compute the search direction.
 Globalization is achieved with an Armijo
backtracking line search; we choose the Armijo constant $c_{\rm armijo}$ in the interval $[10^{-5}, 10^{-4}]$. For a wide class of nonlinear inverse
problems, the number of outer Newton iterations and inner CG
iterations is independent of the mesh size and hence parameter
dimension \cite{Heinkenschloss93}.  This is a consequence of using Newton's method, the compactness of the Hessian (of the data misfit
term), and preconditioning with the inverse regularization
operator. We note that the resulting preconditioned Hessian is a compact
perturbation of the identity, for which Krylov subspace methods
exhibit mesh-independent iterations \cite{CampbellIpsenKelleyEtAl96}.

\begin{algorithm}[t!]
\begin{algorithmic}
\STATE $i \gets 0$
\STATE Given $\vec{m}_0$ solve the forward problem \eqref{eq:abstract_fwd} to obtain $\vec{u}_0$
\STATE Given $\vec{m}_0$ and $\vec{u}_0$ compute the cost functional $\mathcal{J}_0$ using \eqref{eq:optpb}
\WHILE{$i < $ max\_iter}
\STATE Given $\vec{m}_i$ and $\vec{u}_i$ solve the adjoint problem \eqref{eq:abstract_adj} to obtain $\vec{p}_i$
\STATE Given $\vec{m}_i$, $\vec{u}_i$ and $\vec{p}_i$ evaluate the gradient $\vec{g}_i$ using  \eqref{eq:abstract_gradient}
\IF{$\|\vec{g}_i\| \leq \tau$ }
\STATE {\bf break}
\ENDIF
\STATE Given $\vec{m}_i$, $\vec{u}_i$ and $\vec{p}_i$ define a linear operator $\mat{H}_i$ that implements the Hessian action \eqref{eq:abstract_hessianAction}
\STATE Using conjugate gradients, find a search direction $\widehat{\vec{m}}_i$ such that
\vspace{-2mm}
\begin{equation*}
\| \mat{H}_i \widehat{\vec{m}}_i + \vec{g}_i \| \leq \eta_i \| \vec{g}_i \| , \text{ with }  \eta_i = \left( \frac{\|\vec{g}_i\|}{\|\vec{g}_0\|} \right)^{\frac{1}{2}}
\end{equation*}
\vspace{-2mm}
\STATE $j \gets 0$, $\alpha^{(0)} \gets 1$
\WHILE{$j < $ max\_backtracking\_iter}
\STATE Set $\vec{m}^{(j)} = \vec{m}_i + \alpha^{(j)} \widehat{\vec{m}}_i$
\STATE Given $\vec{m}^{(j)}$ solve the forward problem \eqref{eq:abstract_fwd}  to obtain $\vec{u}^{(j)}$
\STATE Given $\vec{m}^{(j)}$ and $\vec{u}^{(j)}$ compute the cost $\mathcal{J}^{(j)}$ using \eqref{eq:optpb}
\IF{ $\mathcal{J}^{(j)} < \mathcal{J}_i + \alpha^{(j)} c_{\rm armijo} \, \vec{g}_i^T \widehat{\vec{m}}_i$}
\STATE $\vec{m}_{i+1} \gets  \vec{m}^{(j)}$, $\mathcal{J}_{i+1} \gets \mathcal{J}^{(j)}$
\STATE {\bf break}
\ENDIF
\STATE $\alpha^{(j+1)} \gets \alpha^{(j)} / 2, \quad j \gets j+1 $
\ENDWHILE
\IF{$  j = $ max\_backtracking\_iter} 
\STATE {\bf break}
\ENDIF
\STATE $i \gets i+1$
\ENDWHILE
\end{algorithmic}
\caption{The inexact Newton-CG algorithm to find the MAP point}
\label{alg:InexactNewtonCG}
\end{algorithm}

\subsection{Low-rank approximation of the Hessian}\label{sec:HessianLowRank}
\label{subsec:lra}

The Hessian (of the negative log-posterior) plays a critical role in inverse problems. First, its spectral properties characterize the degree of ill-posedness. 
Second, the Hessian is the underlining operator for Newton-type optimization algorithms, which are highly desirable when solving inverse problems due to their dimension-independent convergence property.
Third, the inverse of the Hessian locally characterizes the uncertainty in the solution of the inverse problem; under the Laplace approximation, it is precisely the posterior covariance matrix. Unfortunately, after discretization, the Hessian is formally a large, dense matrix; forming each column requires an incremental forward and adjoint solves (see Section \ref{sec:det}). Thus, construction of the Hessian is prohibitive for large-scale problems since its dimension is equal to the dimension of the parameter. To make operations with the Hessian tractable, we exploit the fact that, in many cases, the eigenvalues collapse to zero rapidly, since the data contain limited information about the (infinite-dimensional) parameter field. Thus a low-rank approximation of the data misfit component of the Hessian, $\Hmisfit$, can be constructed.
This can be proven analytically for certain linear forward PDE problems (e.g. advection-diffusion \cite{FlathWilcoxAkcelikEtAl11}, Poisson \cite{Flath13}, Stokes \cite{Worthen12}, acoustics \cite{Bui-ThanhGhattas12, Bui-ThanhGhattas13}, electromagnetics \cite{Bui-ThanhGhattas13}), and demonstrated numerically for more complex PDE problems (e.g. seismic wave propagation \cite{Bui-ThanhGhattasMartinEtAl13,Bui-ThanhBursteddeGhattasEtAl12_gbfinalist}, mantle convection \cite{WorthenStadlerPetraEtAl14},
ice sheet flow \cite{IsaacPetraStadlerEtAl15,PetraMartinStadlerEtAl14}, poroelasticity \cite{HesseStadler14}, and turbulent flow \cite{ChenVillaGhattas19}).
The end result is that manipulations with the Hessian require a number of forward PDE solves that is independent of the parameter and data dimensions.

More specifically, to compute the low-rank factorization of the data misfit component of the Hessian we consider the following generalized symmetric eigenproblem:
\begin{equation}\label{eq:eigenproblem_1}
\Hmisfit \vec{v}_i = \lambda_i \mat{R} \vec{v}_i, \quad \lambda_1 \geq \lambda_2 \geq \ldots \geq \lambda_n,
\end{equation}
where $\mat{R}$ stems from the discretization (with respect to the Euclidean inner product) of the inverse of the prior covariance (i.e., the regularization operator).
We then choose $r \ll n$ such that $\lambda_{r+i}$, $0 < i \leq n-r$, is \emph{small} relative to 1, and we define
\begin{equation*}
\mat{V}_r = \begin{bmatrix} \vec{v}_1, \vec{v}_2,\ldots \vec{v}_r \end{bmatrix} \text{ and } \mat{\Lambda}_r = {\rm diag}([\lambda_1, \lambda_2, \ldots, \lambda_r]),
\end{equation*}
where the matrix $\mat{V}_r$ has $\mat{R}$-orthonormal columns, that is
$\mat{V}_r^T \mat{R} \mat{V}_r = \mat{I}_r$. As in \cite{IsaacPetraStadlerEtAl15}, by using the Sherman-Morrison-Woodbury formula, we write
\begin{align}
 \mat{H}^{-1} = \left(\mat{R} +  \Hmisfit \right)^{-1} = \mat{R}^{-1} - \mat{V}_r \mat{D}_r \mat{V}_r^T + \mathcal{O}\left(\sum_{i=r+1}^n \frac{ \lambda_i }{1 + \lambda_i}\right), \label{eq:sm}
\end{align}
where $\mat{D}_r =  \diag(\lambda_1/(\lambda_1+1), \dots, \lambda_r/(\lambda_r+1)) \in \R^{r\times r}$.
As can be seen from the form of the remainder term above, to obtain an
accurate low-rank approximation of $\mat{H}^{-1}$, we can neglect
eigenvectors corresponding to eigenvalues that are small compared
to~1. This result is used to efficiently apply the inverse and
square-root inverse of the Hessian to a vector, as needed for
computing the pointwise variance and when drawing samples from a
Gaussian distribution with covariance $\mat{H}^{-1}$, as will be shown
in Sections~\ref{sec:samplingpost} and~\ref{sec:varpost},
respectively. Efficient algorithms implemented in \hip for solving eigenproblems using randomized linear algebra methods are described next.

\subsubsection{Randomized algorithm for the generalized eigenvalue problem}

Randomized algorithms for eigenvalue computations have proven to be extremely effective for matrices with rapidly decaying eigenvalues \cite{HalkoMartinssonTropp11}.
For this class of matrices, in fact, randomized algorithms present several advantages compared to Krylov subspace methods.
Krylov subspace methods require sophisticated algorithms to monitor restart, orthogonality, and loss of precision.
On the contrary, randomized algorithm are easy to implement, can be made numerically robust, and expose more opportunities for parallelism since matrix-vector products can be done asynchronously across all vectors. The flexibility in reordering the computation makes randomized algorithms particularly well suited for modern parallel architectures with many cores per node and deep memory hierarchies.

In \hip we apply randomized algorithms to compute the low-rank factorization of the misfit part of the Hessian $\Hmisfit$.
With a change of notation, we write the generalized eigenvalue problem \eqref{eq:eigenproblem_1} as
\begin{equation}\label{eq:abstract_eigenp}
\mat{A} \vec{v} = \lambda \mat{B} \vec{v}
\end{equation}
where $\mat{A} \in \mathbb{R}^{n\times n}$ is symmetric, $\mat{B} \in \mathbb{R}^{n\times n}$ is symmetric positive definite, and $\vec{v} \in \mathbb{R}^{n}$.
Here we present an extension of the randomized eigensolvers in \cite{HalkoMartinssonTropp11}
to the solution of the generalized symmetric eigenproblem \eqref{eq:abstract_eigenp}. Randomized algorithms for generalized symmetric eigenproblems were first introduced in  \cite{SaibabaLeeKitanidis16},  and are revisited here with some modifications.

The main idea behind randomized algorithms is to construct a  matrix $\mat{Q} \in \mathbb{R}^{n \times (r+l)}$ with  $\mat{B}$-orthonormal columns that approximates the range of $\mat{B}^{-1}\mat{A}$. Here, $r$ represents the number of eigenpairs we wish to compute, and $l$ is an \emph{oversampling factor}.
More specifically, we have
\begin{equation}\label{eq:approx_randomized}
\| (\mat{I} - \mat{Q}\mat{Q}^T) \mat{A} \|_{\mat{B}} \leq \epsilon,
\end{equation}
where $\epsilon$ is a random variable whose distribution depends on the generalized eigenvalues of \eqref{eq:abstract_eigenp} with index greater than $r+l$. To construct $\mat{Q}$, we let $\mat{\Omega} \in \mathbb{R}^{n \times (r+l)}$ be a Gaussian random matrix---whose entries are independent identically distributed (\emph{i.i.d.}) standard Gaussian random variables---and we compute a $\mat{B}$-orthogonal basis for the range of $\mat{Y} = \mat{B}^{-1}\mat{A} \mat{\Omega}$ using the so called PreCholQR algorithm (\cite{LoweryLangou14}, see Algorithm \ref{alg:preCholQR}). The main computational cost is the construction of $\mat{Y}$, which requires $(r+l)$ applications of the operator $\mat{A}$, and $(r+l)$ linear solves to apply $\mat{B}^{-1}$. In contrast, the computation of the matrix $\mat{Q}$ with $\mat{B}$-orthonormal columns using PreCholQR requires only an additional $(r+l)$ applications of $\mat{B}$ and $\mathcal{O}(n(r+l)^2)$ dense linear algebra operations for the QR factorization. 
Using \eqref{eq:approx_randomized} it can be shown (see \cite{HalkoMartinssonTropp11}) that
$ \mat{A} \approx (\mat{B}\mat{Q})( \mat{Q}^T \mat{A} \mat{Q})(\mat{B}\mat{Q})^T =  (\mat{B}\mat{Q}) \mat{T} (\mat{B}\mat{Q})^T,$
where we have defined $\mat{T} := \mat{Q}^T \mat{A} \mat{Q}$. Then we compute the eigendecomposion $\mat{T} = \mat{S} \mat{\Lambda} \mat{S}^T$ ($\mat{S}^T \mat{S} = \mat{I}_{r+l}$), and we approximate the $r$ dominant eigenpairs $(\mat{\Lambda}_r, \mat{V}_r)$ of \eqref{eq:abstract_eigenp} by:
\begin{align}
\mat{\Lambda}_r = \mat{\Lambda}(1:r, 1:r), \quad \mat{V}_r = \mat{Q}\mat{S}(:,1:r).
\end{align} 

\begin{figure}
\includegraphics[width=0.32\textwidth]{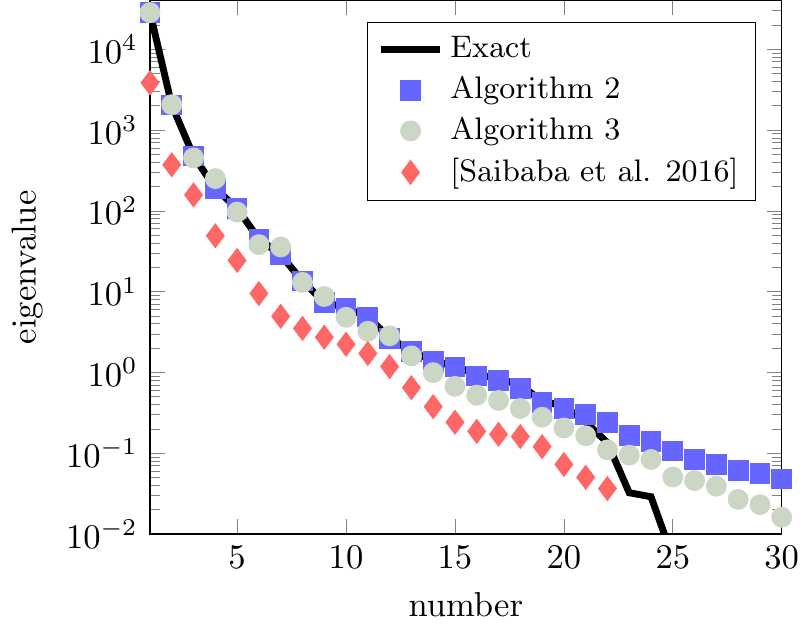}\hfill
\includegraphics[width=0.32\textwidth]{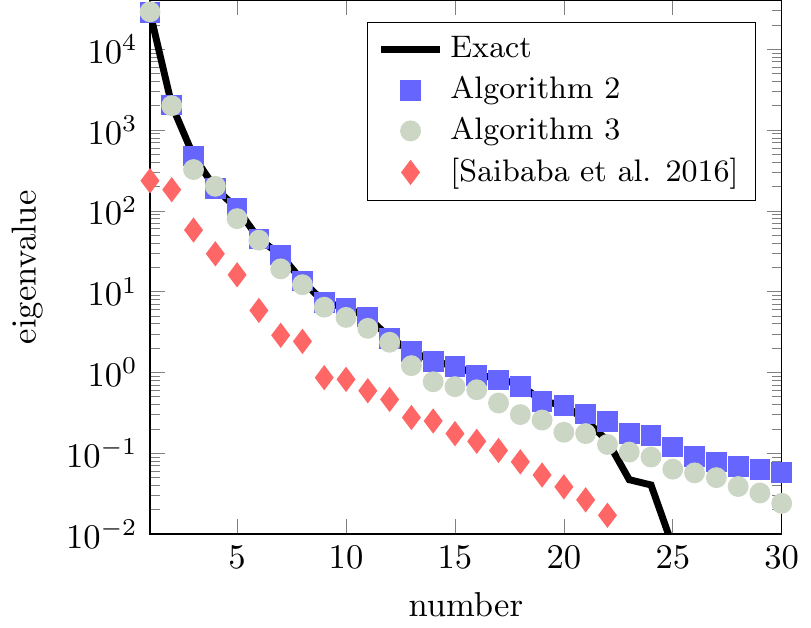}\hfill
\includegraphics[width=0.32\textwidth]{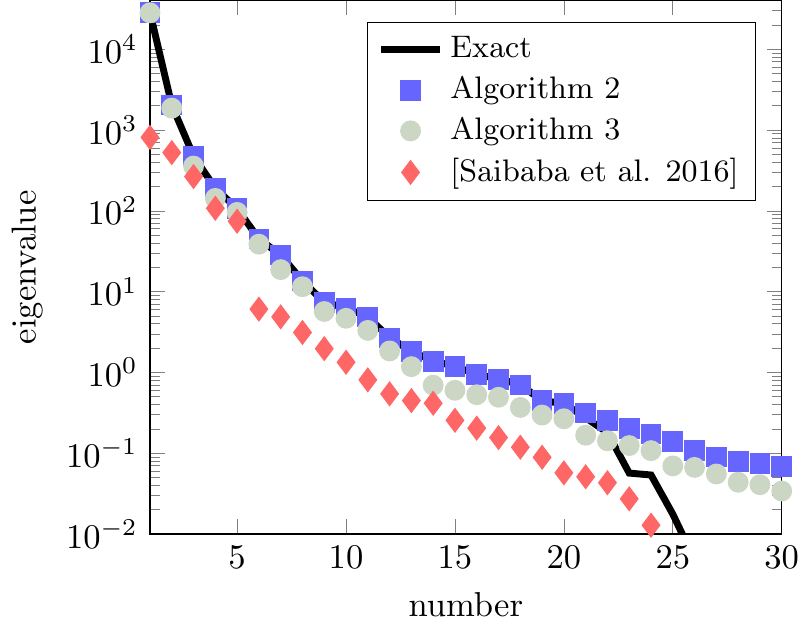}
\caption{ Log-linear plot of first 30 out of 4225 generalized eigenvalues of the data misfit Hessian in Section \ref{subsec:poisson}. A deterministic eigensolver (Exact), Algorithm \ref{alg:doublepass}, Algorithm \ref{alg:singlepass}, and the single pass algorithm in \cite{SaibabaLeeKitanidis16} are compared for different choices of the oversampling parameter $l$ (left $l=5$, center $l=10$, right $l=20$). Note that Algorithm \ref{alg:singlepass} is more accurate for the single pass algorithm in \cite{SaibabaLeeKitanidis16} for all choices of $l$.}
\label{fig:SinglePassComparison}
\end{figure}

Algorithms \ref{alg:doublepass} and \ref{alg:singlepass} summarize the implementation of the double pass and single pass randomized algorithms~\cite{HalkoMartinssonTropp11}.
The main difference between these two algorithms is how the \emph{small} matrix $\mat{T}$ is computed. In the double pass algorithm, $\mat{T}$ is computed directly by performing a second round of multiplication $\mat{A}\mat{Q}$ with the operator $\mat{A}$. In the single pass algorithm, $\mat{T}$ is approximated from the information contained in $\mat{\Omega}$ and $\mat{Y}$. In particular, generalizing the single pass algorithm in \cite{HalkoMartinssonTropp11} to \eqref{eq:abstract_eigenp}, we approximate $\mat{T}$ as the least-squares solution of
\begin{equation}\label{eq:Tsinglepass}
\mat{T} = \argmin_{\mat{X} \in \mathbb{R}^{(r+l)\times(r+l)}, s.s.p.d} \| \mat{X} (\mat{Q}^T \mat{B} \mat{\Omega}) - \mat{Q}^T \mat{B} \mat{Y} \|^2_2.
\end{equation}
For this reason, the single pass algorithm has a lower computational cost compared to the double pass algorithm; however the resulting approximation is less accurate.
We remark that Algorithm \ref{alg:singlepass} is more accurate than the single pass algorithm presented in \cite{SaibabaLeeKitanidis16}. The key difference between the two algorithms is in the definition of $\mat{T}$.  In \cite{SaibabaLeeKitanidis16}, the authors define $\mat{T} = (\mat{Q}^T \mat{B} \mat{\Omega})^{-1} (\mat{Q}^T \mat{B} \mat{Y}) (\mat{Q}^T \mat{B} \mat{\Omega})^{-1}$, while in Algorithm \ref{alg:singlepass} we define $\mat{T}$ as the least-squares solution of \eqref{eq:Tsinglepass}. \figureref~\ref{fig:SinglePassComparison} numerically illustrates the higher accuracy of the proposed approach when computing the first 30 eigenvalues of the data misfit Hessian discussed in Section \ref{subsec:poisson}.

\begin{algorithm}
\begin{algorithmic}
\STATE Let $r$ be the number of eigenpairs to compute and $l$ an oversampling factor
\STATE Let $\mat{\Omega} \in \mathbb{R}^{n \times (r+l)} $ be a Gaussian random matrix
\STATE $\bar{\mat{Y}} \gets \mat{A} \mat{\Omega}$, $\mat{Y} = \mat{B}^{-1} \bar{\mat{Y}}$
\STATE Use \texttt{PreCholQR} to factorize $\mat{Y} = \mat{Q}\mat{R}$ such that $\mat{Q}^T \mat{B} \mat{Q} = \mat{I}_{r+l}$
\STATE $\mat{T} \gets \mat{Q}^T \mat{A} \mat{Q}$
\STATE Compute the eigenvalue decomposition $\mat{T} = \mat{S} \mat{\Lambda} \mat{S}^T$
\STATE Keep the $r$ largest eigenmodes and let $\mat{S}_r \gets \mat{S}(:,1:r)$, $\mat{\Lambda}_r \gets \mat{\Lambda}(1:r,1:r)$
\STATE {\bf Return:} $\mat{V}_{r} \gets \mat{Q} \mat{S}_r$, and $\mat{\Lambda}_r$
\end{algorithmic}
\caption{The double pass randomized algorithm for the solution of the generalized symmetric eigenproblem.}
\label{alg:doublepass}
\end{algorithm}

\begin{algorithm}
\begin{algorithmic}
\STATE Let $r$ be the number of eigenpairs to compute and $l$ an oversampling factor
\STATE Let $\mat{\Omega} \in \mathbb{R}^{n \times (r+l)} $ be a Gaussian random matrix
\STATE $\bar{\mat{Y}} \gets \mat{A} \mat{\Omega}$, $\mat{Y} = \mat{B}^{-1} \bar{\mat{Y}}$
\STATE Use PreCholQR to factorize $\mat{Y} = \mat{Q}\mat{R}$ such that $\mat{Q}^T \mat{B} \mat{Q} = \mat{I}_{r+l}$ and $\bar{\mat{Q}}$ such that $\bar{\mat{Q}}^T \mat{B}^{-1} \bar{\mat{Q}} = \mat{I}_{r+l}$
\STATE Find $\mat{T}$ s.s.p.d such that $\|\mat{T} (\bar{\mat{Q}}^T\mat{\Omega}) - \bar{\mat{Q}}^T \mat{Y}\|^2_2 \rightarrow \min$
\STATE Compute the eigenvalue decomposition $\mat{T} = \mat{S} \mat{\Lambda} \mat{S}^T$
\STATE Keep the $r$ largest eigenmodes and let $\mat{S}_r \gets \mat{S}(:,1:r)$, $\mat{\Lambda}_r \gets \mat{\Lambda}(1:r,1:r)$
\STATE {\bf Return:} $\mat{V}_{r} \gets \mat{Q} \mat{S}_r$, and $\mat{\Lambda}_r$
\end{algorithmic}
\caption{The single pass randomized algorithm for the solution of the generalized symmetric eigenproblem.}
\label{alg:singlepass}
\end{algorithm}

\begin{algorithm}
\begin{algorithmic}
\STATE {\bf Require:} $\mat{Y} \in \mathbb{R}^{n\times (r+l) }$, and $\mat{B} \in \mathbb{R}^{n \times n}$
\STATE $[\mat{Z}, \mat{R}_{\mat{Y}}] \gets {\rm qr}(\mat{Y})$
\STATE $\bar{\mat{Z}} \gets \mat{B} \mat{Z}$
\STATE $\mat{R}_{\mat{Z}} = {\rm chol}(\mat{Z}^T \bar{\mat{Z}})$
\STATE {\bf Return:} $\mat{Q} = \mat{Z} \mat{R}_\mat{Z}^{-1}$, $\bar{\mat{Q}} = \bar{\mat{Z}} \mat{R}_{\mat{Z}}^{-1}$, and $\mat{R} = \mat{R}_{\mat{Z}} \mat{R}_{\mat{Y}}$
\end{algorithmic}
\caption{PreCholQR}
\label{alg:preCholQR}
\end{algorithm}

\subsection{Sampling from large-scale Gaussian random fields}
\label{subsec:sampling}
Sampling techniques play a fundamental role in exploring the posterior distribution and in quantifying the uncertainty in the inferred parameter; for example, Markov chain Monte Carlo methods often use the prior distribution (assumed Gaussian in our settings) or some Gaussian approximation to the posterior to generate proposals for the Metropolis-Hastings algorithm.
In this section, we describe the sampling capabilities implemented in \hip to generate realizations of Gaussian random fields with a prescribed covariance operator $\mathcal{C}$.  Then we describe how the low-rank factorization of the 
data misfit part of the Hessian in Section \ref{sec:HessianLowRank} can be exploited to generate samples from the Laplace approximation of the posterior distribution in \eqref{eq:laplace_approx}. In what follows, we will denote the expected value (mean) of a random vector $\vec{x}$ with the symbol $\ave\left[ \vec{x} \right]$, and its covariance with the symbol ${\rm cov}( \vec{x} ) := \ave\big[ \left(\vec{x} - \ave[\vec{x}]\right)\left(\vec{x} - \ave[\vec{x}]\right)^T\big]$. We will also denote with  $\mat{\Gamma}$ the matrix representation of the covariance operator  $\mathcal{C}$ with respect to the standard Euclidean inner product\footnote{Note that $\mat{\Gamma}$ differs from $\prcov$ and $\postcov$, which are defined in terms of the $\M$-weighted inner product.}.

To sample from a small-scale multivariate Gaussian distribution, it is common to resort to a Cholesky factorization of the covariance matrix $\mat{\Gamma} = \mat{C} \mat{C}^T$. In fact, if $\vec{\eta}$ is a vector of independent identically distributed Gaussian variables $\eta_i$ with zero mean ($\ave\left[ \vec{\eta} \right] = \vec{0}$ ) and unit variance (${\rm cov}(\vec{\eta}) = \mat{I}$), then $\vec{x} = \mat{C} \vec{\eta}$ is such that
$$ {\rm cov}( \vec{x} ) = \ave\left[ \vec{x}\vec{x}^T\right] = \ave\left[ \mat{C} \vec{\eta} \vec{\eta}^T \mat{C}^T \right] = \mat{C} \ave\left[ \vec{\eta} \vec{\eta}^T \right] \mat{C}^T = \mat{C} \mat{C}^T = \mat{\Gamma}.$$
Since an affine transformation of a Gaussian vector is still Gaussian, we have that $\vec{x} \sim \mathcal{N}(\vec{0}, \mat{\Gamma} )$.

This approach is not feasible for large-scale problems since it requires computing the Cholesky factorization of the covariance matrix.
However, note that a decomposition of the form $\mat{\Gamma} = \mat{C} \mat{C}^T$ can be obtained using a matrix $C$ other than the Cholesky factor. In particular, the matrix $\mat{C}$  need not be a triangular, or even square, matrix. 
In Appendix~\ref{sec:rect_decomp}, we exploit this observation and show a scalable sampling technique based on a rectangular decomposition of $\mat{R} = \mat{\Gamma}^{-1}$, for the case when $\mat{R}$ is a finite element discretization of a differential operator. We note that a similar approach was, independently, investigated in \citep{CrociGilesRognesEtAl18} to sample realizations of  white noise by exploiting a rectangular decomposition of the finite element mass matrix. Our approach is more general as it allows for decomposing matrices stemming from finite element discretization of operators other than identity.

\subsubsection{Sampling from the prior}

Sampling from a Gaussian distribution with a prescribed covariance matrix $\mat{\Gamma}$ is a difficult task for large-scale problems.
Different approaches have been investigated, but how to make these algorithms scalable is still an active area of research. In \cite{ParkerFox12}, the authors introduce a conjugate gradient sampler that is a simple extension of the conjugate gradient method for solving linear systems. However, loss of orthogonality in finite arithmetic and the need of a factorized preconditioner limit the efficiency of this sampler for large-scale applications. In \cite{ChowSaad14}, the authors consider a preconditioned Krylov subspace method to approximate the action of $\mat{\Gamma}^{\frac{1}{2}}$ on a generic vector $z$. In \cite{ChenAnitescuSaad11}, the authors discuss a method to compute $f(\mat{\Gamma})\vec{b}$ via least-squares polynomial approximations for a generic matrix function $f(x) = \sqrt{x}$. To this aim the authors approximate the function by a spline of a desired accuracy on the spectrum of $\mat{\Gamma}$ and introduce a weighted inner product to simplify the computation.

\hip implements a new sampling algorithm
that strongly relies on the structure of the covariance matrix and on
the assembly procedure of finite element matrices. In particular, we
restrict ourselves to the class of priors described in Section
\ref{sec:bayes}.  For this class of priors, the inverse of the
covariance matrix admits a sparse representation as a finite element
matrix $\mat{R}$ stemming from the finite element discretization of a
coercive symmetric differential operator. To draw a sample $\vec{x}$ from the
distribution $\mathcal{N}\left(\vec{0}, \mat{R}^{-1}\right)$, we
solve the linear system
\begin{equation}\label{eq:sample_prior}
\mat{R} \vec{x} = \mat{C} \vec{\eta}, \text{ where } \vec{\eta} \sim \mathcal{N}\left(\vec{0}, \mat{I}_q\right),
\end{equation}
where $\mat{C} \in \mathbb{R}^{n \times q}$ ($q \geq n$) is the rectangular factor of $\mat{R}$ described in Appendix \ref{sec:rect_decomp}.
In particular, we have
$$ \ave[\vec{x}\vec{x}^T] = \mat{R}^{-1} \mat{C} \ave[ \vec{\eta}\vec{\eta}^T] \mat{C}^T \mat{R}^{-1} = \mat{R}^{-1}\, \mat{R} \, \mat{R}^{-1} = \mat{R}^{-1},$$
where we exploited the fact that $\ave[ \vec{\eta}\vec{\eta}^T] = \mat{I}_q$ and $\mat{C}\mat{C}^T = \mat{R}$.
This method is particularly efficient at large-scale since: \emph{i}) the
matrix $\mat{C}$ is sparse and can be computed efficiently by
exploiting the finite element assembly routine; \emph{ii}) the dominant cost is
the solution of a linear system with coefficient matrix $\mat{R}$, for which efficient and
scalable methods are available (e.g. conjugate gradients with algebraic
multigrid preconditioner); and \emph{iii}) the stochastic dimension of
$\vec{\eta}$ also scales linearly with the size of the problem.

\subsubsection{Sampling from the Laplace approximation of the posterior}
\label{sec:samplingpost}

To sample from the Laplace approximation of the posterior we assume that the posterior covariance operator can be expressed as a low-rank update of the prior covariance, i.e., in the form of equation \eqref{eq:sm}. This assumption is often verified for many inverse problems as we discussed in Section \ref{sec:HessianLowRank}.
Then given a sample from the prior distribution $\vec{x} \sim \mathcal{N}\left( \vec{0}, \mat{R}^{-1}\right)$, a sample from the Laplace approximation of the posterior $\mathcal{N}\left( \mat{0}, (\Hmisfit + \mat{R})^{-1}\right)$ can be computed as

\begin{equation}\label{eq:sample_post}
\vec{y} = \left( \mat{I}_n - \mat{V}_r \mat{S}_r \mat{V}_r^T \mat{R} \right)\vec{x}, 
\end{equation}
where
  $\mat{S}_r = \mat{I}_r - (\mat{\Lambda}_r + \mat{I}_r)^{-\frac{1}{2}} =
  \diag(1-1/\sqrt{\lambda_1+1}, \dots, 1-1/\sqrt{\lambda_r+1}) \in \R^{r\times r}$.
This can be verified by the following calculation,
\begin{alignat*}{2}
  {\rm cov}( \vec{y} ) = \ave \left[ \vec{y} \vec{y}^T \right]
  &= \ave \left[ \left( \mat{I}_n - \mat{V}_r \mat{S}_r \mat{V}_r^T \mat{R} \right) \vec{x}\vec{x}^T \left( \mat{I}_n - \mat{V}_r \mat{S}_r \mat{V}_r^T \mat{R} \right)^T \right]\\
  &= \left( \mat{I}_n - \mat{V}_r \mat{S}_r \mat{V}_r^T \mat{R} \right) \ave\left[ \vec{x}\vec{x}^T \right] \left( \mat{I}_n - \mat{V}_r \mat{S}_r \mat{V}_r^T \mat{R} \right)^T\\
  &=\left( \mat{I}_n - \mat{V}_r \mat{S}_r \mat{V}_r^T \mat{R} \right) \mat{R}^{-1} \left( \mat{I}_n - \mat{V}_r \mat{S}_r \mat{V}_r^T \mat{R} \right)^T \\
  &= \mat{R}^{-1} - \mat{V}_r (2\mat{S}_r - \mat{S}_r^2) \mat{V}_r^T \\
  &= \mat{R}^{-1} - \mat{V}_r \mat{D}_r \mat{V}_r^T \approx \mat{H}^{-1},
\end{alignat*}
where we have used the definition of $\mat{D}_r$, the $\mat{R}$-orthogonality of the eigenvectors
matrix $\mat{V}_r$ (i.e., $\mat{V}_r^T \mat{R} \mat{V}_r = \mat{I}_r$), and the fact that
\begin{align*}
  2 \left(1 - \frac{1}{ \sqrt{1 + \lambda_i} }\right)
  \! - \! \left(1 - \frac{1}{ \sqrt{1 + \lambda_i} }\right)^2
  \!=\!1 - \left[ \left(1 - \frac{1}{ \sqrt{1 + \lambda_i} } \right)
    \! - \! 1 \right]^2 = \frac{\lambda_i}{1+\lambda_i}.
\end{align*}

\subsection{Pointwise variance of Gaussian random fields}
\label{subsec:var}
Consider a Gaussian random field $m \sim \mathcal{N}\left(0, \mathcal{C} \right)$ and its discrete counterpart $\vec{m} \sim \mathcal{N}\left(\vec{0}, \mat{Q}^{-1} \right)$, where $\mat{Q}$ is the precision matrix. Here $\mathcal{C} = \mathcal{C}_{\rm pr}$ and $\mat{Q} = \mat{R}$ if we are interested in the prior distribution; $\mathcal{C} = \mathcal{C}_{\rm post}$ and $\mat{Q} = \mat{H}$ for the posterior distribution. Then we define the pointwise variance of $m$ as the field $\sigma^2(x)$ such that
$$ \sigma^2(x) = {\rm Var}[ \ipar(x) ], \quad \forall x \in \D.$$
In this section, we present an efficient numerical method to compute a finite element approximation $\sigma_h^2$ of $\sigma^2(x)$. As shown in \cite{Bui-ThanhGhattasMartinEtAl13},  the diagonal of $Q^{-1}$ (i.e., ${\rm diag}(\mat{Q}^{-1})$) is the vector corresponding to the coefficients of the expansion of $\sigma_h^2$ in the finite element basis.
A na\"ive approach would require solution of $n$ linear systems with $\mat{Q}$. This is not feasible for large-scale problems: even for the case when an optimal (i.e.,  $\mathcal{O}(n)$) solver for $\mat{Q}$ is available, the complexity is at least $\mathcal{O}(n^2)$ operations. In what follows we discuss stochastic estimators and probing methods to efficiently estimate the pointwise variance of the prior distribution and we explore how the low-rank representation of the data misfit component of the Hessian can be efficiently exploited to compute the pointwise variance of the Laplace approximation of posterior distribution.

\subsubsection{Pointwise variance of the prior}

Estimating the pointwise variance of the prior reduces to the well studied problem of estimating the diagonal of the inverse of a matrix $\mat{R}$. Recall that in our case,  $\mat{R}$ arises from finite element discretization of an elliptic differential operator. Two commonly used methods to solve this task are the stochastic estimator in \cite{BekasKokiopoulouSaad07} and the probing method in \cite{TangSaad12}. 
Specifically, the  unbiased stochastic estimator for the diagonal of the inverse of $\mat{R}$ in \cite{BekasKokiopoulouSaad07} reads
\begin{equation}\label{eq:stoc_diag_estimator}
{\rm diag}(\mat{R}^{-1}) \approx \left[ \sum_{j=1}^s \vec{z}_j \odot \vec{w}_j \right] \oslash \left[ \sum_{j=1}^s \vec{z}_j \odot \vec{z}_j \right],
\end{equation}
where $\vec{w}_j$ solves $\mat{R} \vec{w}_j = \vec{z}_j$ and $\vec{z}_j$ are random \emph{i.i.d.} vectors.
Here $\odot$ and $\oslash$ represent the element-wise multiplication and division operators of vectors, respectively. The convergence of the method is independent of the size of the problem, but convergence is in general slow.
The probing method in \cite{TangSaad12}, on the other hand, leads to faster convergence in the common situation in which $\mat{R}^{-1}$ exhibits a decay property, i.e., the entries far away from the diagonal are small. Probing vectors are determined by applying some coloring algorithm to the graph $\mathcal{G}$ whose adjacency matrix is the sparsity pattern of some power $k$ of $\mat{R}$. More specifically, there are as many probing vectors as the number of colors in the graph $\mathcal{G}$ and the probing vector associated to color $i$ is the binary vector whose non-zero entries correspond to the nodes of $\mathcal{G}$ colored with color $i$. The higher the power $k$, the more accurate will be the estimation, but also the more expensive due to the increased number of colors (and, therefore, of probing vectors to solve for).
The main shortcoming of this approach is that it is not mesh independent, i.e., as we refine the mesh we need to increase the power $k$ and therefore the number of probing vectors.

To overcome these difficulties, \hip implements, in addition to the methods mentioned above, a novel approach based on a randomized eigendecomposition of $\mat{R}^{-1}$, taking advantage of the fact that $\mat{R}$ is the discretization of an elliptic differential operator. Specifically, we write
\begin{equation}\label{eq:prior_mvar}
{\rm diag}(\mat{R}^{-1}) \approx \left[ \sum_{i=1}^{\overline{r}} \mu_i \overline{{\vec{v}}}_i \odot \overline{{\vec{v}}}_i \right],
\end{equation}
where $\{(\mu_i, \overline{\vec{v}}_i)\}_1^{\overline{r}}$ denote the approximation of the $\overline{r}$ dominant eigenpairs of the matrix $\mat{R}^{-1}$ obtained by using the double pass randomized eigensolver (Algorithm \ref{alg:doublepass}) with $l=0$, $\mat{A} = \mat{R}^{-1}$, and $\mat{B}$ the identity matrix.
The main advantage of this approach is that, thanks to the rapid decay of the eigenvalues $\mu_i$ ($ \mat{R}^{-1}$ is compact), $\overline{r}$ is much smaller than the number of samples $s$ necessary for the stochastic estimator in \eqref{eq:stoc_diag_estimator} to achieve a given accuracy and that, in contrast to the probing algorithm, it is independent of the mesh size.

To illustrate  the convergence properties of the proposed method, we estimate the marginal variance of the Gaussian prior distribution in \eqref{eq:PoissonPrior} for the model problem in Section \ref{subsec:poisson}.
\figureref~\ref{fig:comparison_diag_estimators}-a shows the superior accuracy of our method when compared to the stochastic estimator in \cite{BekasKokiopoulouSaad07} for a given number of covariance operator applies. \figureref~\ref{fig:comparison_diag_estimators}-b numerically demonstrates the mesh independence of the proposed method.

\def \pos {0.43\columnwidth}
\begin{figure}
  \begin{tikzpicture}
\node (1) at (0*\pos-0.7*\pos, 0*\pos){\includegraphics[width=.4\textwidth]{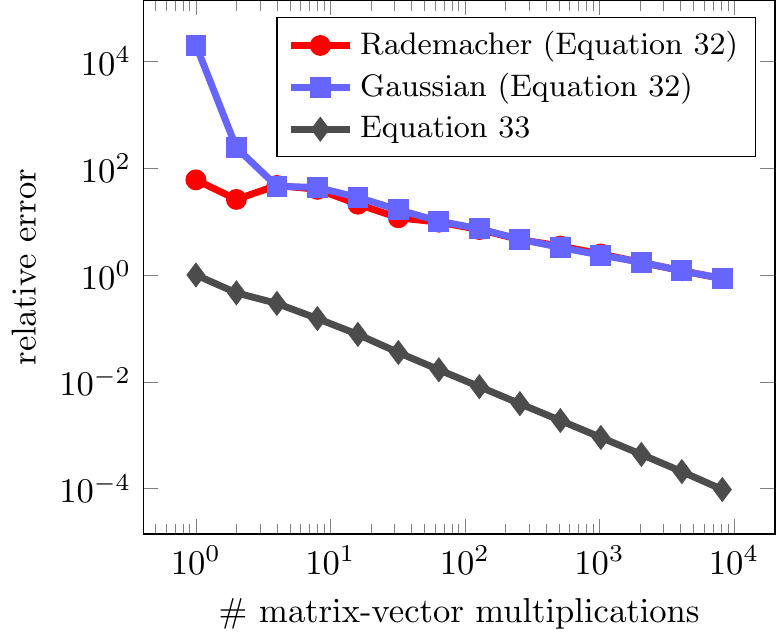}};
 \node (2) at (0.3*\pos, 0*\pos){\includegraphics[width=.4\textwidth]{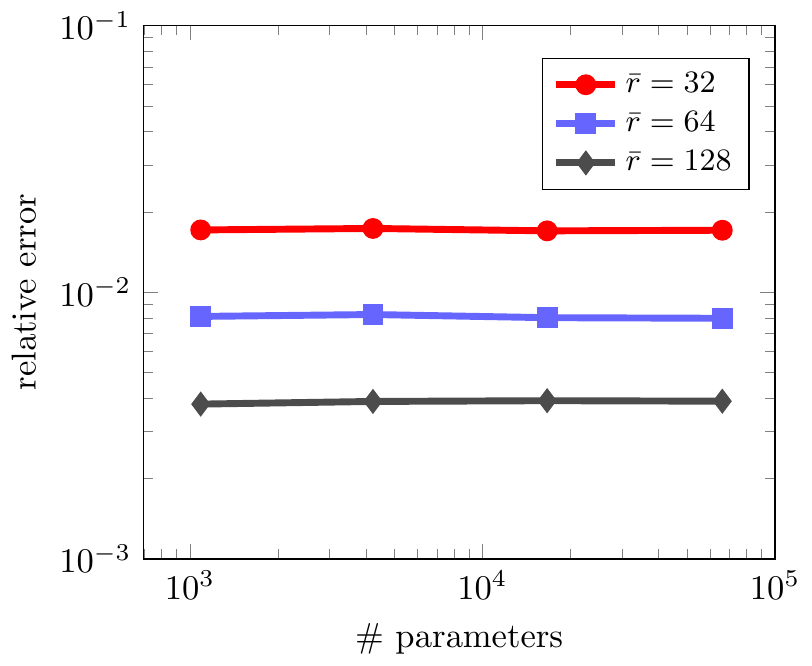}};
\node at (0.0*\pos-0.95*\pos, -0.205*\pos) {\sf a)};
\node at (1.0*\pos-0.95*\pos, -0.21*\pos) {\sf b)};
\end{tikzpicture}
\caption{The $L^2(\Omega)$ relative error for the estimation of the
  marginal variance for the prior distribution in
  Section~\ref{subsec:poisson} as a function of number of covariance
  operator applies (a) and number of parameters (b).}
\label{fig:comparison_diag_estimators}
\end{figure}

\subsubsection{Pointwise variance of the posterior}
\label{sec:varpost}
We resort to the low-rank representation of the data misfit component of the Hessian and the Woodbury formula to obtain the approximation
\begin{equation}\label{eq:post_mvar}
 {\rm diag}( \mat{H}^{-1} ) \approx {\rm diag}( \mat{R}^{-1} - \mat{V}_r \mat{D}_r \mat{V}_r^T) = {\rm diag}( \mat{R}^{-1}) - {\rm diag}(\mat{V}_r \mat{D}_r \mat{V}_r^T).
\end{equation}
In \hip, the first term is approximated using \eqref{eq:prior_mvar}, while the data-informed correction ${\rm diag}(\mat{V}_r \mat{D}_r \mat{V}_r^T)$ can be explicitly computed in $\mathcal{O}(n)$ operations as follows
$$ {\rm diag}(\mat{V}_r \mat{D}_r \mat{V}_r^T) = \sum_{i=1}^r \left[ \left( \frac{\lambda_i}{1+\lambda_i} \vec{v}_i\right) \odot \vec{v}_i \right]. $$

\section{Model problems}
\label{sec:models}

In this section we apply the inversion methods discussed in previous
sections to two model problems: inversion for the log coefficient
field in an elliptic partial differential equation and inversion for
the initial condition in a time-dependent advection-diffusion
equation. The main goal of this section is to illustrate the
deterministic inversion and linearized Bayesian analysis capabilities
of \hip for the solution of these two representative types of inverse
problems.  The numerical results showed below were obtained 
using \hip version 2.3.0 and \Fe 2017.2. A Docker image \cite{Merkel14} containing
the preinstalled software and examples can be downloaded at
\url{https://hub.docker.com/r/hippylib/toms}. 
For a line-by-line explanation of the source code for these
two model problems, we refer the reader to the Python Jupyter
notebooks available at \url{https://hippylib.github.io/tutorial_v2.3.0/}.

\subsection{Coefficient field inversion in a Poisson PDE problem}
\label{subsec:poisson}

In this section, we study the inference of the log coefficient field
$\ipar$ in a Poisson partial differential equation from pointwise
state observations.  In what follows we describe the forward and
inverse problems setup, present the prior and the likelihood
distributions for the Bayesian inverse problem, and derive
the expressions for the gradient and Hessian action
using the standard Lagrangian approach as described in Section~\ref{sec:InfDimIP}.
The forward model is formulated as
follows
\begin{equation}\label{equ:poi}
  \begin{split}
    -\grad \cdot (\Exp{m} \grad u) &= f \quad \text{ in }\D, \\
    u  &= g \quad \text{ on } \GD, \\
    \Exp{m} \grad{u} \cdot \vec{n} &= h \quad \text{ on } \GN,
  \end{split}
\end{equation}
where $\D \subset \R^d$ ($d=2,3$) is an open bounded domain with
sufficiently smooth boundary $\Gamma = \GD \cup \GN$, $\GD \cap \GN =
\emptyset$.  Here, $u$ is the state variable, $f\in L^2(\D)$ is a
source term, and $g\in H^{1/2}(\GD)$ and $h\in L^2(\GN)$ are Dirichlet
and Neumann boundary data, respectively.  We define the spaces,
\begin{align*}
    \Vg = \{ v \in H^1(\D) : \restr{v}{\GD} = g\}, \quad
    \V =  \{ v \in H^1(\D) : \restr{v}{\GD} = 0\},
\end{align*}
where $H^1(\D)$ is the Sobolev space of functions whose derivatives
are in $L^2(\D)$.  Then, the weak form of \eqref{equ:poi} reads as
follows: find $u \in \Vg$ such that
\begin{equation}\label{eq:poisson_weak}
\ip{\Exp{m} \grad{u}}{\grad{p}} = \ip{f}{p} + \ip{h}{p}_{\GN}, \quad
\forall p \in \V.
\end{equation}
Here $\ip{\cdot}{\cdot}$ and $\ip{\cdot}{\cdot}_{\GN}$ denote the
standard inner products in $L^2(\D)$ and $L^2(\GN)$, respectively.

\subsubsection{Prior and noise models}
\label{sec:prob1_prior_and_noise}

We take the prior as a Gaussian distribution $\mathcal{N}(\iparpr,
\Cprior)$, with mean $\iparpr$ and covariance $\Cprior =
\mathcal{A}^{-2}$ following~\cite{Stuart10}.  $\mathcal{A}$
is a differential operator with domain $\mathcal{M} := H^1(\D)$ and action
\begin{equation}\label{eq:PoissonPrior}
  \mathcal{A}\ipar = \left\{ 
\begin{array}{ll}  
- \gamma \nabla \cdot \left(  \boldsymbol{\Theta} \grad \ipar \right) + \delta \ipar & \text{ in } \mathcal{D}\\
 \boldsymbol{\Theta} \grad \ipar \cdot \boldsymbol{n} + \beta \ipar  & \text{ on } \partial \mathcal{D},
 \end{array}
\right.
\end{equation}
where $\beta \propto \sqrt{\gamma\delta}$ is the optimal Robin coefficient derived in~\cite{DaonStadler18,RoininenHuttunenLasanen14}
to minimize boundary artifacts, and $\boldsymbol{\Theta}$ is an s.p.d.\! anisotropic tensor of the form
$$ \boldsymbol{\Theta} =
\begin{bmatrix}
\theta_1 \sin(\alpha)^2 & (\theta_1-\theta_2) \sin(\alpha) \cos{\alpha} \\
(\theta_1-\theta_2) \sin(\alpha) \cos{\alpha} & \theta_2 \cos(\alpha)^2
\end{bmatrix}. $$
In \figureref~\ref{fig:prior}
(left), we show the prior mean $\iparpr$ and three random draws from
the prior distribution with $\gamma = 0.1$, $\delta = 0.5$, $\alpha = \frac{\pi}{4}$,
$\theta_1 =2$, $\theta_2 = 0.5$.

Next, we specify the log-likelihood (data misfit) functional. We
denote with $\obs\in \R^q$ the vector of (noisy) pointwise
observations of the state $u$ at $q = 50$ random locations uniformly
distributed in $\D_{\rm obs} := [0.1, 0.9] \times [0.1, 0.5]$
($\D_{\rm obs} \subset \D$).  That is,
\[\obs = \mathcal{B} u + \vec{\eta},
\]
where $\mathcal{B}:
\Vg \mapsto \mathbb{R}^q$ is a linear observation operator,
a sum of delta functions to be specific, that extracts measurements
from $u$.  The measurement noise vector $\vec{\eta}$ is a multivariate
Gaussian variable with mean $\vec{0}$ and covariance
$\mat{\Gamma}_{\rm noise} = \sigma^2 \mat{I}$, where $\sigma = 0.01$,
and $\mat{I} \in \mathbb{R}^{q\times q}$.

\def \pos {0.43\columnwidth}
\begin{center}
\begin{figure}[tb]\centering
\begin{tikzpicture}
  \node (1) at (0*\pos-0.5*\pos, 0*\pos){\includegraphics[height=.23\textwidth]{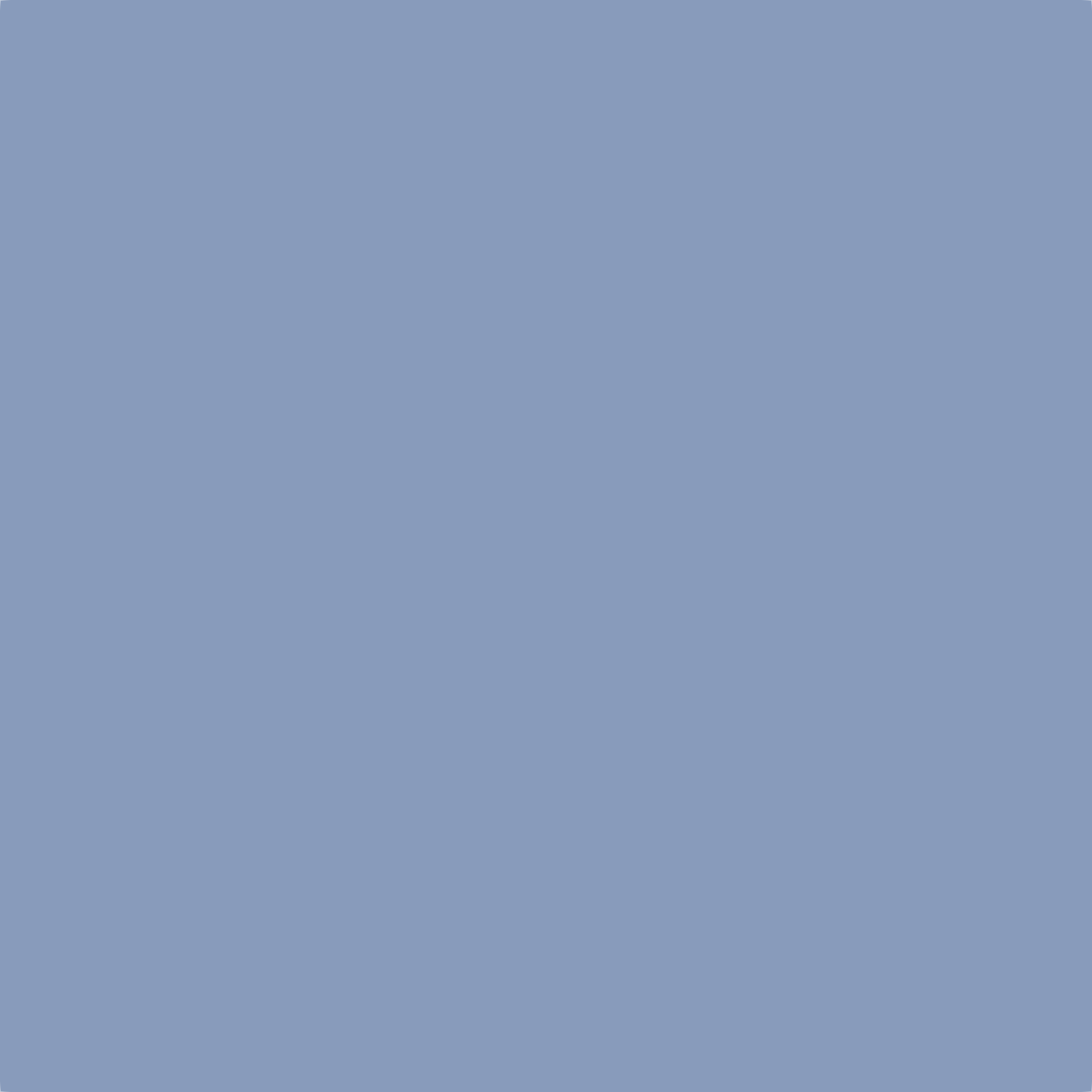}};\hspace{.15cm}
  \node at (-0.5*\pos-0.25*\pos, -0.2*\pos) {\sf a)};

  \node (2) at (0.05*\pos, 0*\pos){\includegraphics[height=.23\textwidth]{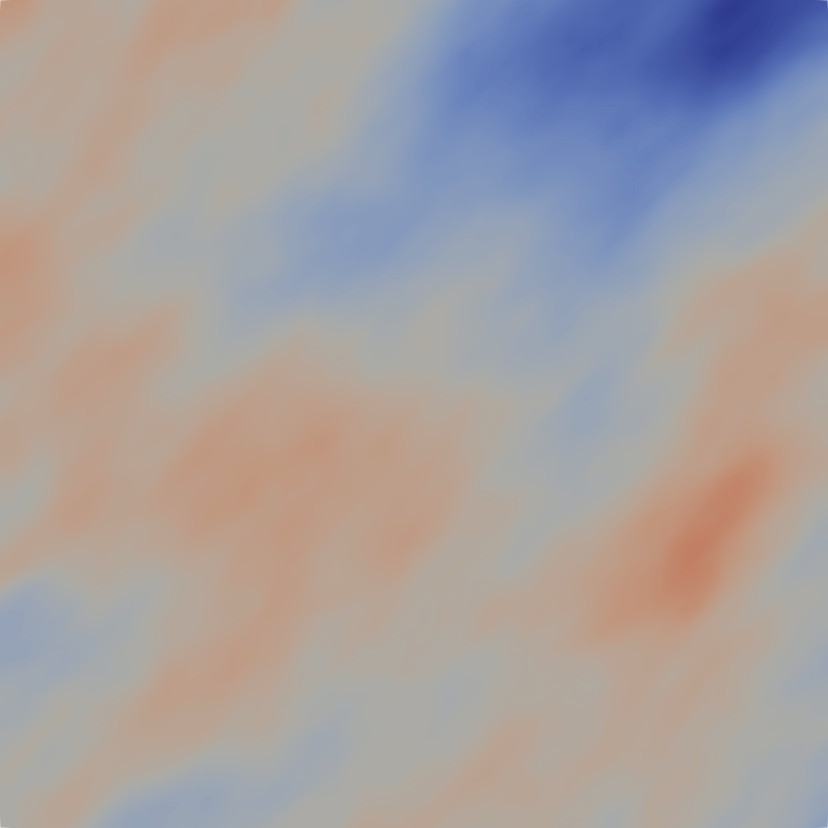}};\hspace{.15cm}
  \node at (0.05*\pos-0.25*\pos,  -0.2*\pos) {\sf b)};

  \node (3) at (0.6*\pos, 0*\pos){\includegraphics[height=.23\textwidth]{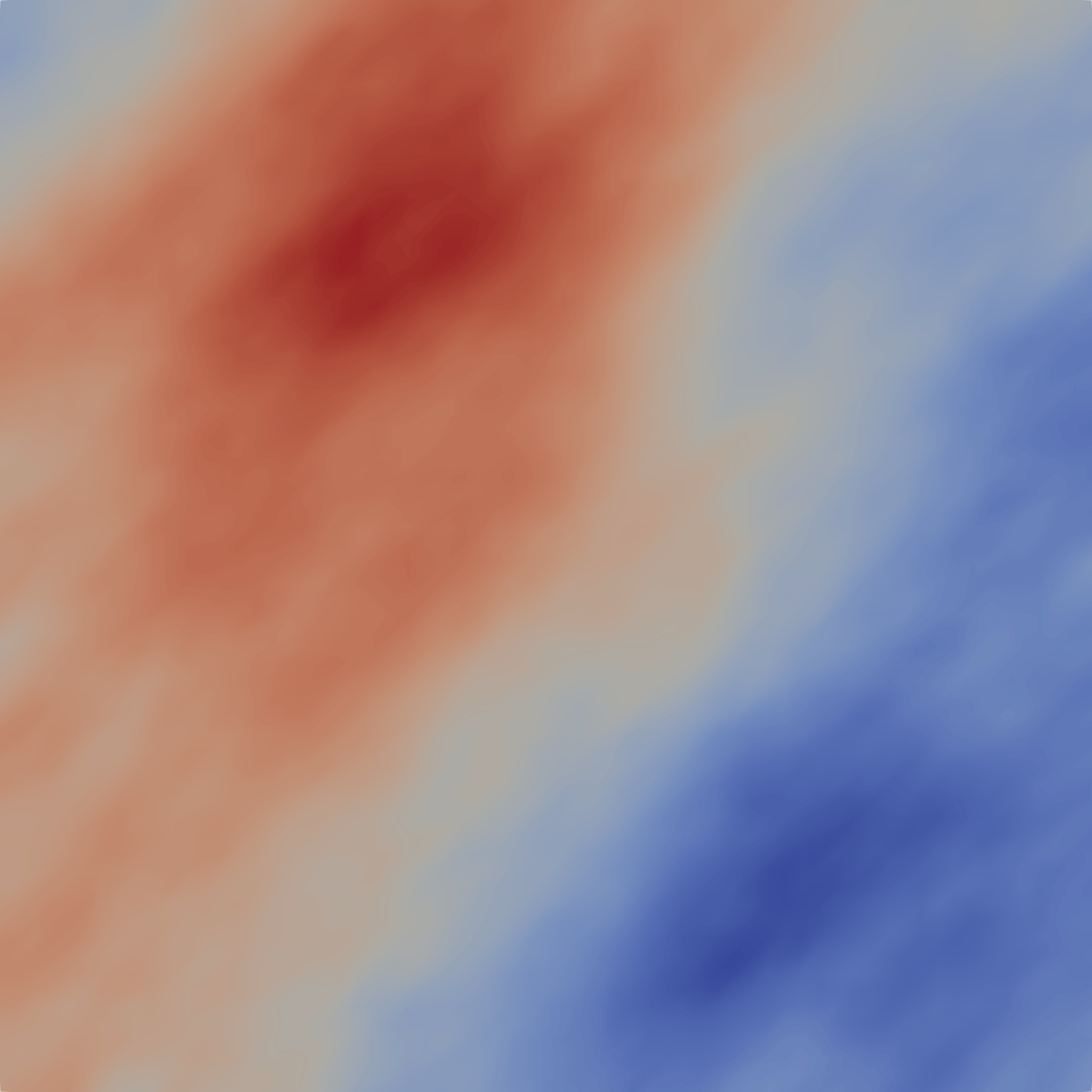}};\hspace{.15cm}
  \node at (0.6*\pos-0.25*\pos, -0.2*\pos) {\sf c)};

  \node (4) at (1.15*\pos, 0*\pos){\includegraphics[height=.23\textwidth]{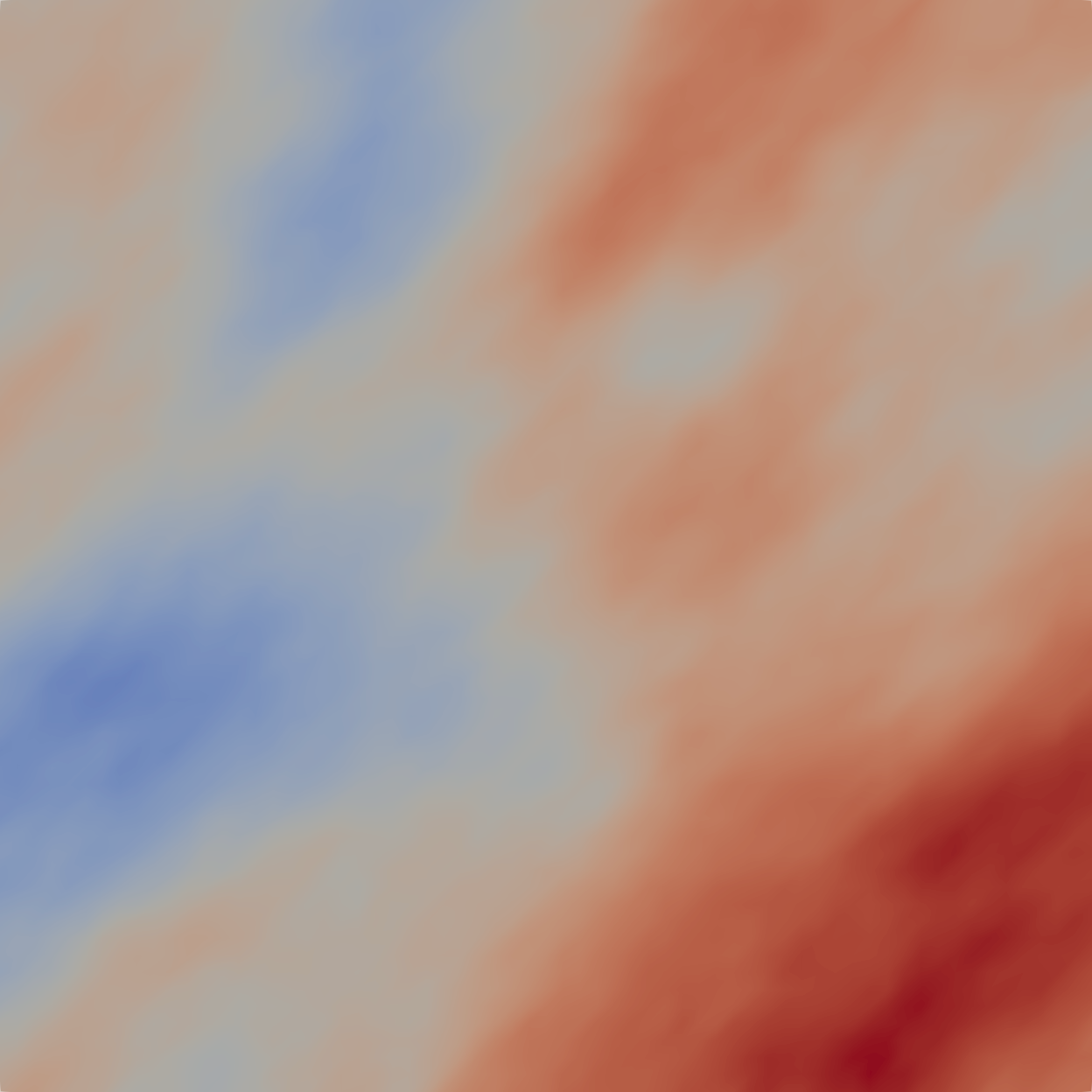}};
  \node at (1.15*\pos-0.22*\pos, -0.2*\pos) {\sf d)};
\end{tikzpicture}
\caption{Prior mean $\iparpr$ (a), and samples drawn from the prior
  distribution (b)--(d).}
\label{fig:prior}
\end{figure}
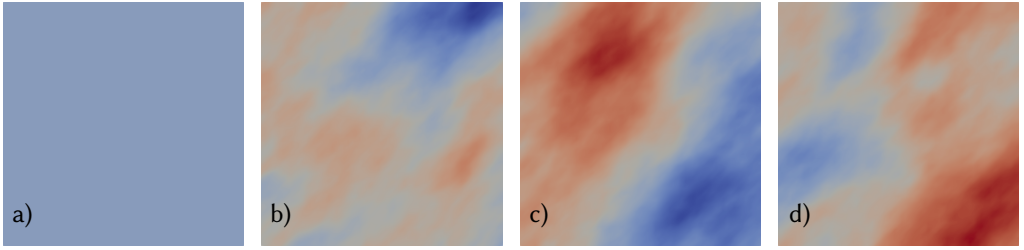
\end{center}

\subsubsection{The MAP point}
\label{subsec:map}
To find the MAP point we solve the following variational nonlinear
least-squares optimization problem
\begin{equation}\label{eq:cost-elliptic}
\begin{aligned}
& \min_{\ipar \in \iparspace} \! \J (\ipar) := \frac{1}{2} \left\| \mathcal{B}u(\ipar) -\obs
  \right\|^2_{\ncov^{-1}}
  + \frac{1}{2} \left\|{\ipar - \iparpr}\right\|^2_{\Cprior^{-1}},
\end{aligned}
\end{equation}
where the state variable $u$ is the solution to~\eqref{equ:poi},
$\iparpr$ is the prior mean of the log coefficient field $\ipar$, and
$\obs\in \R^q$ is a given data vector.  To solve this optimization
problem we use the inexact Newton-CG algorithm in Algorthim
\ref{alg:InexactNewtonCG}, which requires gradient and Hessian
information. These are automatically computed by \hip applying symbolic 
differentiation to the variational form of the forward problem \eqref{eq:poisson_weak};
we also refer to Appendix~\ref{app-sec:MAP} where the gradient and Hessian-apply expressions
are derived using the Lagrangian formalism. We note
that the use of CG to solve the resulting Newton system does not
require computing the Hessian operator by itself but only its action
in a given direction.

\subsubsection{Numerical results}
For the forward Poisson problem~\eqref{equ:poi}, no source term,
(i.e., $f = 0$) and no normal flux on $\GN:=\{0,1\}\times (0,1)$
(i.e., the homogeneous Neumann condition $\Exp{\ipar} \grad u \cdot
\vec{n} = 0$ on $\GN$) are imposed.  Dirichlet conditions are
prescribed on the top and bottom boundaries, in particular $u = 1$ on
$(0,1)\times \{1\}$ and $u = 0$ on $(0,1)\times \{0\}$. This Dirichlet
part of the boundary is denoted by $\GD:=(0,1)\times\{0,1\}$.  In
\figureref~\ref{fig:forwardprob}, we show the \emph{true} parameter field
used in our numerical tests, and the corresponding state field.
We used quadratic finite elements to discretize the state and adjoint
variables and linear elements for the parameter. The degrees of
freedom for the state and parameter were 16641 and 4225,
respectively.

\def \pos {0.49\columnwidth}
\begin{figure}[t]\centering
  \begin{tikzpicture}
    \node (1) at (0.0*\pos-0.75*\pos, 0.0*\pos){
      \includegraphics[height=.35\textwidth, trim=0 0 0 0, clip=true]{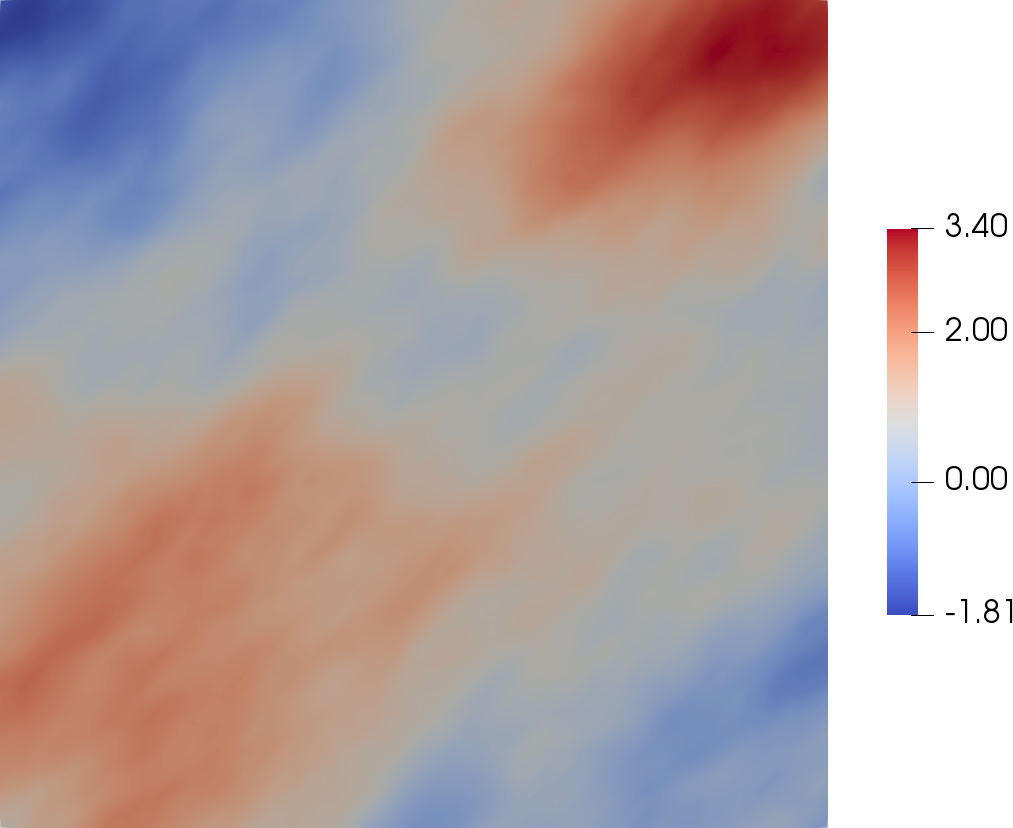}
    };
    \node (2) at (0.5*\pos-0.25*\pos, 0.0*\pos){
      \includegraphics[height=.35\textwidth, trim=0 0 0 0, clip=true]{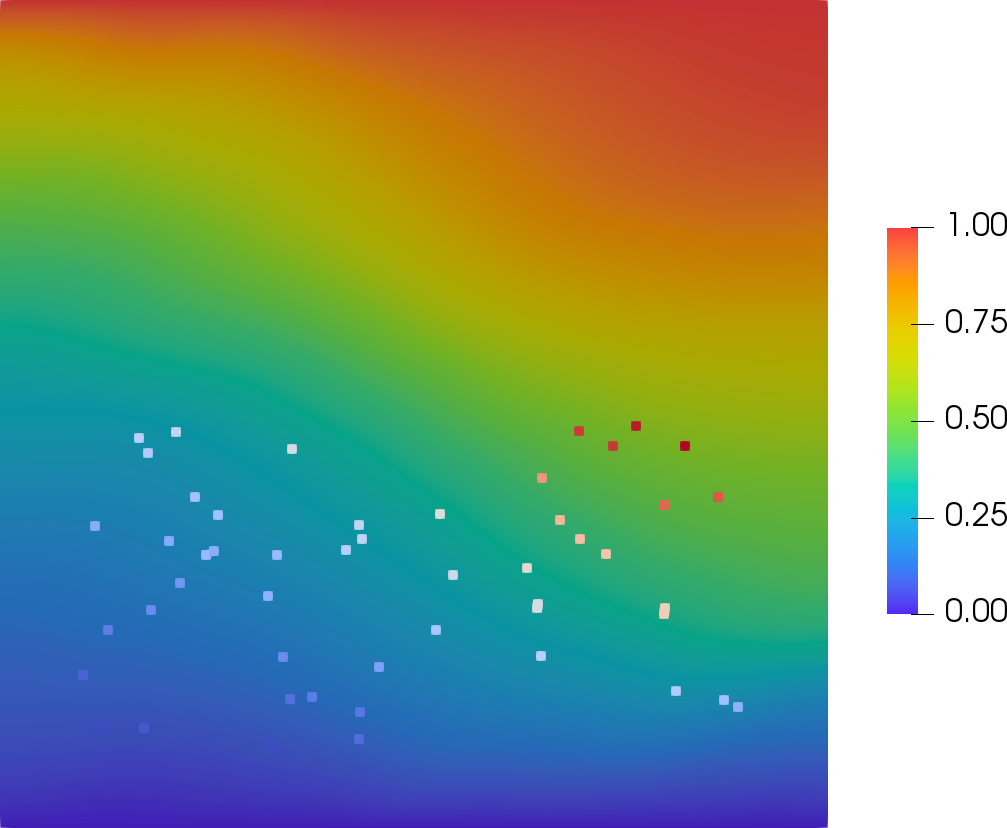}
    };
    \node at (0.0*\pos-1.15*\pos, -0.3*\pos) {\sf a)};
    \node at (0.7*\pos-0.83*\pos, -0.3*\pos) {\sf b)};
  \end{tikzpicture}

  \caption{The true parameter field $\ipart$ (left)
    and the state $u$ obtained by solving the
    forward PDE with $\ipart$ (right). The squares in (b) represent
    locations of the $q=50$ randomly chosen observation points and their
    color corresponds to the observed noisy data~$\obs$.}
\label{fig:forwardprob}
\end{figure}

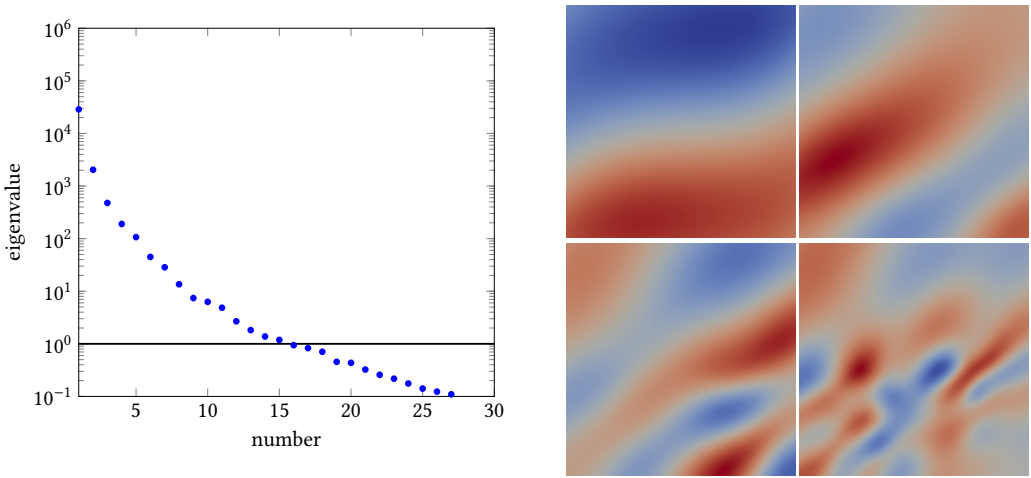
\begin{figure}[t]
  \centering
  \resizebox{\columnwidth}{!}{\input{tikz/spec-and-evecs-possion.tex}}
  \caption{Left: Log-linear plot of first 30 out of 4225 eigenvalues
    of the data misfit Hessian for the generalized eigenproblem
    \eqref{eq:eigenproblem_1}. The low-rank based approximation
    captures the dominant, data-informed portion of the spectrum. The
    eigenvalues are truncated at 0.07. Right: Prior-orthogonal
    eigenvectors of the data misfit Hessian corresponding (from left
    to right) to the 1st, 3rd, 8th, and 27th eigenvalues. Note that
    eigenvectors corresponding to smaller eigenvalues are increasingly
    more oscillatory (and thus inform smaller length scales of the
    parameter field) but are also increasingly less informed by the
    data.}
\label{fig:spectrum-evecs}
\end{figure}

Next we study the spectrum of the data misfit Hessian evaluated at the
MAP point.  \figureref~\ref{fig:spectrum-evecs} (left) shows a logarithmic
plot of the eigenvalues of the generalized symmetric eigenproblem
\[
\Hmisfit \vec{v}_i = \lambda_i \mat{R} \vec{v}_i, \quad \lambda_1 \geq
\lambda_2 \geq \ldots \geq \lambda_n,
\]
where $\Hmisfit$ and $\mat{R}$ stems from the discretization (with
respect to the Euclidean inner product of the data misfit Hessian and
prior precision (cf Eq. \eqref{eq:eigenproblem_1}). This plot shows
that the spectrum decays rapidly. As seen in~\eqref{eq:sm}, an
accurate low-rank based approximation of the inverse Hessian can be
obtained by neglecting eigenvalues that are small compared to 1. Thus,
retaining around $r=30$ eigenvectors out of 4225 (i.e., the dimension
of parameter space) appears to be sufficient. We stress that $r$ is
strictly less than the number $q=50$ of observation, reflecting
redundancy in the data.  We note that the cost of obtaining this
low-rank based approximation, measured in the number of forward and
adjoint PDE solves, is $2(r+l)$, where $r+l$ is the number of random
vectors. Here, $l = 20$ is an oversampling parameter used to ensure
the accurate computation of the most significant
eigenvalues/eigenvectors. The corresponding retained eigenvectors are
those modes in parameter space that are simultaneously well-informed
by the data and assigned high probability by the prior.
\figureref~\ref{fig:spectrum-evecs} (right) displays several of these
eigenvectors.

\begin{figure}[t!]
  \centering
  \begin{tikzpicture}
    \node (l) at (0,0.5){
      \includegraphics[width=0.35\columnwidth]{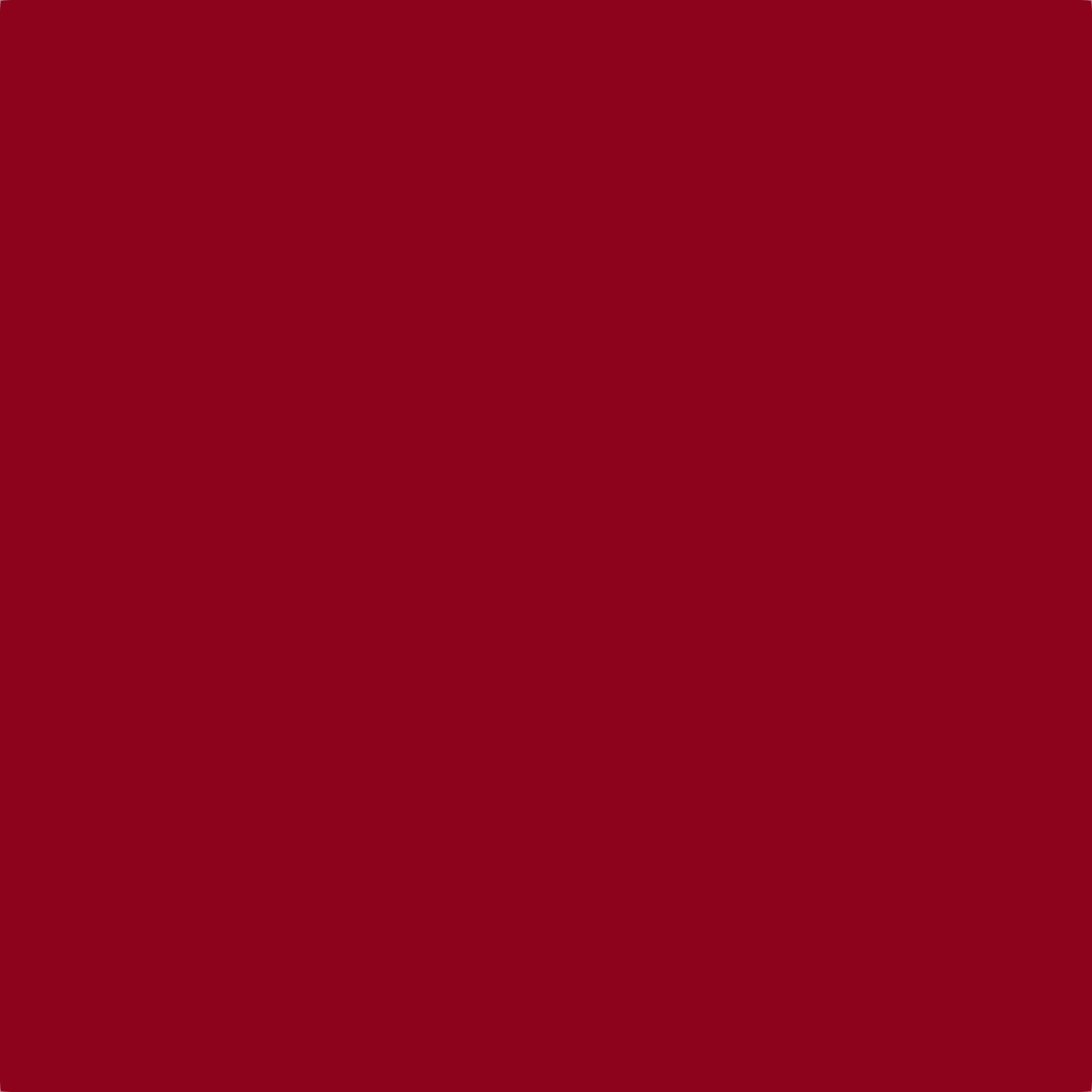}
      \hspace{0.1in}
      \includegraphics[width=0.35\columnwidth]{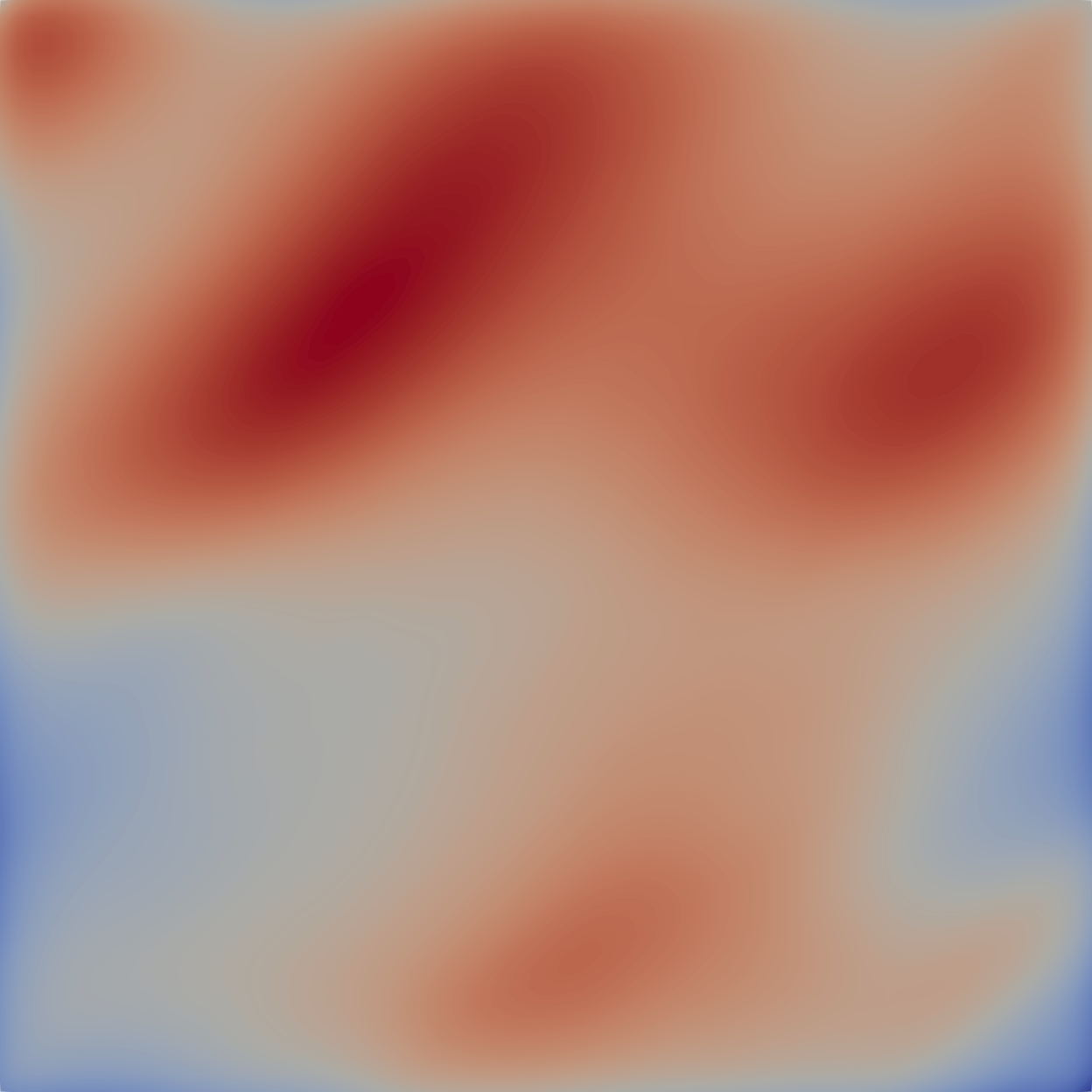}};
    \node (r) at (5.8,0.5){
      \includegraphics[width=0.05\columnwidth]{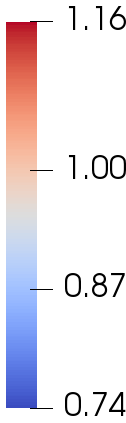}};
  \end{tikzpicture}
\caption{Pointwise variance of the prior distribution (left) and the
  Laplace approximation of the posterior distribution (right). Note
  that uncertainty is mostly reduced in the lower half of the domain
  where data is measured.}
\label{fig:variance}
\end{figure}

\figureref~\ref{fig:variance} depicts the prior and posterior pointwise
variances computed using \eqref{eq:prior_mvar} and
\eqref{eq:post_mvar} with $\overline{r} = 300$ and $r=50$,
respectively.  One observes that the uncertainty is vastly reduced in
the bottom half of the domain, which is expected given that
observations are present only on the lower half of $\D$.  In
\figureref~\ref{fig:post} we show the MAP point (a) and samples from the
Laplace approximation~\eqref{eq:lapl_discrete} of the posterior
probability density function (b)-(d). These samples were obtained by
first computing samples from the prior distribution---shown in
\figureref~\ref{fig:prior}---according to~\eqref{eq:sample_prior}, and
then applying~\eqref{eq:sample_post}.  The variance reduction between 
posterior samples in \figureref~\ref{fig:post} and prior samples in \figureref~\ref{fig:prior}
reflects the information gained from the data in solving the inverse
problem. In addition, we note that the MAP point resembles the truth
better in the lower half of the domain where data are available.
The presence (or absence) of data also affects the
posterior samples, in fact, we observe higher variability in the
upper half of the domain, where there is no data.

\def \pos {0.43\columnwidth}
\begin{center}
\begin{figure}[tb]\centering
\begin{tikzpicture}
  \node (1) at (0*\pos-0.5*\pos, 0*\pos){\includegraphics[height=.23\textwidth]{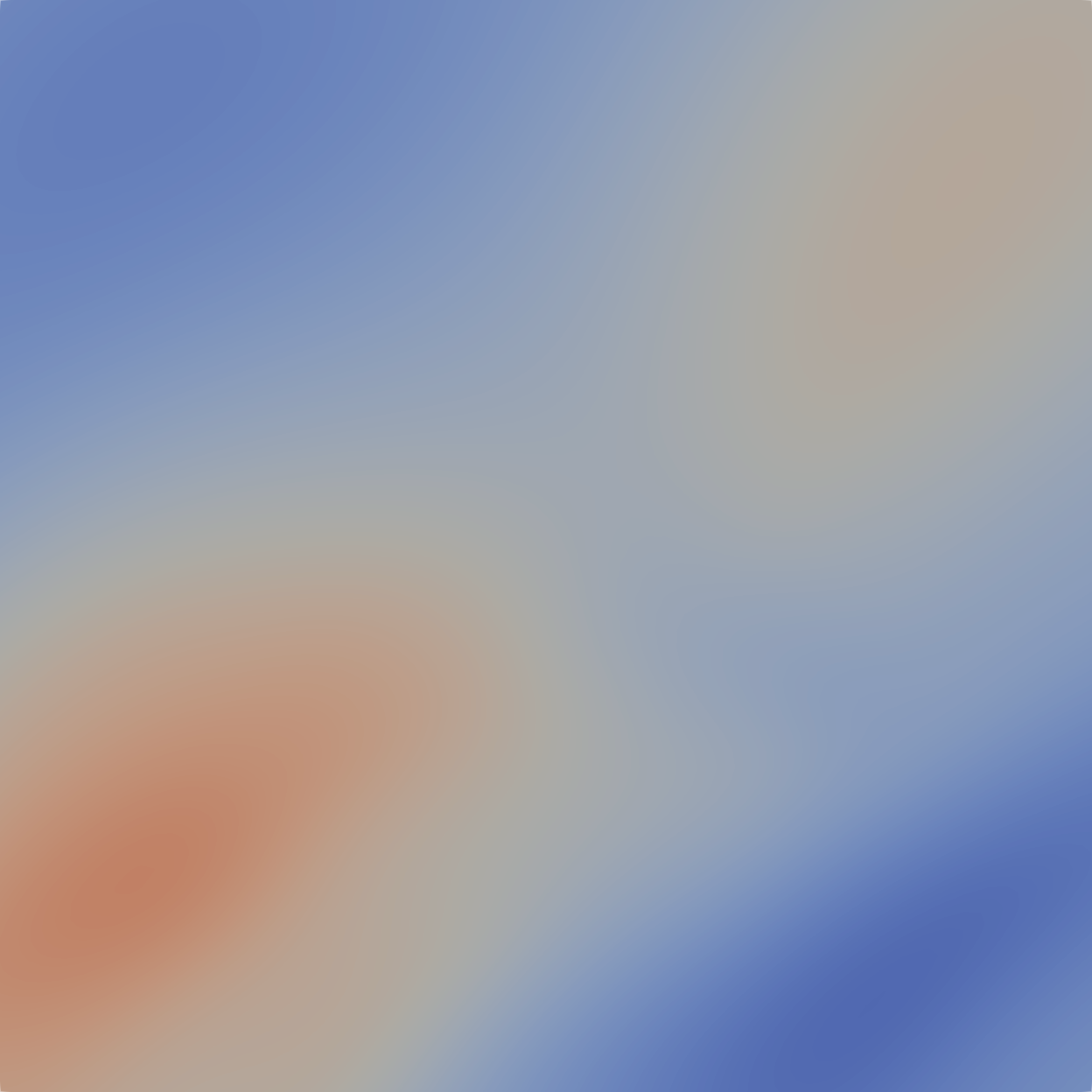}};\hspace{.15cm}
  \node at (-0.5*\pos-0.25*\pos, -0.2*\pos) {\sf a)};

  \node (2) at (0.05*\pos, 0*\pos){\includegraphics[height=.23\textwidth]{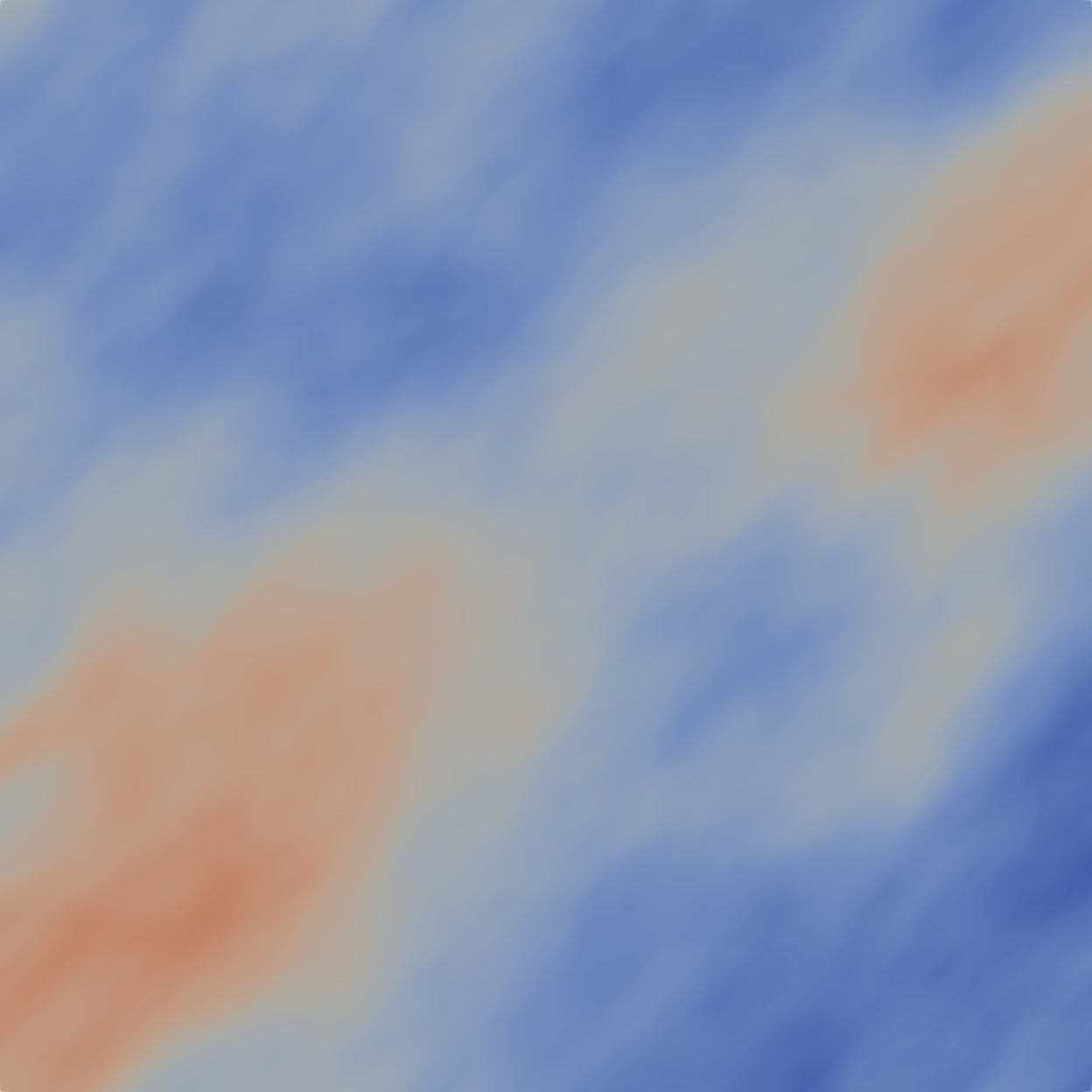}};\hspace{.15cm}
  \node at (0.05*\pos-0.25*\pos,  -0.2*\pos) {\sf b)};

  \node (3) at (0.6*\pos, 0*\pos){\includegraphics[height=.23\textwidth]{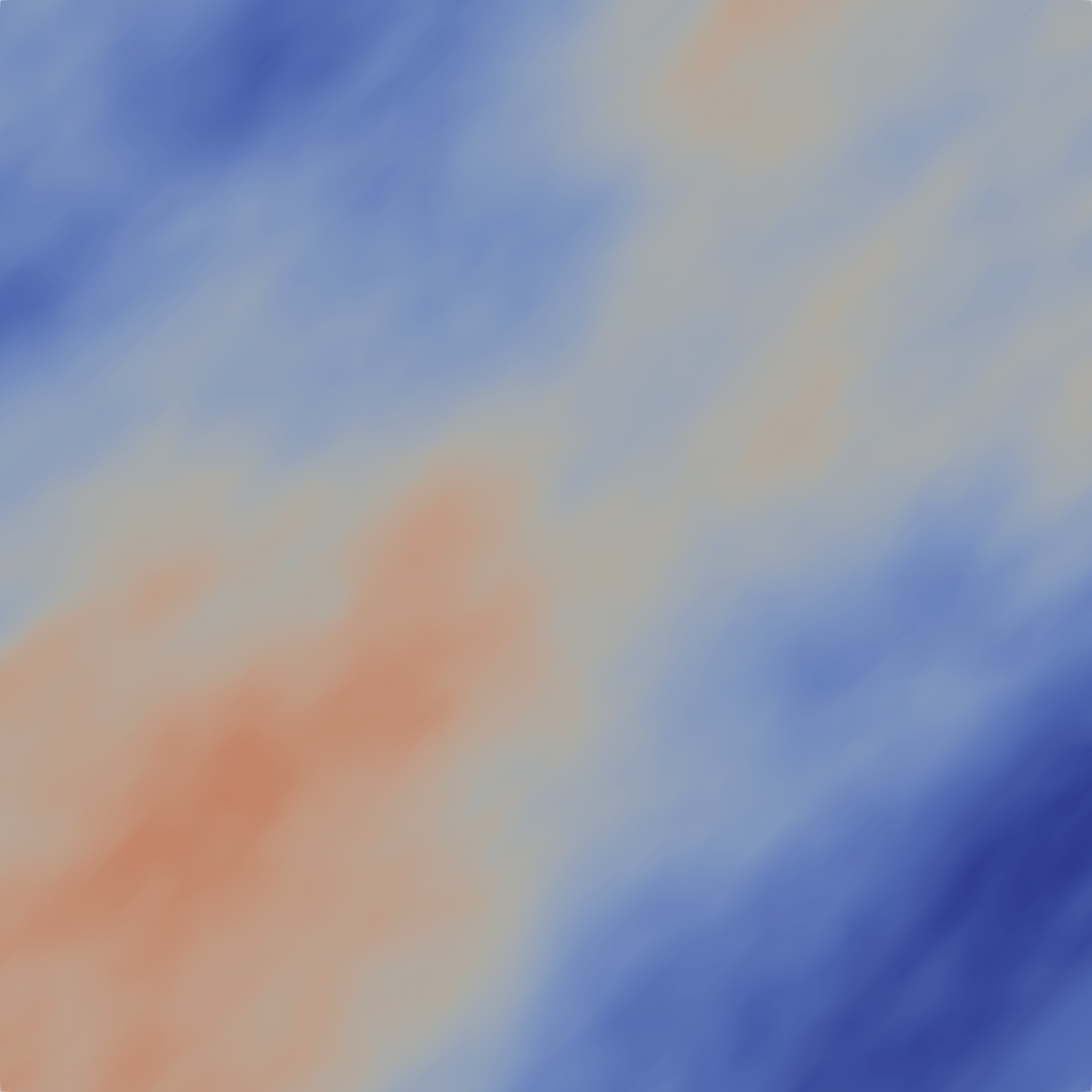}};\hspace{.15cm}
  \node at (0.6*\pos-0.25*\pos, -0.2*\pos) {\sf c)};

  \node (4) at (1.15*\pos, 0*\pos){\includegraphics[height=.23\textwidth]{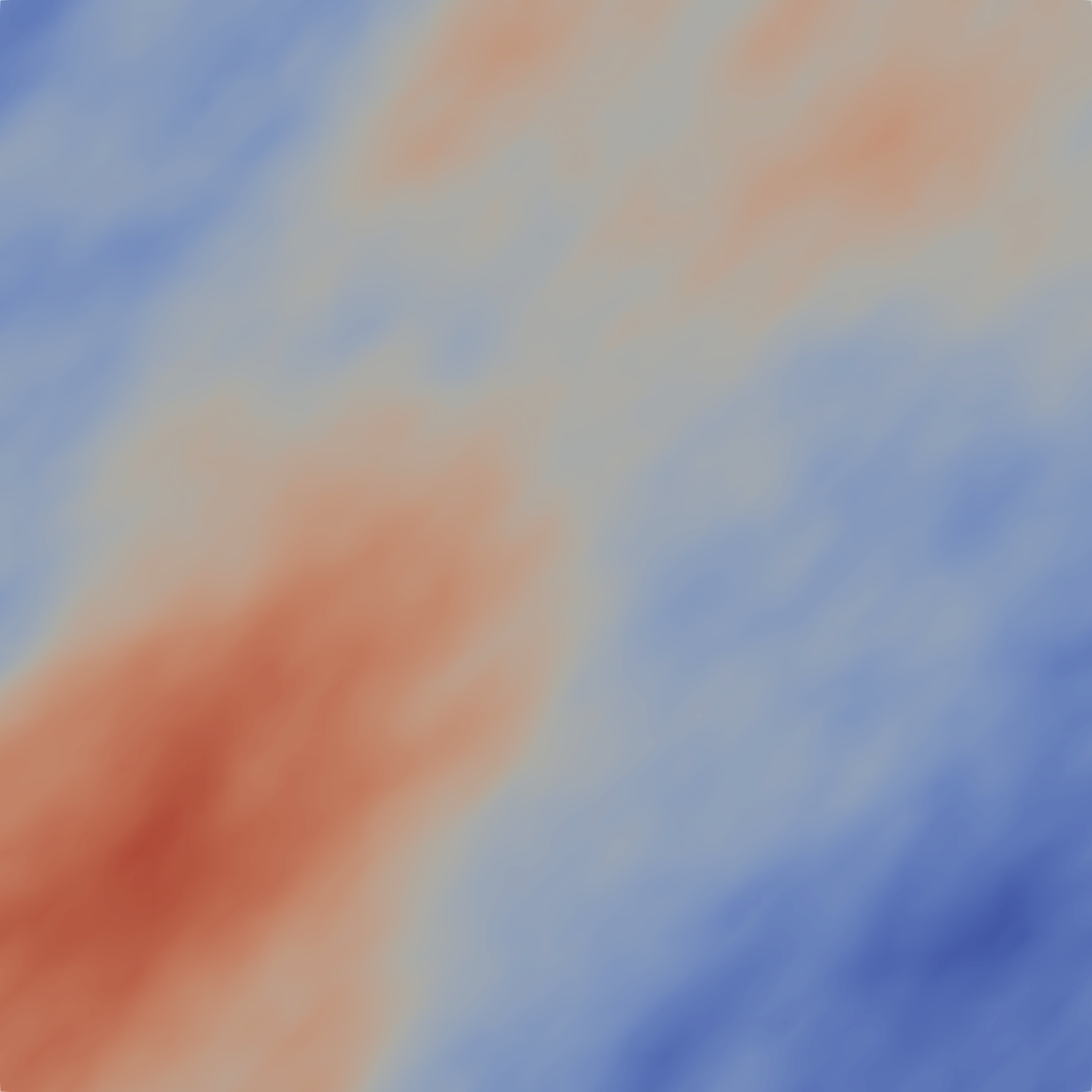}};
  \node at (1.15*\pos-0.22*\pos, -0.2*\pos) {\sf d)};

\end{tikzpicture}
\caption{The MAP point (a) and samples drawn from the Laplace approximation of
  posterior distribution (b)--(d).}
\label{fig:post}
\end{figure}
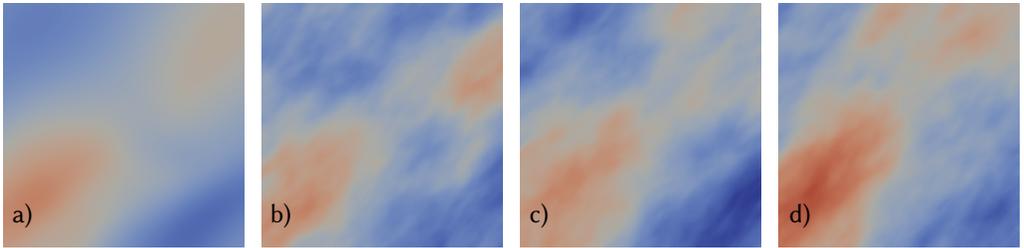
\end{center}

\subsection{Inversion for the initial condition in an advection-diffusion PDE}
\label{sec:advdiff}
\begin{figure}[t]\centering
  \includegraphics[width=0.5\columnwidth]{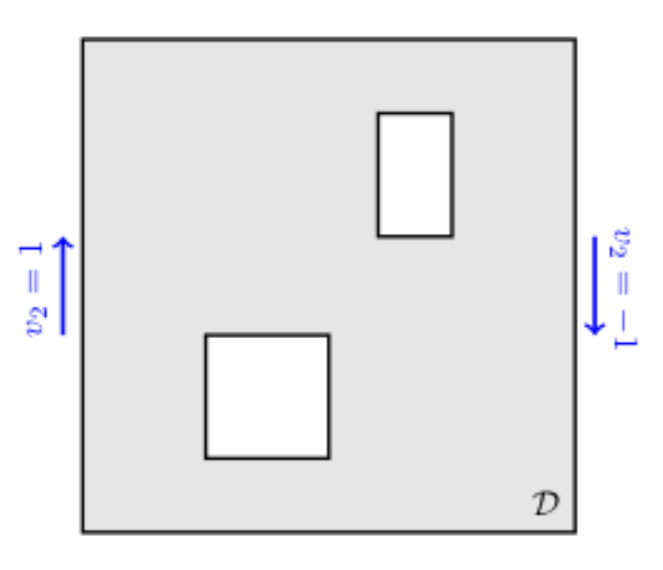}
  \includegraphics[width=0.43\columnwidth]{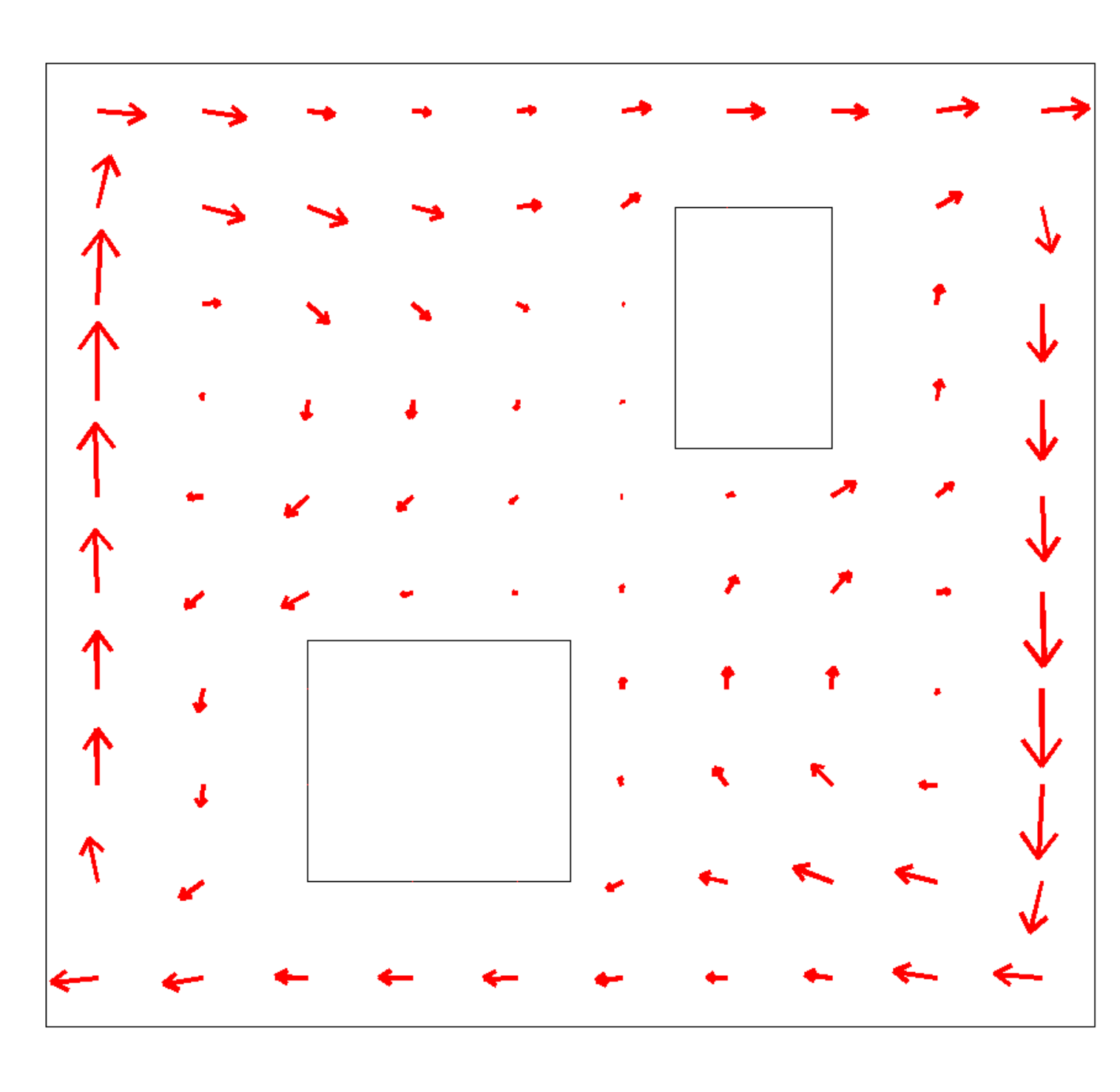}
  \caption{Left: Sketch of domain for the advective-diffusive inverse
    transport problem~\eqref{eq:ad} showing imposed velocities used to
    compute the velocity field $\vv$. Right: The velocity field
    computed from the solution of the Navier-Stokes equations and
    subsampled on a coarser grid for visualization purposed.}
  \label{fig:ad_domain}
\end{figure}

Here we consider a time-dependent advection-diffusion equation for
which we seek to infer an unknown initial contaminant field from pointwise
measurements of its concentration. The problem description below closely follows
the one in \cite{PetraStadler11}.
The PDE in the parameter-to-observable map models diffusive transport
in a domain $\D \subset \R^d$, which is depicted in
\figureref~\ref{fig:ad_domain}.  The domain boundaries $\del \D$ include the
outer boundaries as well as the internal boundaries of the rectangles,
which represent buildings.  The parameter-to-observable map $\iFF$ maps an
initial condition $\ipar \in L^2(\D)$ to pointwise spatiotemporal
observations of the concentration field $u(\vec{x}, t)$ through solution of the advection-diffusion equation
 given by
\begin{equation}\label{eq:ad}
  \begin{aligned}
    u_t - \kappa\Delta u + \vec{v}\cdot\nabla u &= 0 & \quad&\text{in
    }\D\times (0,T), \\
    u(\cdot, 0) &= \ipar  &&\text{in } \D , \\
    \kappa\nabla u\cdot \vec{n} &= 0 &&\text{on } \partial\D \times (0,T).
  \end{aligned}
\end{equation}
Here, $\kappa > 0$ is a diffusivity coefficient, and $T > 0$ is the
final time of observations. 
The velocity field $\vec{v}$, shown in \figureref~\ref{fig:ad_domain}
(right), is computed by solving the steady-state Navier-Stokes equations
for a two dimensional flow with Reynolds number 50 and
boundary conditions as in \figureref~\ref{fig:ad_domain} (left); see \cite{PetraStadler11} 
for details.
The time evolution of the state variable $u$ from a given initial
condition $\ipar$ is illustrated in
\figureref~\ref{fig:add-diff-state-attrueandmap} (top).

To derive the weak formulation of \eqref{eq:ad}, we define the spaces
\begin{align*}
    \mathcal{V} :=  \{ v \in H^1(\D), \mbox{ for each } t \in (0, T)\}, \text{ and }  \iparspace := H^1(\D).
\end{align*}
Then, the weak form of~\eqref{eq:ad} reads as follows: Find $u \in \mathcal{V}$ such that
\begin{equation}\label{eq:ad_weak}
\int_0^T\!\int_\D (u_t + \vec v\cdot\nabla u)p\,d\x\,dt +
\int_0^T\!\int_\D\kappa\nabla u \cdot \nabla p\,d\x\,dt + \int_\D
(u(\vec{x}, 0) - \ipar )p_0\,d\x = 0, 
\end{equation}
$\forall p \in \mathcal{V},\, p_0 \in \iparspace$.  Above, the initial
condition $u(\vec{x}, 0) = \ipar$ is imposed weakly by means of the
test function $p_0 \in \iparspace$.

\subsubsection{The noise and prior models}
\label{subsubsec:ad-prior-noise}

We consider the problem of inferring the initial condition $\ipar$ in
\eqref{eq:ad} from pointwise noisy observations $\boldsymbol{d}_i \in
\mathbb{R}^{n_t}$ ($i=1,\ldots,n_s$) of the state $u$ at $n_s$
discrete time samples $t_i$ in interval $[T_1, T] \subset [0, T]$, and
$n_t$ locations in space.  We assume that the observations
$\boldsymbol{d}_i$ are perturbed with \emph{i.i.d.} Gaussian additive
noise with variance $\sigma^2$.  To construct the prior, we assume a
constant mean and define $\Cprior := \mathcal{A}^{-2} = (-\gamma \gbf
\Delta + \delta I)^{-2}$, equipped with Robin boundary conditions
$\gamma \nabla m \cdot\boldsymbol{n} + \beta m$ on $\partial \D$.  The
parameters $\gamma, \delta > 0$ control the correlation length and
variance of the prior operator; here we take $\gamma = 1$ and $\delta
= 8$. The Robin coefficient $\beta$ is chosen as in
\cite{DaonStadler18,RoininenHuttunenLasanen14} to reduce boundary
artifacts.

\subsubsection{The MAP point}
\label{subsubsec:ad-map}

To compute the MAP point we minimize the negative log-posterior,
defined in general in~\eqref{eq:map_continuum}, which---for Gaussian
prior and noise---is analogous to the least-squares functional
minimized in the solution of a deterministic inverse problem. For this
particular problem, this reads
\begin{align}\label{eq:ad_cost}
& \mathcal{J}(\ipar) := \frac{1}{2\sigma^2} \sum_{i=1}^{n_s}\int_{T_1}^T(\mathcal{B}u-\boldsymbol{d}_i)^2 \delta_{t_i} \,dt + \frac{1}{2} \left\|
  \Acal(\ipar - \iparpr \right)\|^2_{L^2(\D)},
\end{align}
where $\mathcal{B}:\V \mapsto \mathbb{R}^{n_t}$ is the interpolation
operators at the observation locations, $ \delta_{t_i}$ is the Dirac
delta functions at the observation time sample $t=t_i$ ($i = 1,
\ldots, n_s$), and $\sigma^2$ represents the noise level in the
observations $\boldsymbol{d}_i$, here taken $2.45 \times 10^{-7}$, and
$\iparpr =0$ is the prior mean. We use the conjugate gradient method to solve this (linear) inverse problem. The
derivation of the gradient and Hessian-apply is given in
Appendix~\ref{app:ad-map}.

\subsubsection{Numerical results}

Next, we present numerical results for the initial condition inverse
problem. The discretization of the forward and adjoint problems uses
an unstructured triangular mesh,
Galerkin finite elements with piecewise-quadratic globally-continuous polynomials, 
and an implicit Euler method for the time discretization.
Galerkin Least-Squares stabilization of the convective
term~\cite{HughesFrancaHulbert89a} is added to ensure stability of the discretization.
The space-time dimension of the state variable is 433880 (10847 spatial degree of freedom times 40 time steps),
and the dimension of the parameter space is 10847. The data dimension $q$ is 1200 with $n_t = 80$ measurement locations
and $n_s = 15$ time samples.

To illustrate properties of the forward problem,
\figureref~\ref{fig:add-diff-state-attrueandmap} show three snapshots in
time of the field $u$, using the advective velocity $\vec{v}$ from
\figureref~\ref{fig:ad_domain} with the ``true'' initial condition (top
row) and its MAP point estimate (bottom row), respectively. Next we
study the numerical rank of the prior-preconditioned data misfit
Hessian. Note that due to linearity of the parameter-to-observable map
$\iFF$, the prior-preconditioned data misfit Hessian is independent of
$\ipar$. \figureref~\ref{fig:ad_diff-spectrum-evecs} (left) shows a
logarithmic plot of the truncated spectra of the prior-preconditioned
data misfit Hessians for several observation time horizons. This plot
shows that the spectrum decays rapidly and, as expected, the decay is
faster when the observation time horizon is shorter (i.e., there are
fewer observations).  As seen in~\eqref{eq:sm}, an accurate low-rank
based approximation of the inverse Hessian can be obtained by
neglecting eigenvalues that are small compared to 1. Thus, retaining
around 70 eigenvectors out of 10847 appears to be sufficient for the
target problem with spatial and temporal observation points in the
interval $[1,4]$. These eigenvalues and the corresponding
prior-orthogonal eigenvectors (see the right panel in \figureref~
\ref{fig:ad_diff-spectrum-evecs}) were computed using the double pass
algorithm Algorithm \ref{alg:doublepass} with $r = 50$ and
oversampling parameter $l=10$.

\begin{figure}[ht]
  \centering
  \begin{tikzpicture}
    \node (l) at (0,0.5){
      \includegraphics[width=0.35\columnwidth]{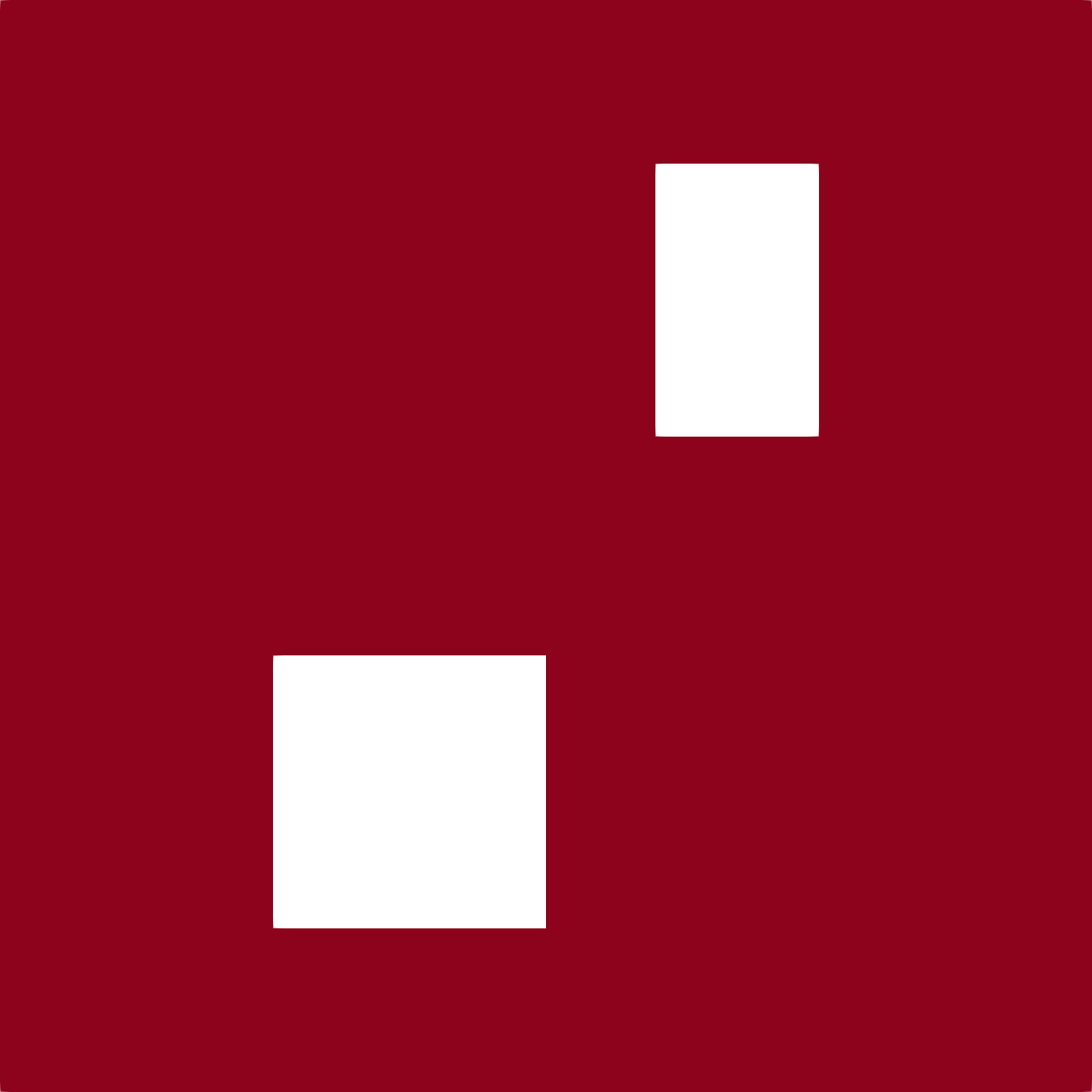}
      \hspace{0.1in}
      \includegraphics[width=0.35\columnwidth]{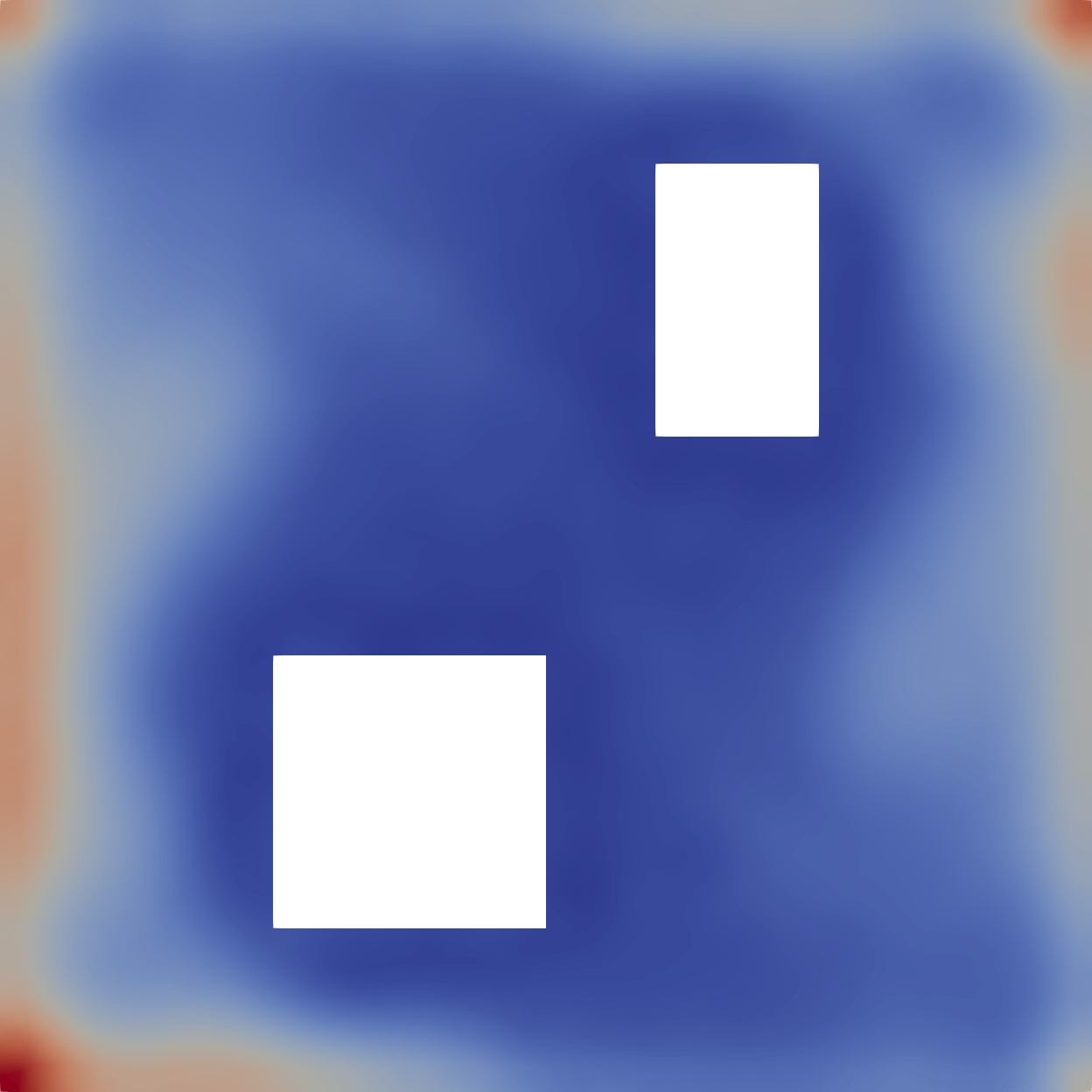}};
    \node (r) at (6.1,0.5){
      \includegraphics[width=0.1\columnwidth]{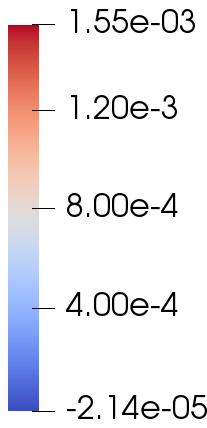}};
  \end{tikzpicture}
\caption{This figure shows the pointwise variance of the prior (left)
  and posterior (right) distributions.}
\label{fig:ad_diff-variance}
\end{figure}

Due to the linearity of the parameter-to-observable map and the choice
of a Gaussian prior and noise model, the posterior distribution is
also Gaussian whose mean coincides with the MAP point and the
covariance with the inverse of the Hessian evaluated at the MAP point.
Thus, for this problem, the Laplace approximation is the posterior
distribution.  \figureref~\ref{fig:ad_diff-variance} depicts the prior
and posterior pointwise variances. This figure shows that the
uncertainty is reduced everywhere in the domain, and that the
reduction is greatest near to the observations (at the boundaries of
the interior rectangles).  In \figureref~\ref{fig:add-diff-prior-post} we
show samples (of the initial condition) from the prior and from the
posterior, respectively. The difference between the two sets of
samples reflects the information gained from the data in solving the
inverse problem. The small differences in the parameter field $\ipar$
across the posterior samples (other than near the external boundaries)
demonstrate that there is small variability in the inferred
parameters, reflecting large uncertainty reduction.

\def \pos {0.43\columnwidth}
\begin{center}
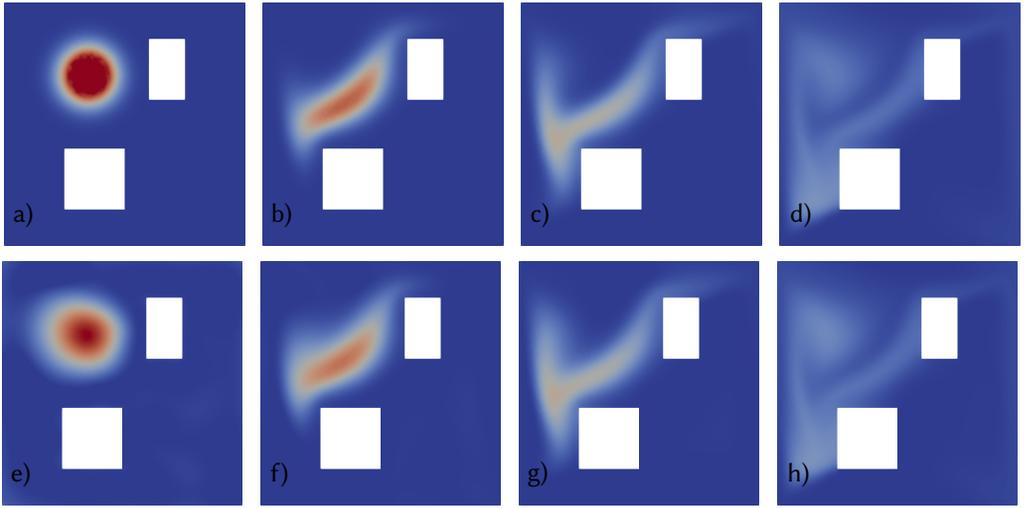
\begin{figure}[tb]\centering
\begin{tikzpicture}
  \node (1) at (0*\pos-0.5*\pos, 0*\pos){\includegraphics[height=.23\textwidth]{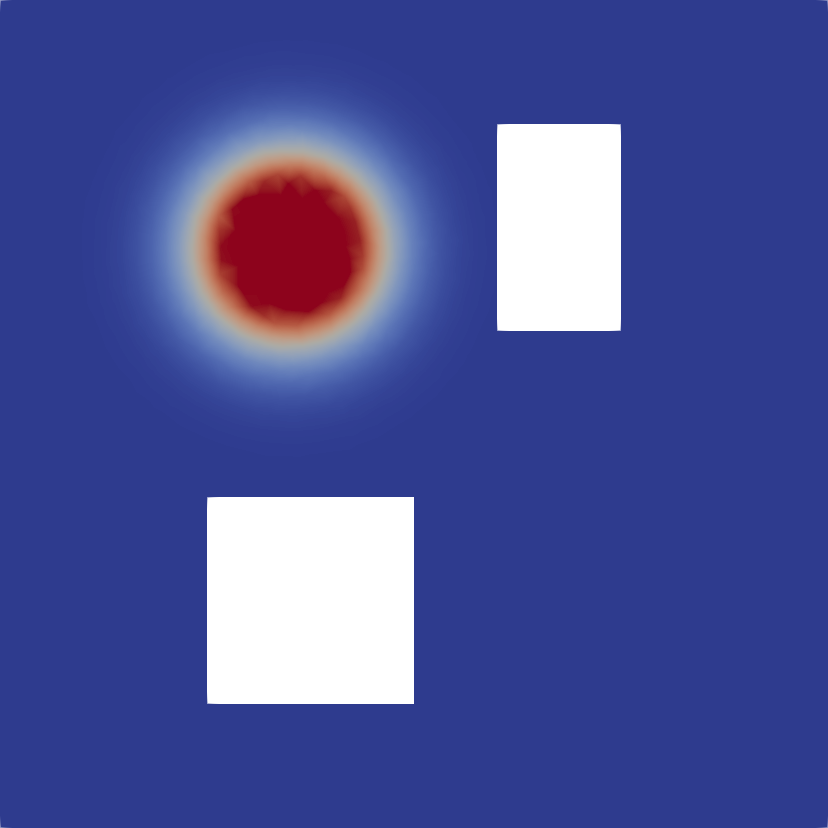}};\hspace{.15cm}
  \node at (-0.5*\pos-0.25*\pos, -0.2*\pos) {\sf a)};

  \node (2) at (0.05*\pos, 0*\pos){\includegraphics[height=.23\textwidth]{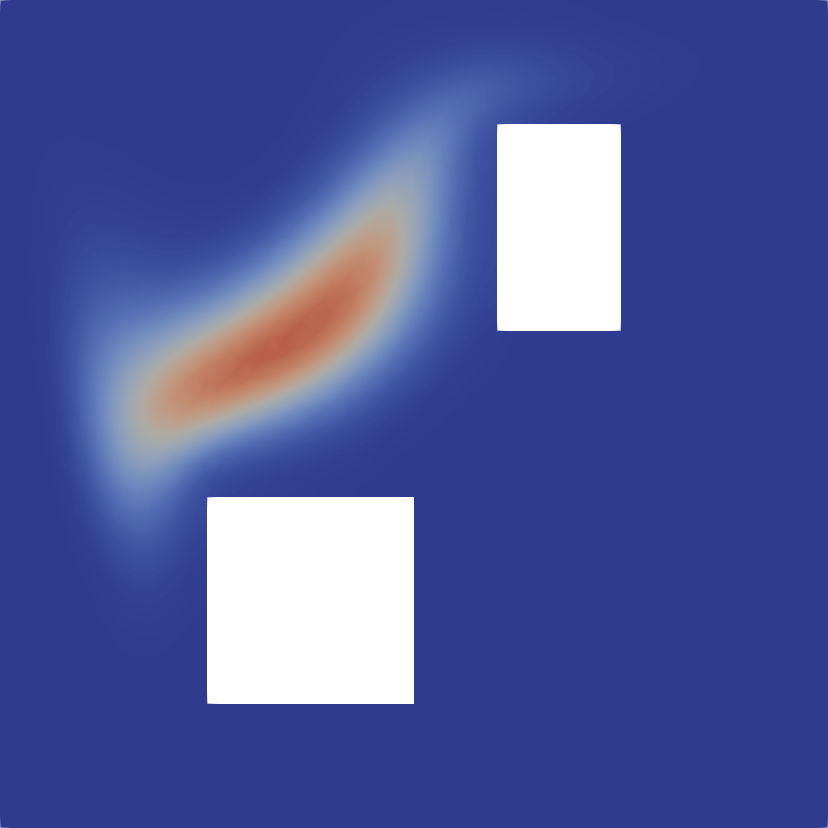}};\hspace{.15cm}
  \node at (0.05*\pos-0.25*\pos,  -0.2*\pos) {\sf b)};

  \node (3) at (0.6*\pos, 0*\pos){\includegraphics[height=.23\textwidth]{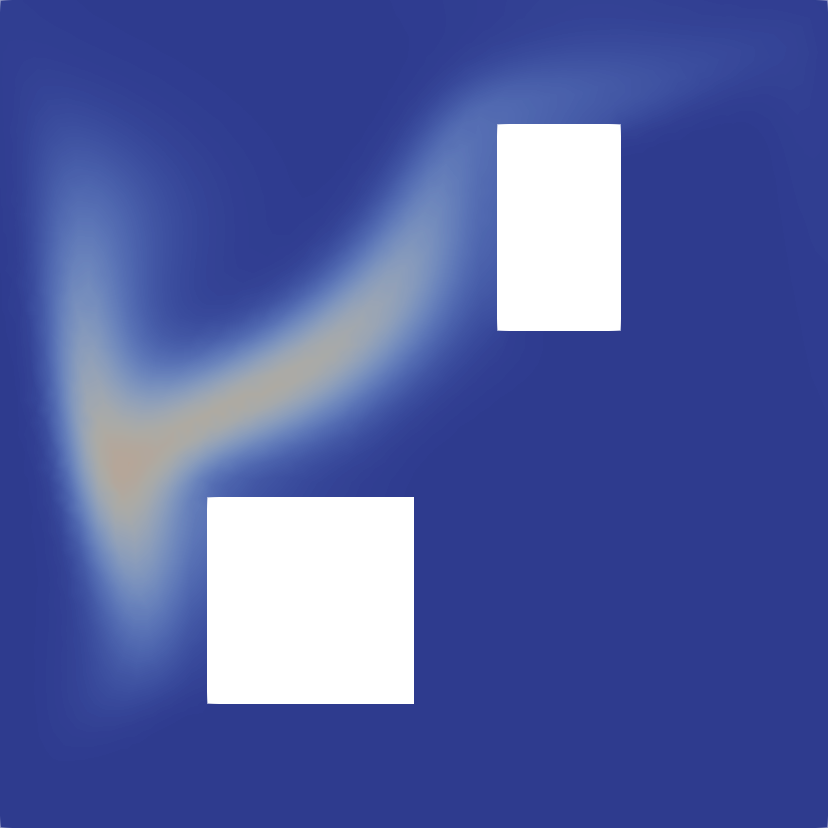}};\hspace{.15cm}
  \node at (0.6*\pos-0.25*\pos, -0.2*\pos) {\sf c)};

  \node (4) at (1.15*\pos, 0*\pos){\includegraphics[height=.23\textwidth]{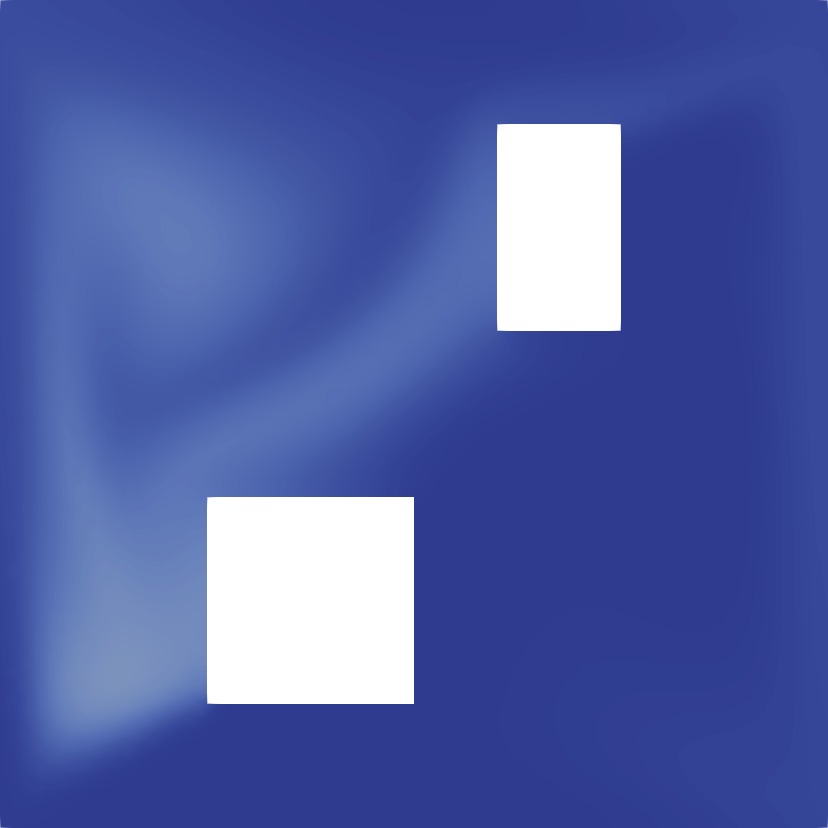}};
  \node at (1.15*\pos-0.22*\pos, -0.2*\pos) {\sf d)};

  \node (5) at (0*\pos-0.58*\pos, -0.57*\pos){\includegraphics[height=.23\textwidth]{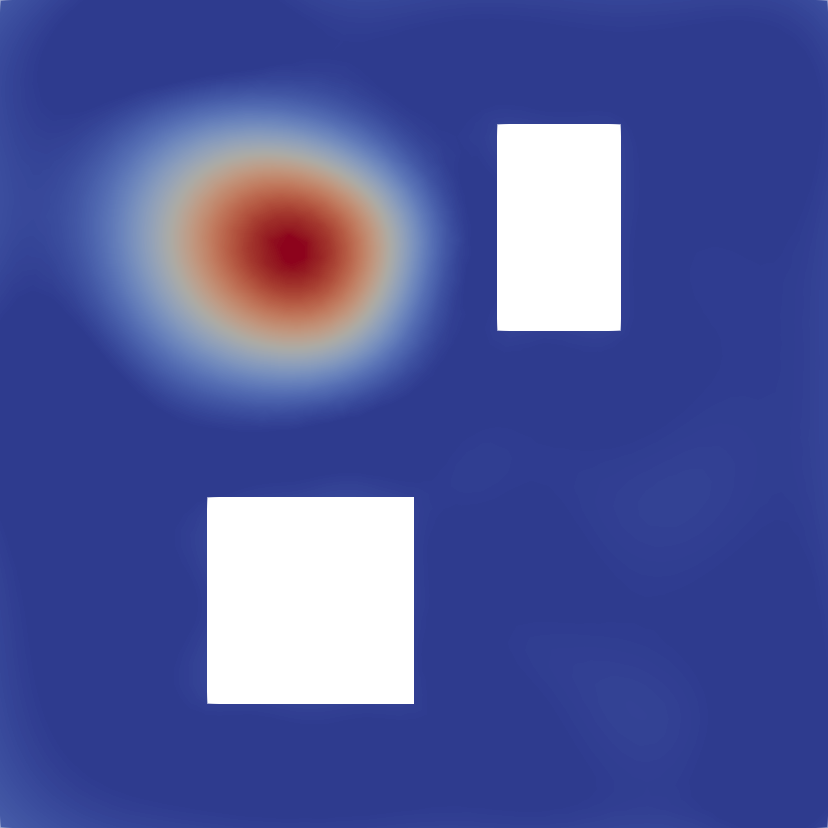}};\hspace{.15cm}
  \node at (-0.58*\pos-0.25*\pos, -0.77*\pos) {\sf e)};

  \node (6) at (-0.03*\pos, -0.57*\pos){\includegraphics[height=.23\textwidth]{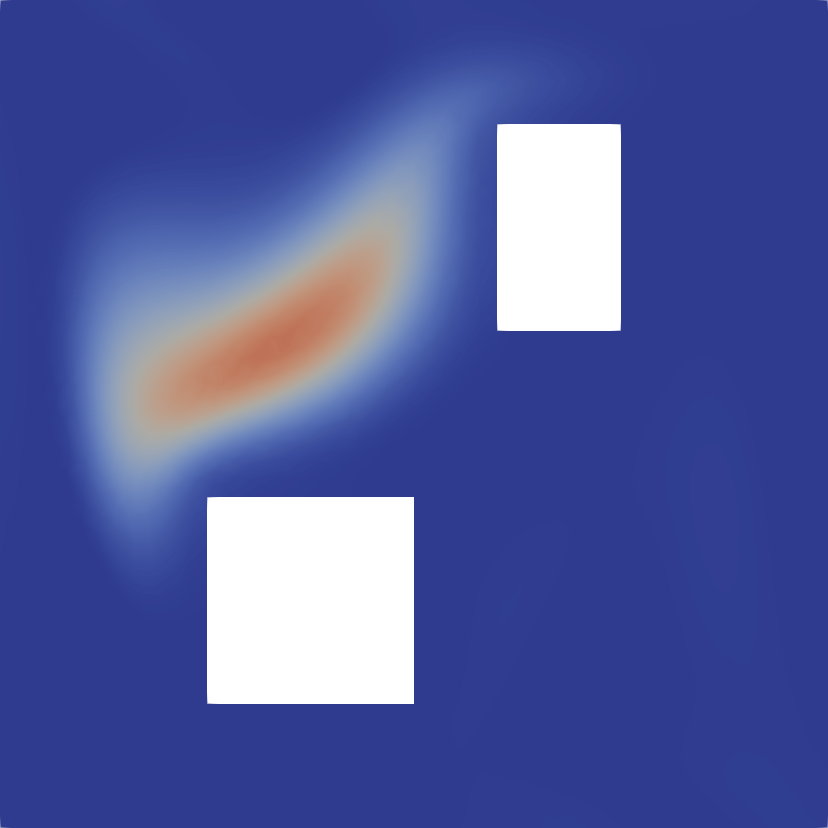}};\hspace{.15cm}
  \node at (-0.03*\pos-0.25*\pos,  -0.77*\pos) {\sf f)};

  \node (7) at (0.52*\pos, -0.57*\pos){\includegraphics[height=.23\textwidth]{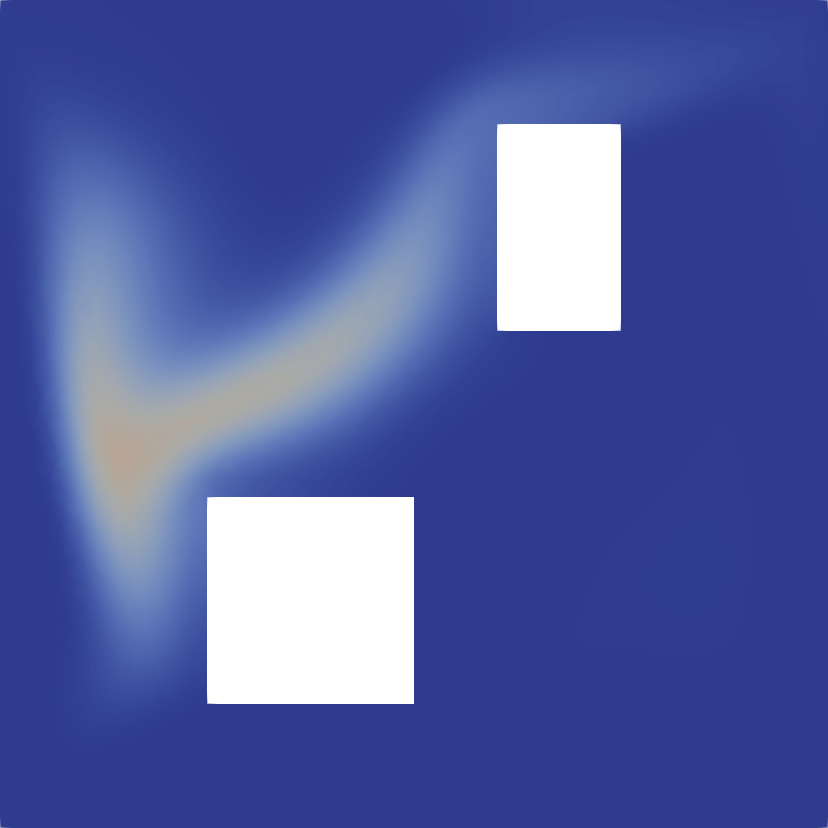}};\hspace{.15cm}
  \node at (0.52*\pos-0.25*\pos, -0.77*\pos) {\sf g)};

  \node (8) at (1.07*\pos, -0.57*\pos){\includegraphics[height=.23\textwidth]{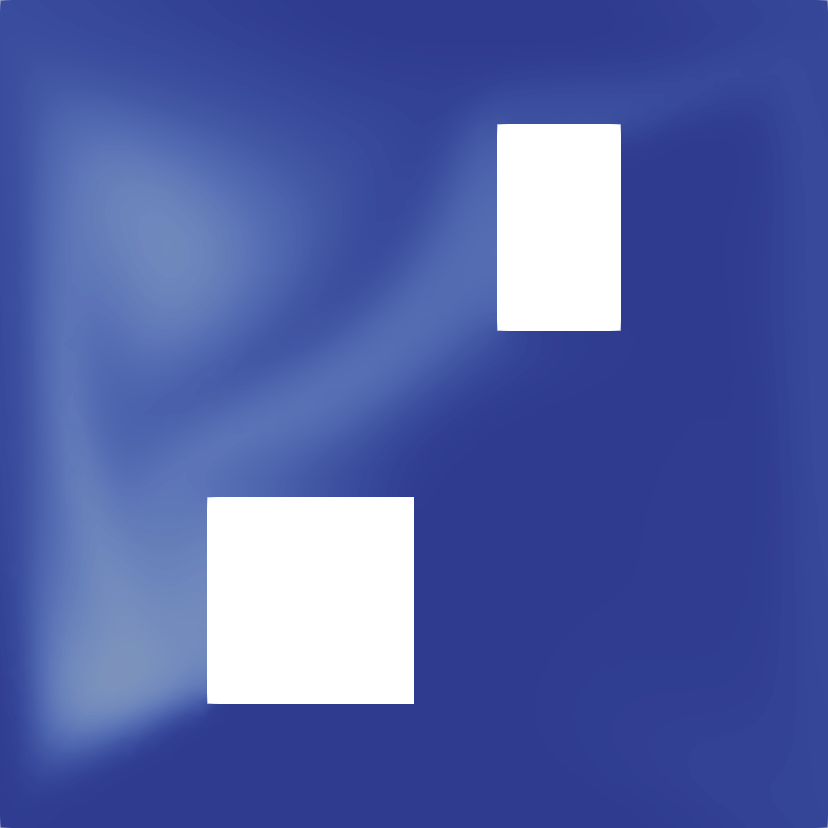}};
  \node at (1.07*\pos-0.22*\pos, -0.77*\pos) {\sf h)};
\end{tikzpicture}
\caption{Forward advective-diffusive transport estimate of the inverse
  solution at initial time t = 0 (a, d), at t = 1 (b, e), t = 2 (c, f)
  and at final time t = 4 (d) with the ``true'' (top) and MAP (bottom)
  as initial conditions.}
\label{fig:add-diff-state-attrueandmap}
\end{figure}
\end{center}

\begin{figure}[t]
  \centering
  \resizebox{\columnwidth}{!}{\input{tikz/spec-and-evecs-ad_diff.tex}}
  \caption{Left: Log-linear plot of the truncated spectrum of
    prior-preconditioned data misfit Hessian for observation times
    (sampled every 0.2 time units) in the intervals $[1, 4]$ (blue),
    $[2, 4]$ (red), and $[3, 4]$ (green). The low-rank approximation
    captures the dominant, data-informed portion of the spectrum. The
    eigenvalues are truncated at around $0.06$. Right:
    Prior-orthogonal eigenvectors of the prior-preconditioned data
    misfit Hessian corresponding (from left to right) to the 1st, 4th,
    7th, and 60th eigenvalues.  Eigenvectors corresponding to smaller
    eigenvalues are increasingly more oscillatory (and thus inform
    smaller length scales of the initial concentration) but are also
    increasingly less informed by the data.}
\label{fig:ad_diff-spectrum-evecs}
\end{figure}

\def \pos {0.43\columnwidth}
\begin{center}
\begin{figure}[tb]\centering
\begin{tikzpicture}
  \node (1) at (0*\pos-0.5*\pos, 0*\pos){\includegraphics[height=.23\textwidth]{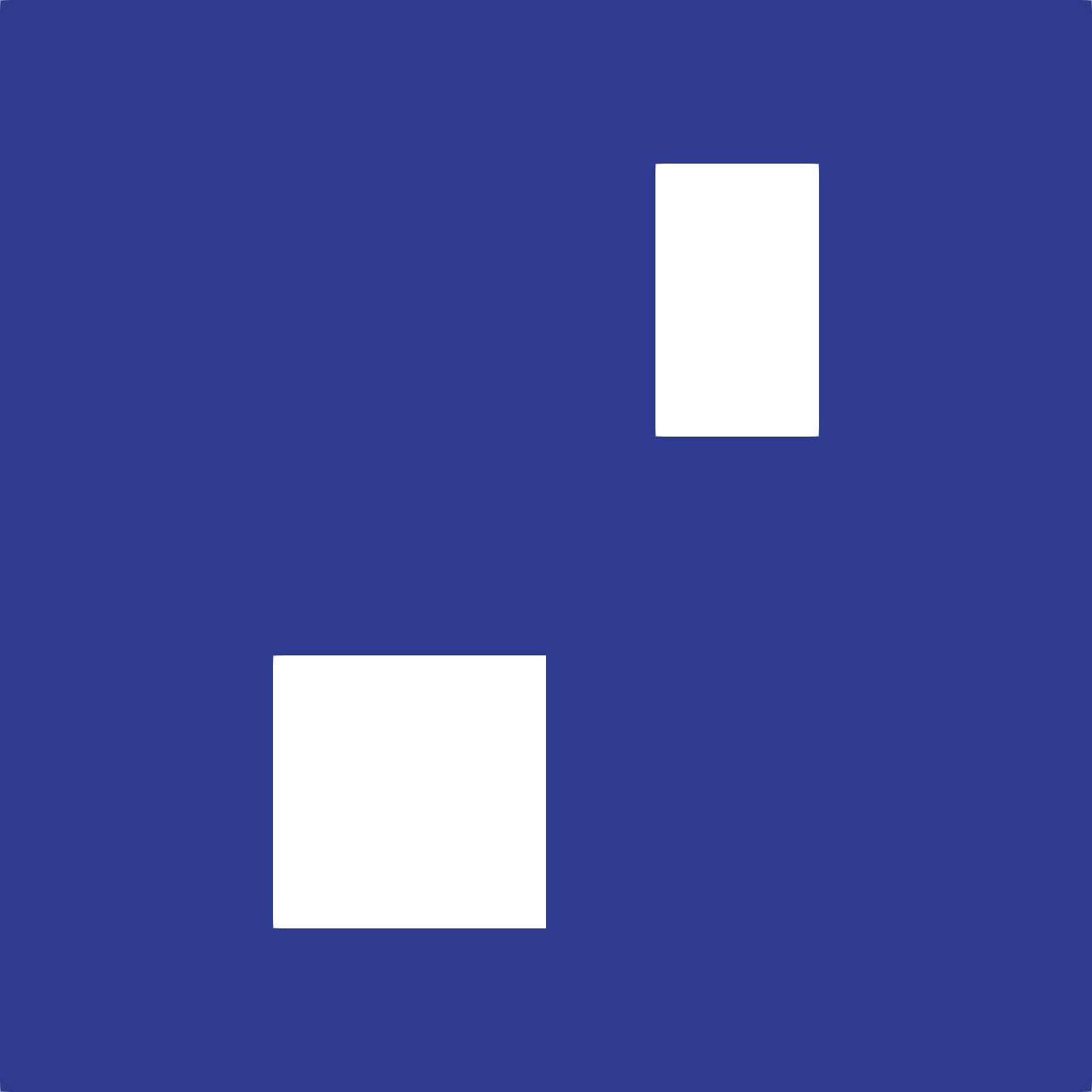}};\hspace{.15cm}
  \node at (-0.5*\pos-0.25*\pos, -0.2*\pos) {\sf a)};

  \node (2) at (0.05*\pos, 0*\pos){\includegraphics[height=.23\textwidth]{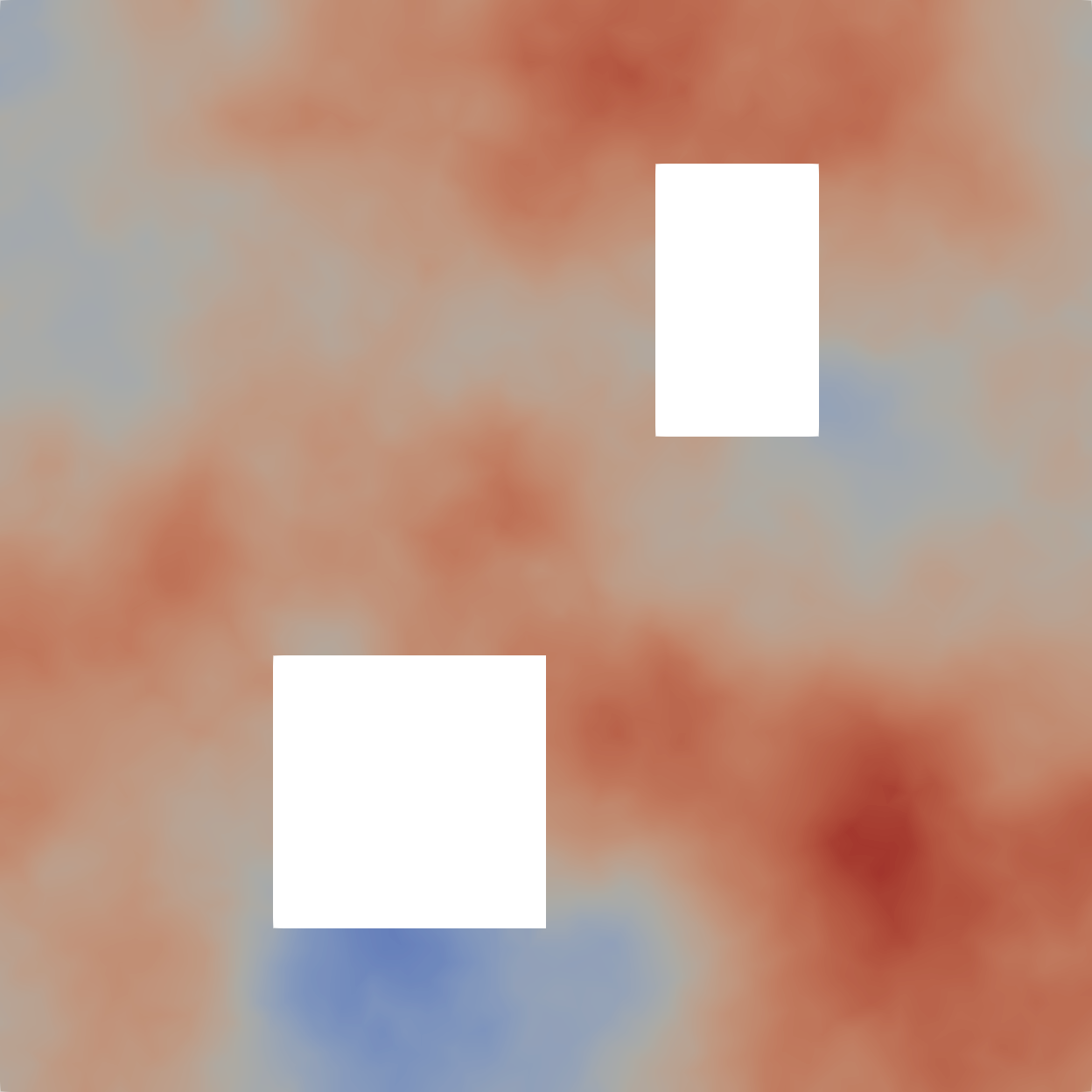}};\hspace{.15cm}
  \node at (0.05*\pos-0.25*\pos,  -0.2*\pos) {\sf b)};

  \node (3) at (0.6*\pos, 0*\pos){\includegraphics[height=.23\textwidth]{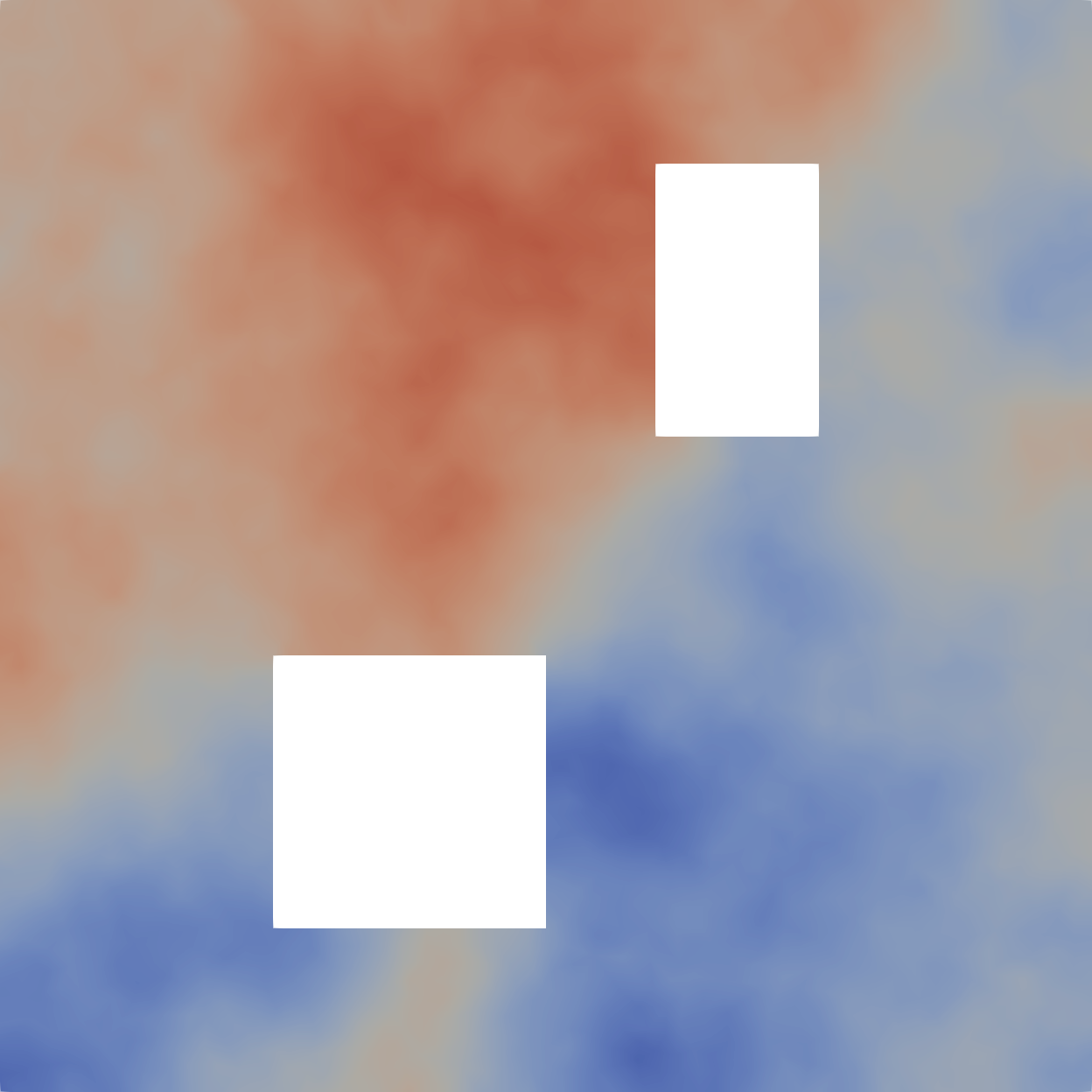}};\hspace{.15cm}
  \node at (0.6*\pos-0.25*\pos, -0.2*\pos) {\sf c)};

  \node (4) at (1.15*\pos, 0*\pos){\includegraphics[height=.23\textwidth]{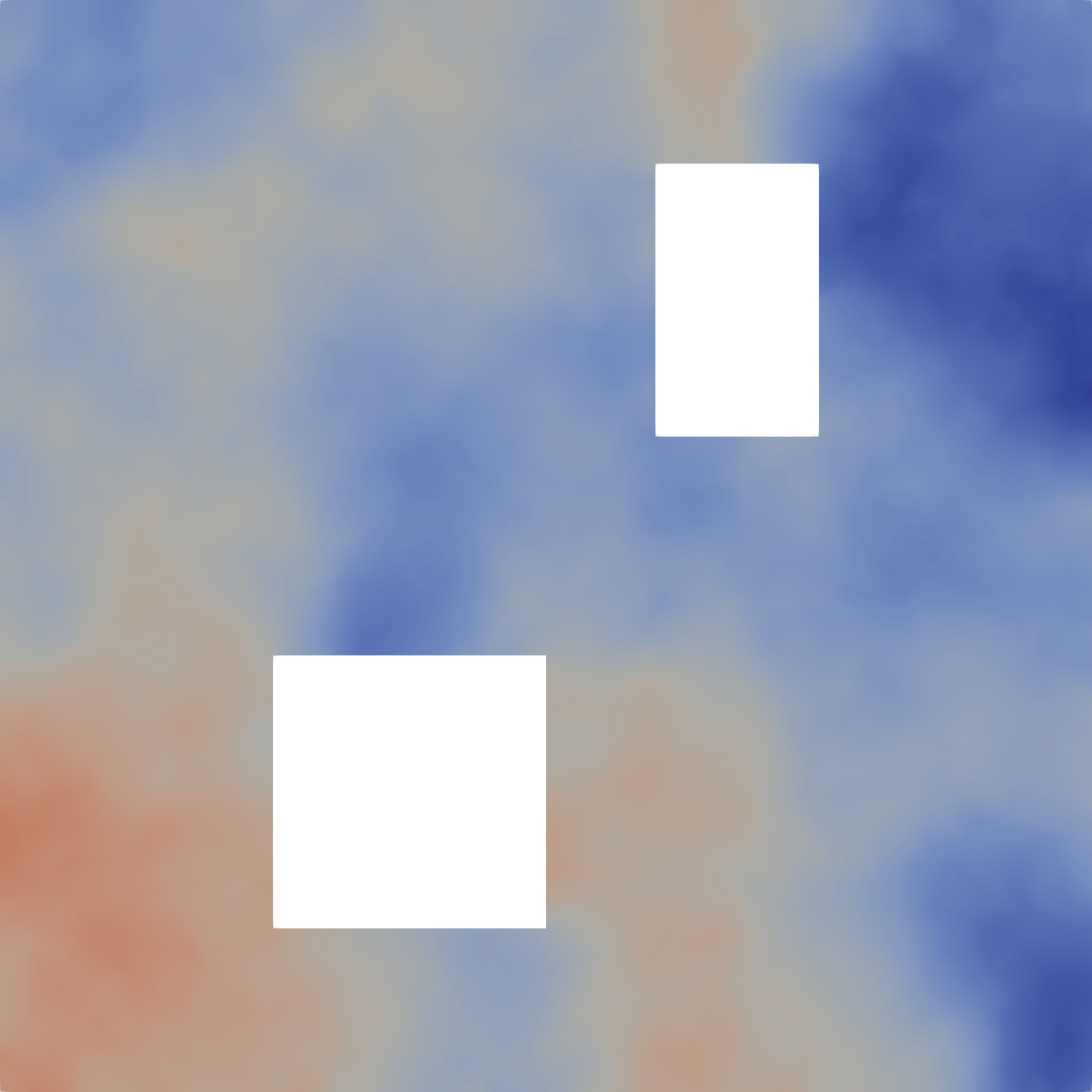}};
  \node at (1.15*\pos-0.22*\pos, -0.2*\pos) {\sf d)};

  \node (5) at (0*\pos-0.58*\pos, -0.57*\pos){\includegraphics[height=.23\textwidth]{extraplots/ad_diff-new/map}};\hspace{.15cm}
  \node at (-0.58*\pos-0.25*\pos, -0.77*\pos) {\sf e)};

  \node (6) at (-0.03*\pos, -0.57*\pos){\includegraphics[height=.23\textwidth]{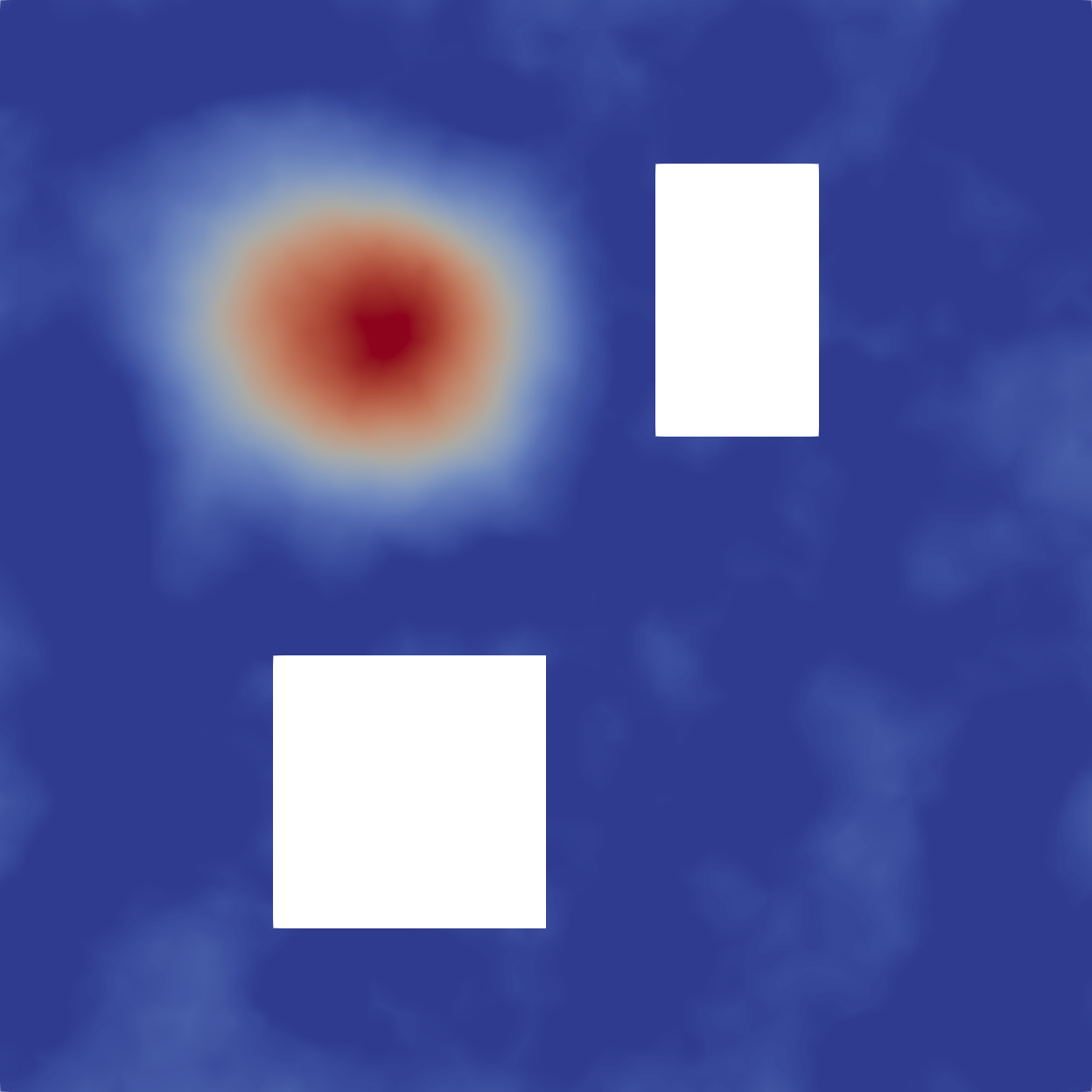}};\hspace{.15cm}
  \node at (-0.03*\pos-0.25*\pos,  -0.77*\pos) {\sf f)};

  \node (7) at (0.52*\pos, -0.57*\pos){\includegraphics[height=.23\textwidth]{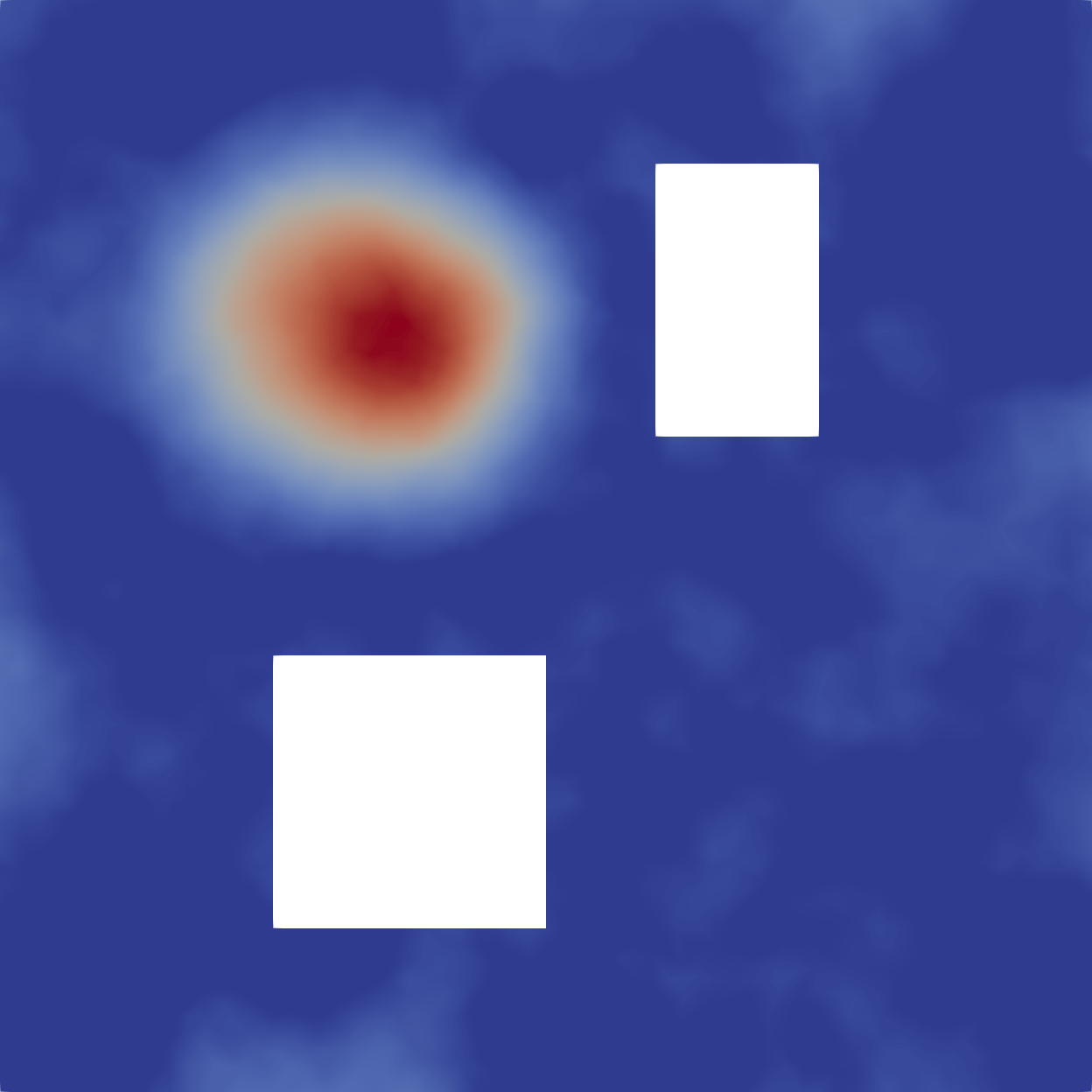}};\hspace{.15cm}
  \node at (0.52*\pos-0.25*\pos, -0.77*\pos) {\sf g)};

  \node (8) at (1.07*\pos, -0.57*\pos){\includegraphics[height=.23\textwidth]{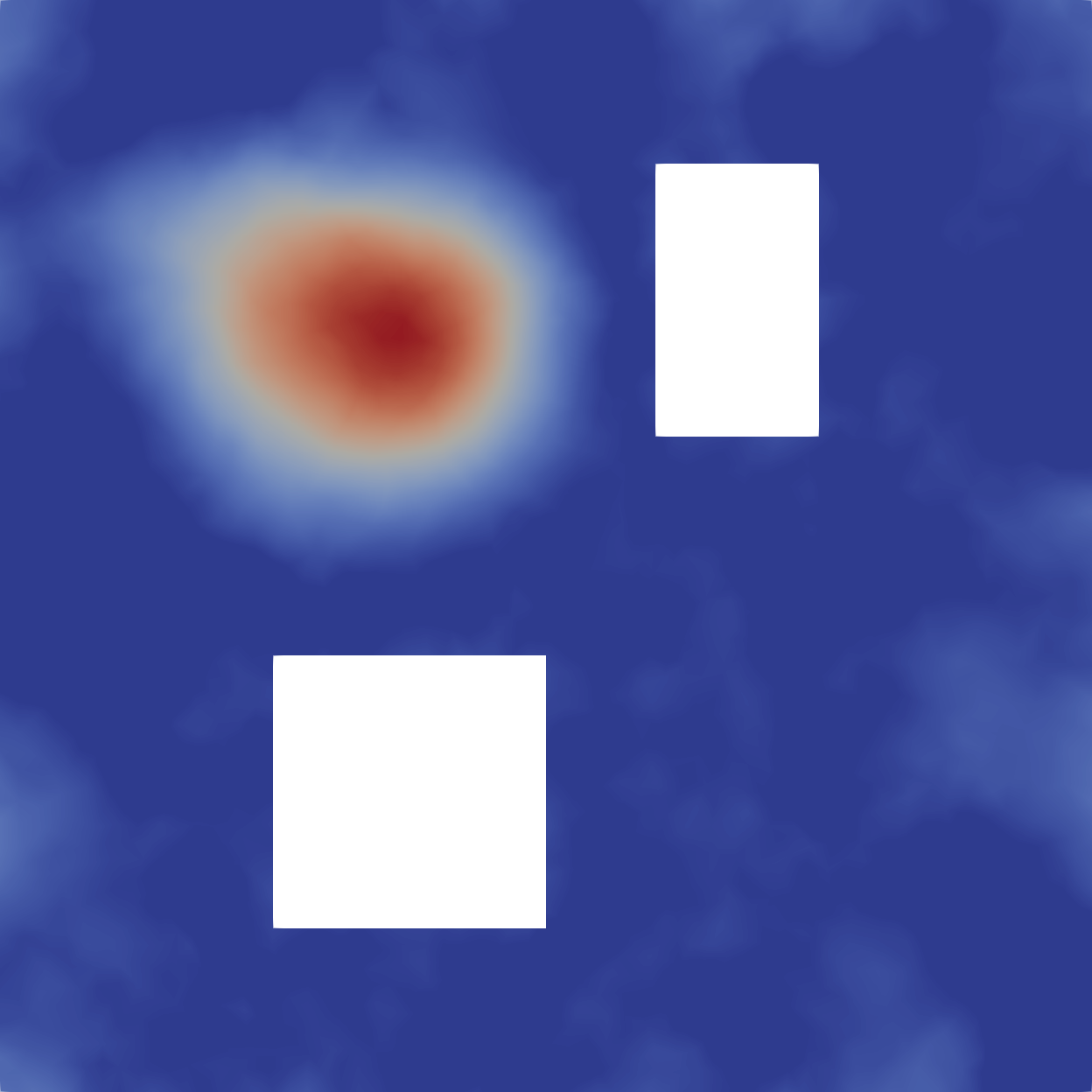}};
  \node at (1.07*\pos-0.22*\pos, -0.77*\pos) {\sf h)};
  
\end{tikzpicture}
\caption{Top: Prior mean initial concentration $\iparpr$ (a), and
  samples drawn from the prior distribution (b)--(d). Bottom: The MAP
  point (e), and samples drawn from the posterior
  distribution (f)--(h).}
\label{fig:add-diff-prior-post}
\end{figure}
\end{center}

\section{Conclusions}\label{sec:conclusions}
We have presented an extensible software framework
for large-scale
deterministic and linearized Bayesian inverse problems governed by
partial differential equations. The main advantage of this framework
is that it exploits the structure of the underlying
infinite-dimensional PDE based parameter-to-observable map, in
particular the low effective dimensionality, which leads to scalable
algorithms for carrying out the solution of deterministic and
linearized Bayesian inverse problems. By scalable, we mean that the
cost---measured in number of (linearized) forward (and adjoint)
solves---is independent of the state, parameter, and data
dimensions. The cost depends only on the number of modes in parameter
space that are informed by the data.
We have described the main algorithms implemented in \hip, namely the
inexact Newton-CG method to compute the MAP point
(Section~\ref{subsec:NCG}), randomized eigensolvers to compute the low
rank approximation of the Hessian evaluated at the MAP point (Section~\ref{subsec:lra}), and
algorithms for sampling and computing the pointwise variance from large-scale
Gaussian random fields (Sections~\ref{subsec:sampling}
and~\ref{subsec:var}). To illustrate their use, we applied these methods to two model problems:
inversion for the log coefficient field in a Poisson equation, and inversion for the initial condition in a
time-dependent advection-diffusion equation.

The contributions of our work are as follows. On the algorithm side,
our framework incorporates modifications of state-of-the-art
algorithms to ensure consistency with infinite-dimensional
settings and a novel square-root-free implementation of the low-rank
approximation of the Hessian, sampling strategies, and pointwise variance field computation. On the software
side, we created a library for the solution of deterministic and linearized Bayesian inverse problems that allows researchers who are familiar
with variational methods to solve inverse problems under uncertainty
even without possessing expertise in all of the necessary numerical
optimization and statistical aspects.  Our framework provides dimension-independent algorithms for
finding the maximum a posteriori (MAP) point, constructing a low-rank based
approximation of the Hessian and its inverse at the MAP, sampling from the prior and
posterior distributions, and computing pointwise variance fields. \hip is easily extensible, that is, if a user can
express the forward problem in variational form using \Fe, \hip effortlessly allows solving the inverse problem, exploring and testing
various priors, observation operators, noise covariance models, etc.

The framework presented here relies
on a second order Taylor expansion of the negative log-likelihood with respect to the uncertain parameter centered 
 at the MAP point, which leads to the Laplace approximation of the posterior distribution.
Ultimately one
would like to relax this approximation and fully explore the resulting
non-Gaussian distributions.  Ongoing work includes the implementation of
scalable, robust, Hessian-based MCMC methods capitalizing on \hip's
capabilities to build local \emph{Laplace} approximations of the posterior
based on gradient and Hessian information as described here.

\begin{acks}
The authors would like to thank Georg Stadler for providing very helpful comments on an
earlier version of this paper.

This work was supported by the U.S. National Science Foundation, Software Infrastructure for
Sustained Innovation (SI2: SSE \& SSI) Program under grants
ACI-1550593, and ACI-1550547; the Defense Advanced Research Projects
Agency, Enabling Quantification of Uncertainty in Physical Systems
  Program under grant W911NF-15-2-0121; the US National Science
  Foundation, Division Of Chemical, Bioengineering, Environmental, \&
  Transport Systems under grant CBET-1508713; and the Air Force Office
  of Scientific Research, Computational Mathematics program under
  grant FA9550-12-1-81243.
\end{acks}

\bibliographystyle{ACM-Reference-Format}
\bibliography{ccgo,local,ymarz}

\appendix

\section{Gradient and Hessian actions computation for the inverse problem governed by the Poisson PDE}
\label{app-sec:MAP}
In what follows, we apply the technique outlined in
Section~\ref{sec:det} and derive expressions for the gradient and
Hessian actions of the cost functional $\mathcal{J}(\ipar)$ defined
in~\eqref{eq:cost-elliptic}. The Lagrangian functional for this
optimization problem is given by
\begin{equation}\label{eq:model:L}
  \LI^\mathcal{G}(u,m,p):= \J(\ipar) + \ip{\Exp{m}\grad u}{\grad p}
  - \ip{f}{p} - \ip{p}{h}_{\GN},
\end{equation}
where the last three terms stem from the variational form
\eqref{eq:poisson_weak} of the forward problem \eqref{equ:poi}.
The formal Lagrange multiplier method
\cite{Troltzsch10} requires that, at a minimizer of
\eqref{eq:cost-elliptic}, variations of the Lagrangian functional
with respect to $p$ and $u$ vanish, which yields to solving the forward and adjoint problems
\begin{subequations}
  \begin{align}
    \ip{\Exp{m} \grad u}{\grad \ut{p}} -
    \ip{f}{\ut{p}} - \ip{\ut{p}}{h}_{\GN} & = 0, \quad \forall \ut{p} \in \mathcal{V}_0;  \label{eq:firststate} \\
    \ip{\Exp{m} \grad \ut u}{\grad p} 
    +\ip{\B^* \ncov^{-1}(\B u - \obs)}{\ut{u}} &= 0,  \quad \forall \ut{u} \in \mathcal{V}_0.  \label{eq:firstadj}
  \end{align}
\end{subequations}
The strong form of the forward problem is given in \eqref{equ:poi}, while the strong form of
the adjoint problem reads
\begin{equation}\label{equ:poi-adj}
  \begin{split}
    -\grad \cdot (\Exp{m} \grad p) &= \B^* \ncov^{-1}(\B u - \obs) \quad \text{ in }\D, \\
    p  &= 0 \quad \text{ on } \GD, \\
    \Exp{m} \grad{p} \cdot \vec{n} &= 0 \quad \text{ on } \GN.
  \end{split}
\end{equation}
Finally, the gradient of the cost functional~\eqref{eq:cost-elliptic}
is given in weak form by
\begin{align}
    \left( \mathcal{G}(m), \tilde m\right) = \ip{m - \iparpr}{\ut{m}}_{\Cprior^{-1}} + \ip{\ut{m} \Exp{m}\grad
      u}{\grad p}, \quad \forall \ut{m} \in \iparspace, \label{eq:firstcontrol}
\end{align}
where $u$ and $p$ are solutions to the forward and adjoint problems \eqref{eq:firststate}-\eqref{eq:firstadj},
respectively \cite{Troltzsch10,BorziSchulz12}. In strong form this reads
\begin{align}
    \mathcal{G}(\ipar)  = 
\left\{
\begin{array}{ll}    
    \Cprior^{-1} (\ipar - \iparpr) + \Exp{m}\left({\grad
      u}\cdot{\grad p}\right) & \text{in }\D,\\
      \gamma \left(\Theta \nabla \ipar\right) \cdot\boldsymbol{n} + \beta \ipar & \text{on } \partial \D.
\end{array}
\right.
\label{eq:firstcontrolstrongform}
\end{align}
We note that to evaluate the gradient for a given parameter $\ipar$,
one needs to solve the forward problem~\eqref{eq:firststate} for $u$,
and then given $\ipar$ and $u$ solve the adjoint problem for $p$. This
evaluation of the gradient costs one forward and one adjoint PDE
solve.

Next, we derive the expression of the Hessian action following
Section~\ref{sec:det}.
The second order Lagrangian functional in this case reads
\begin{align*}
\LI^\mathcal{H}(u,m,p; \hat{u}, \hat{m}, \hat{p}) &:=
\left( \mathcal{G}(m), \hat{m} \right) \\ {}& +\ip{\Exp{m} \grad u}{\grad \hat{p}} -
\ip{f}{\hat{p}} - \ip{\hat{p}}{h}_{\GN} \\ {}&+ \ip{\Exp{m} \grad \hat
  u}{\grad p} + \ip{\B^* \ncov^{-1}(\B u - \obs)}{\hat{u}}.
\end{align*}
To obtain the action of the Hessian in a direction $\hat m$ we take
the variation of $\LI^\mathcal{H}$ with respect to $m$, namely
\begin{align}
  \! \left( \tilde m, \H(m) \hat m\right)
  & = \ip{\ut{m} \Exp{m} \grad \hat u}{\grad p}
  \!+ \! \ip{\hat m}{\ut{m}}_{\Cprior^{-1}} \! \\
  {} & + \!\ip{\ut{m} \hat m \Exp{m} \grad u}{\grad p}
  \!+ \! \ip{ \ut{m} \Exp{m} \grad u}{\grad \hat p},  \quad \forall \ut{m} \in \iparspace, \label{eq:incrementals3}
\end{align}
where as before $u$ and $p$ are the solutions of the forward and
adjoint problems in \eqref{eq:firststate} and \eqref{eq:firstadj}, respectively, and the $\hat u$ and $\hat p$ are the solutions of the
incremental forward and adjoint problems, respectively. These equations are given by
\begin{equation}\label{eq:incrementals1}
    \ip{ \Exp{m} \grad \hat u}{\grad \ut p} + \ip{ \hat m \Exp{m}\grad u}{\grad
      \ut p} = 0, \quad \forall \ut p \in \V, 
\end{equation}
and
\begin{equation}\label{eq:incrementals2}
    \ip{ \B^* \ncov^{-1} \B \hat u}{\ut u} + \ip{ \hat m \Exp{m}\grad \ut u}{\grad p}
    + \ip{\Exp{m} \grad \ut u}{\grad \hat p} = 0, \quad \forall \ut u \in \V.
\end{equation}

Once we have the gradient and Hessian action expressions, we can apply
Algorithm~\ref{alg:InexactNewtonCG} to solve the optimization problem
given by~\eqref{eq:cost-elliptic}.

\section{Gradient computation for inverse problem governed by the advection-diffusion PDE}
\label{app:ad-map}

To derive an expression for the gradient of $\mathcal{J}(\ipar)$ in
\eqref{eq:ad_cost}, we define the Lagrangian functional
\begin{align*}
\mathscr{L}^\mathcal{G}(u,\ipar,p,p_0):=\mathcal{J}(\ipar) &+\int_0^T\!\int_\D (u_t +
\vec{v}\cdot\nabla u)p\,d\x\,dt\\ & +
\int_0^T\!\int_\D\kappa\nabla u \cdot \nabla p\,d\x\,dt + \int_\D
(u(\vec{x}, 0) - \ipar )p_0\,d\x,
\end{align*}
where $p \in \mathcal{V}$ and $p_0 \in \iparspace$ are the Lagrangian
multiplier, i.e., the adjoint variables, for, respectively, the
advection-diffusion PDE and initial condition in the forward problem
\eqref{eq:ad}.  Expressions needed to compute the gradient of
\eqref{eq:ad_cost} are obtained by setting variations of the
Lagrangian $\mathscr{L}^\mathcal{G}$ with respect to $p$, $p_0$ and
$u$ to zero. Variations with respect to $p$ and $p_0$ recover the
variational form \eqref{eq:ad_weak} of the forward problem.  The
variation with respect to $u$ in an arbitrary direction $\tilde u$
yields
\begin{align*}
\frac{1}{\sigma^2} \sum_{i=1}^{n_s}\int_{T_1}^T(\mathcal{B}u-\boldsymbol{d}_i)\tilde u \,\delta_{t_i} \,dt +
\int_0^T\int_\D (\tilde u_t + \vec{v}\cdot\nabla \tilde u)p\,d\x\,dt\\
+ \int_0^T\int_\D\kappa\nabla \tilde u \cdot \nabla
p\,d\x\,dt + \int_\D \tilde u(\vec{x}, 0)p_0\,d\x = 0 & \quad \forall \tilde u \in \mathcal{V}.
\end{align*}
Integration by parts in time for the term $\tilde u_t p$ and in space
for $(\vec{v}\cdot\nabla\tilde u) p$ and $\kappa\nabla \tilde u \cdot
\nabla p$ results in
\begin{alignat*}{2}
\frac{1}{\sigma^2} \sum_{i=1}^{n_s}\int_{T_1}^T(\mathcal{B}u & -\boldsymbol{d}_i)\tilde u \, \delta_{t_i} \,dt -
\int_0^T\!\int_\D (\tilde p_t + \nabla
\cdot(\vec{v}p)+\kappa\Delta p)\tilde u\,d\x\,dt + \int_\D \tilde u(\vec{x}, T)p(\vec{x}, T) \\
&- \tilde u(\vec{x}, 0)p(\vec{x}, 0) + \tilde u(\vec{x}, 0)p_0 \,d\x +
\int_0^T\int_{\partial\D}(\vec{v}p + \kappa\nabla p)\cdot \vec{n}\tilde
u \, d\x\,dt = 0, 
\end{alignat*}
$\tilde u \in \mathcal{V}$. This implies $p_0=p(\vec{x}, 0)$ and leads
to the strong form of the adjoint problem,
\begin{equation}\label{eq:ad:adj}
  \begin{aligned}
    -p_t - \nabla \cdot (p \vec{v}) - \kappa\Delta p  &=  -\frac{1}{\sigma^2}\mathcal{B}^*\sum_{i=1}^{n_s}\left(\mathcal{B}u - \boldsymbol{d}_i\right)\delta_{t_i}
&\quad&\text{ in
    }\D\times (0,T),\\
    p(\cdot, T) &= 0 &&\text{ in } \D,  \\
    (\vec{v}p+\kappa\nabla p)\cdot \vec{n} &=  0 &&\text{ on }
   \partial\D  \times (0,T).
  \end{aligned}
\end{equation}
Note that \eqref{eq:ad:adj} is a final value problem, since $p$ is
specified at $t=T$ rather than at $t=0$. Thus, \eqref{eq:ad:adj} is solved backwards in time, which amounts to the solution of an
advection-diffusion equation with velocity $-\vec{v}$.

Finally, the variation of the Lagrangian with respect to the initial condition
$\ipar$ in a direction $\tilde{\ipar}$ gives the weak form of the gradient of the
cost functional $\mathcal{J}(\ipar)$,
\begin{equation}\label{eq:ad:control_weak}
\left( \mathcal{G}(\ipar), \tilde \ipar\right) =
\int_\D \big(\mathcal{A}(\ipar -  \iparpr) \big) \big(\mathcal{A}\tilde \ipar\big) - p(\vec{x}, 0)\tilde \ipar\,d\x,  \quad \forall \tilde \ipar \in \iparspace,
\end{equation}
where we used $p_0 = p(\vec{x}, 0)$. The strong form of the gradient
expression then reads
\begin{equation}\label{eq:ad:control}
\mathcal{G}(\ipar) = \left\{
\begin{array}{rl}
\mathcal{A}^2(\ipar -  \iparpr) - p(\vec{x}, 0) & \text{ in } \D\\
\gamma \nabla m \cdot\boldsymbol{n} + \beta m & \text{ in } \partial\D.
\end{array}
\right.
\end{equation}
Thus, the adjoint of the parameter-to-observable map $\iFF^*$ is
defined by setting $\iFF^*\obs = p(\vec{x}, 0)$.  We note that the
gradient expression in \eqref{eq:ad:control} is linear in $\ipar$,
since $\iadj$ depends linearly on $\istate$ through the solution of
the adjoint problem \eqref{eq:ad:adj}, and $\istate$ depends linearly
on $m$ through the solution of the forward problem \eqref{eq:ad}.
Elimination of the forward and adjoint equations for a given $\ipar$
gives the action of the linear operator in \eqref{eq:ad:control} in
the direction of that $\ipar$. Thus, to compute the MAP point, we use
CG to set $\mathcal{G}(\ipar) = 0$ with the
inverse of the regularization operator $\mathcal{A}^2$ as
preconditioner.

\section{Finite element assembly and rectangular decompositions}\label{sec:rect_decomp}

In this section, we describe a technique based on rectangular
decompositions of finite element matrices to efficiently generate large
scale samples from priors of the form described in
Section~\ref{sec:bayes}. A similar approach has been also proposed in \citep{CrociGilesRognesEtAl18} to 
generate samples of white noise.
The method described here is more general  as it can be applied to matrices stemming from
finite element discretization of any differential operator and not only mass matrices.
More specifically, we present a finite
element assembly procedure to compute a rectangular decomposition of
the form
\begin{equation}\label{eq:rect_decomp}
 \mat{A} = \mat{C} \mat{C}^T,
\end{equation}
for any symmetric positive definite finite element matrix $\mat{A}$. 
More specifically, consider the finite element assembly procedure for a generic symmetric positive definite bilinear form $a(u_h,v_h)$ on $\iparspace_h$, a finite-dimensional subspace of $\iparspace \subseteq L^2(\mathcal{D})$.
The entry $(i,j)$ of the matrix $\mat{A}$, which stems from finite element discretization of the bilinear form $a(u_h,v_h)$, are given by
$$ A_{i,j} = a\left( \phi_i, \phi_j \right), \quad i,j = 1, \ldots, n, $$
where $\{ \phi_i \}_{i=1}^n$ is the finite element basis of the space $\iparspace_h$.
In the finite element assembly procedure we first compute the element matrices $\mat{A}_e$, which correspond to the restriction of the bilinear form $a$ to each element $e$ in the mesh. Then, using the global-to-local mapping of the degrees of freedom (dof) $\mat{G}_e$, the global matrix $\mat{A}$ is computed by summing all of the local contributions as follows:
\begin{equation}
\mat{A} = \sum_e \mat{G}_e^T \mat{A}_e \mat{G}_e = \sum_e  \mat{G}_e^T \mat{B}^T \mat{D}_e \mat{B} \mat{G}_e.
\end{equation}
Here we have written the element matrix $\mat{A}_e = \mat{B}^T \mat{D}_e \mat{B}$ as the product of the element-independent dof-to-quadrature point basis evaluation matrix $\mat{B}$ and the (block) diagonal matrix $\mat{D}_e \in \mathbb{R}^{q \times q}$ at the quadrature points, where $q$ denotes the total number of quadrature nodes over all elements, which scales linearly with the number of elements in the mesh.\\
A rectangular decomposition of $\mat{A}$ can then be explicitly constructed from the matrices $\mat{G}_e$, $\mat{B}$, and $\mat{D}_e$ as follows.
For each element $e$ of the mesh we define the matrix $\mat{C}_e = \mat{G}_e^T \mat{B}^T \mat{D}_e^{\frac{1}{2}}$. Since, for any two elements $e_i$ and $e_j$ ( $i\neq j$) in the mesh, the sets of quadrature nodes relative to the elements $e_i$ and $e_j$ are disjoint, we have that
\begin{equation}\label{eq:Ce_ortho}
 \mat{C}_{e_i} \mat{C}^T_{e_j} = \delta_{ij} \mat{A}_{e_i}.
\end{equation}
Then the rectangular matrix $\mat{C} \in \mathbb{R}^{n \times q}$ defined as
\begin{equation}
 \mat{C} = \sum_e \mat{G}_e^T \mat{B}^T \mat{D}_e^{\frac{1}{2}} = \sum_e \mat{C}_e
\end{equation}
satisfies \eqref{eq:rect_decomp}. In fact, thanks to \eqref{eq:Ce_ortho} we have
$$ \mat{C}\mat{C}^T = \left( \sum_e \mat{C}_e \right)\left( \sum_e \mat{C}_e \right)^T = \sum_e \mat{C}_e \mat{C}_e^T = \sum_e \mat{G}_e^T \mat{B}^T \mat{D}_e \mat{B} \mat{G}_e = \mat{A}. $$

\end{document}

%% file: defjcp14.txt
\usepackage{booktabs}
\usepackage[mathscr]{euscript}
\usepackage{color}
\usepackage{graphicx}
\usepackage{algorithmic,algorithm}
\usepackage{tikz,pgfplots}
\usepackage{upgreek}
\usepackage{multirow}
\usepackage{yfonts}
\usepackage{mathtools}
\usepackage[normalem]{ulem}
\usepackage{latexsym}
\usepackage{url}
\usepackage{amsmath,amssymb,amsbsy,amsfonts}
\usepackage{enumitem}

\definecolor{utorange}{rgb}{0.8,0.33,0.}
\definecolor{themec}{RGB}{51,108,121}
\definecolor{darkred}{rgb}{.6,.1,.1}
\definecolor{darkblue}{rgb}{.1,.1,.9}
\definecolor{greenback}{rgb}{.19,.94,.13}
\definecolor{orange}{rgb}{.76,.39,.13}
\definecolor{grass}{rgb}{.19,.64,.13}
\definecolor{sierp}{RGB}{209,28,209}
\definecolor{bgorange}{rgb}{1.,.95,.78}
\definecolor{grassgreen}{RGB}{92,135,39}
\definecolor{thinbox}{rgb}{.7,.8,1.}

\renewcommand{\vec}[1]{{\mathchoice
                     {\mbox{\boldmath$\displaystyle{#1}$}}
                     {\mbox{\boldmath$\textstyle{#1}$}}
                     {\mbox{\boldmath$\scriptstyle{#1}$}}
                     {\mbox{\boldmath$\scriptscriptstyle{#1}$}}}}

\newcommand{\ip}[2]{{\left\langle {#1}, {#2} \right\rangle}}
\newcommand{\mip}[2]{\left\langle{#1}, {#2}\right\rangle_{\!\scriptscriptstyle{\text{M}}}}
\newcommand{\mat}[1]{\mathbf{{#1}}}
\newcommand\restr[2]{{
  \left.\kern-\nulldelimiterspace 
  {#1}\vphantom{\big|} \right|_{#2}}}

\DeclareMathOperator*{\argmin}{argmin} 
\DeclareMathOperator{\diag}{diag} 

\newcommand{\R}{\mathbb{R}}

\newcommand{\BB}{\mat{B}}

\newcommand{\RR}{\mat{R}}

\newcommand{\B}{\mathcal{B}}

\newcommand{\D}{\mathcal{D}}

\newcommand{\G}{\mathcal{G}}
\newcommand{\J}{\mathcal{J}}

\newcommand{\bit}{\begin{itemize}}
\newcommand{\eit}{\end{itemize}}
\newcommand{\bdm} {\begin{displaymath}}
\newcommand{\edm} {\end{displaymath}}

\newcommand{\ave}{\mathsf{E}}
\newcommand{\GM}[2]{\mathcal{N}\!\left( {#1}, {#2}\right)}

\newcommand{\iparspace}{\mathcal{M}}

\newcommand{\obs}{\vec{d}}
\newcommand{\ipar}{m}
\newcommand{\istate}{u}
\newcommand{\iadj}{p}
\newcommand{\iparpr}{m_{\text{pr}}}
\newcommand{\iparmap}{m_{\scriptscriptstyle\text{MAP}}}

\newcommand{\M}{\mat{M}}
\newcommand{\ff}{\mathcal{F}} 

\newcommand{\FF}{{\ensuremath{\mat{F}}}}

\newcommand{\like}{\pi_{\mbox{\tiny like}}}

\newcommand{\post}{\pi_{\mbox{\tiny post}}}

\newcommand{\mupost}{\mu_{\mbox{\tiny post}}}
\newcommand{\muprior}{\mu_{\mbox{\tiny prior}}}
\newcommand{\mulaplace}{\hat{\mu}_{\mbox{\tiny post}}}

\newcommand{\ncov} {\bs{\Gamma}_{\!\mbox{\tiny noise}} }
\newcommand{\prcov} {\bs{\Gamma}_{\!\mbox{\tiny prior}} }
\newcommand{\postcov} {\bs{\Gamma}_{\!\mbox{\tiny post}} }

\newcommand{\Cprior}{\mathcal{C}_{\text{\tiny{prior}}}}
\newcommand{\Cpost}{\mathcal{C}_{\text{\tiny{post}}}}

\newcommand{\map} {{m}_{\mbox{\tiny MAP}} }

\newcommand{\mpr} {{\vec{m}}_{\mbox{\tiny pr}} }

\newcommand{\Hmisfit}{ \mat{H}_{\mbox{\tiny misfit}} }
\newcommand{\tHmisfit}{\tilde{\mat{H}}_{\mbox{\tiny misfit}} }
\newcommand{\tH}{\tilde{\mat{H}} }

\newcommand{\adjMacroMM}[1]{{#1}^*}

\newcommand{\Badj}{\adjMacroMM{\BB}}

\newcommand{\mc}[1]{\mathcal{#1}}

\newcommand{\gbf}[1]{\boldsymbol{#1}}
\newcommand{\bs}[1]{\ensuremath{\boldsymbol{#1}}}

\renewcommand{\vec}[1]{\gbf{#1}}

\newcommand{\vv}{\ensuremath{\vec{v}}}

\usepackage{xspace}
\newcommand{\software}[1]{\texttt{#1}\xspace}
\newcommand{\hip}{\texttt{hIPPYlib}\xspace}
\newcommand{\bhip}{\texttt{{\bf hIPPYlib}}\xspace}
\newcommand{\Fe}{\texttt{FEniCS}\xspace}
\newcommand{\pet}{\texttt{PETSc}\xspace}

\newcommand{\del}{\partial}

\newcommand{\iFF}{\mathcal{F}}
\newcommand{\Acal}{\mc{A}}

\newcommand{\GN}{\ensuremath{\Gamma_{\!\!N}}}
\newcommand{\Exp}[1]{e^{{#1}}}
\newcommand{\GD}{\ensuremath{\Gamma_{\!\!D}}}
\newcommand{\Vg}{\mathscr{V}_{\!g}}
\newcommand{\V}{\mathscr{V}_{\!\scriptscriptstyle{0}}}

\newcommand{\LI}{\mathscr{L}}
\renewcommand{\H}{\mathcal{H}}

\newcommand{\ipart}{m_{\scriptscriptstyle\text{true}}}
\newcommand{\grad}{\nabla}
\newcommand{\ut}[1]{\ensuremath{\tilde{#1}}}

\newcommand{\LRp}[1]{\left( #1 \right)}

\newcommand{\LRc}[1]{\left\{ #1 \right\}}

\newcommand{\xx}{\ensuremath{\boldsymbol x}}

\newcommand{\x}{\xx}

\newcommand{\half} {\ensuremath{\frac{1}{2}}}
\newcommand{\nor}[1]{\left\| #1 \right\|}
\newcommand{\equaldef}{{:=}}
\newcommand{\iFFadj}{\adjMacroMM{\iFF}}

%% file: tikz/spec-and-evecs-possion.tex
\begin{tikzpicture}
  \begin{axis}[width=8cm, height=7cm, scale only axis, 
      xlabel = number, ylabel=eigenvalue, xmin=1.0, xmax=30,
      ymin=0.1, ymax=1e6, ymode = log, legend style={font=\small,nodes=right}, legend pos= north east]
    \addplot [color=black,line width=1pt] table[x=nr,y=eig]{extraplots/poisson-new/ones.dat};
    \addplot [color=blue, mark=*, only marks, mark size=1pt] table[x=nr,y=eig]{extraplots/poisson-new/eigenvalues.dat};
  \end{axis};
  \node[anchor=south west] at (9, -1.75)
       {
         {\renewcommand{\arraystretch}{1} 
           \begin{tabular}{c@{}c}
             \includegraphics[scale=.1]{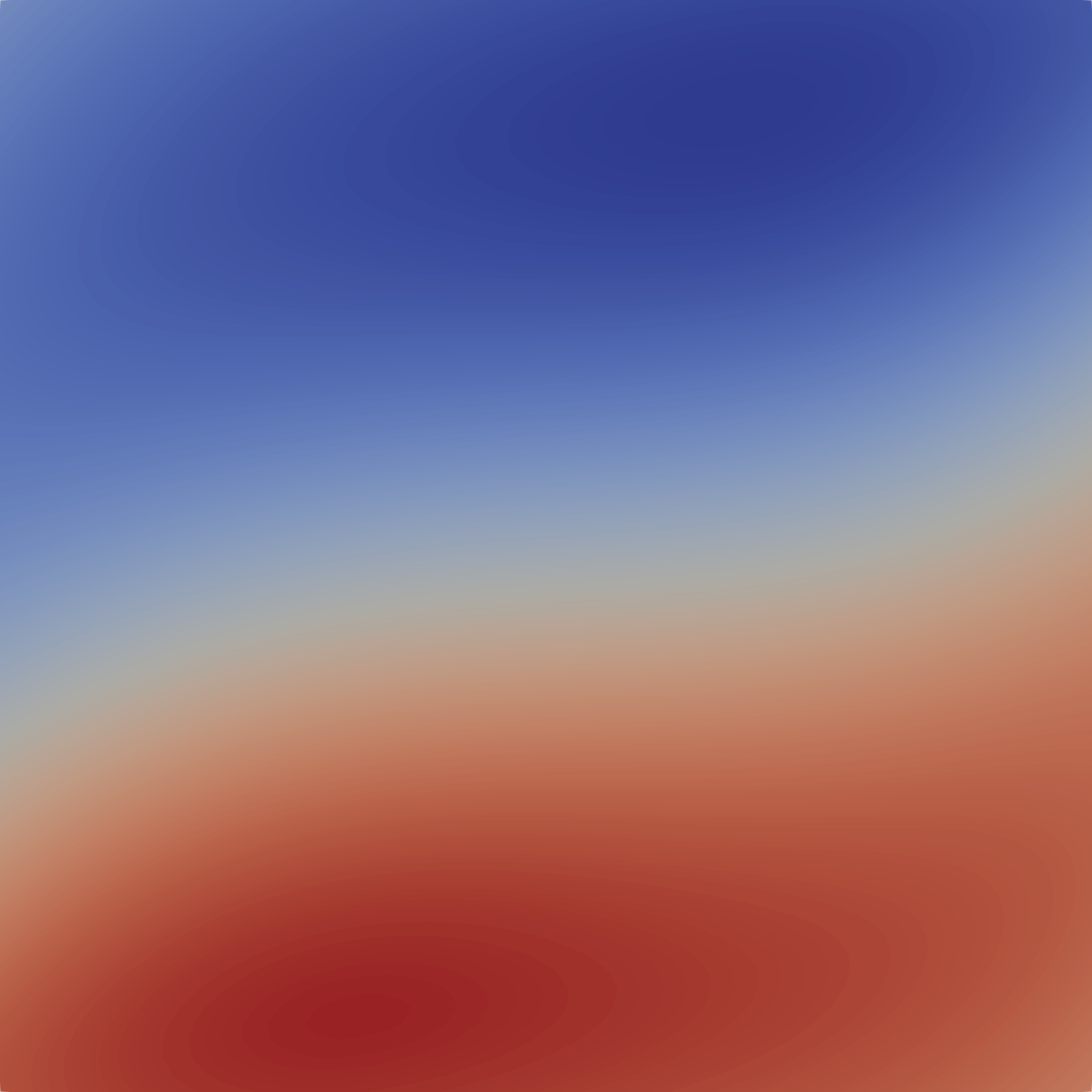} \hfill
             \includegraphics[scale=.1]{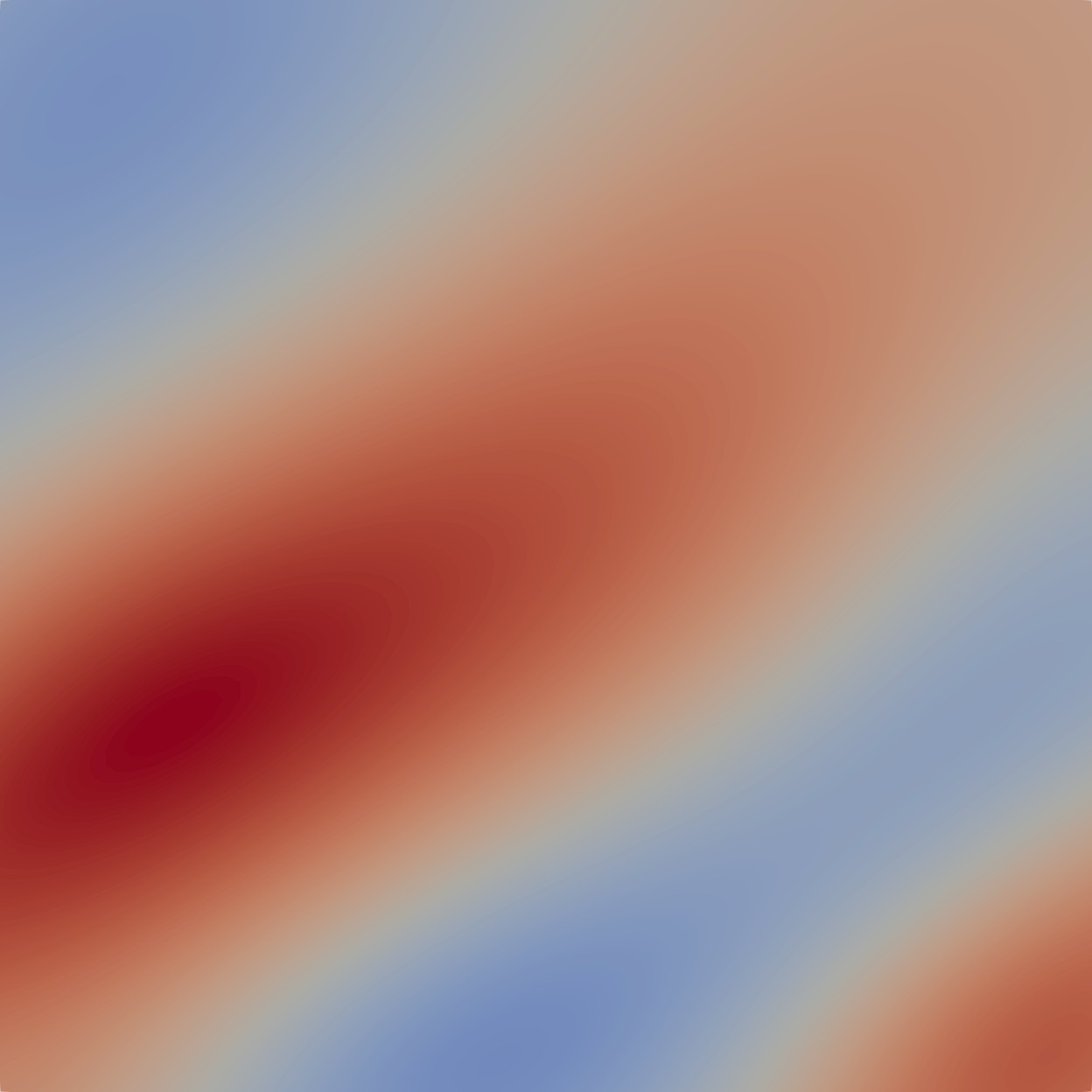} \\
             \includegraphics[scale=.1]{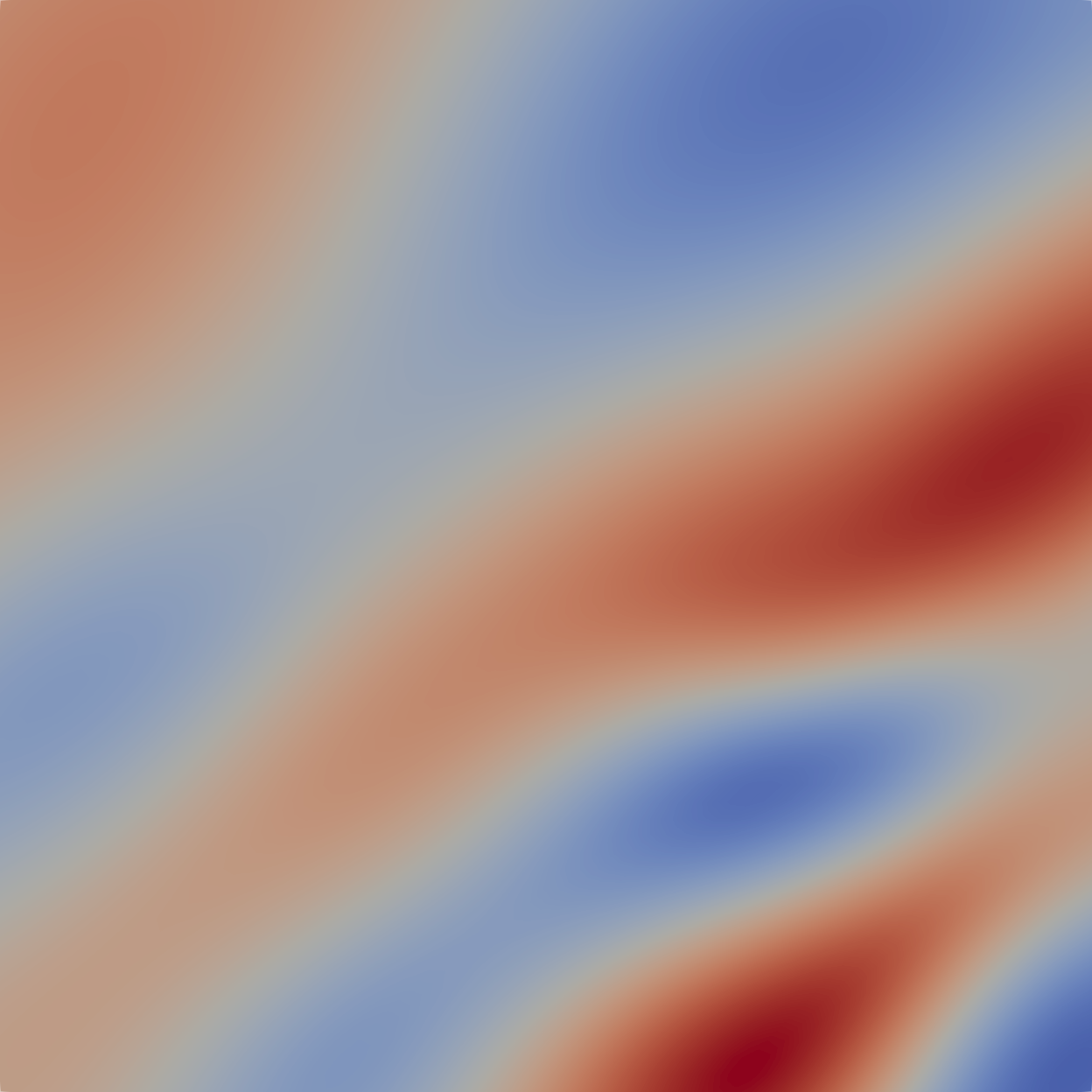} \hfill
             \includegraphics[scale=.1]{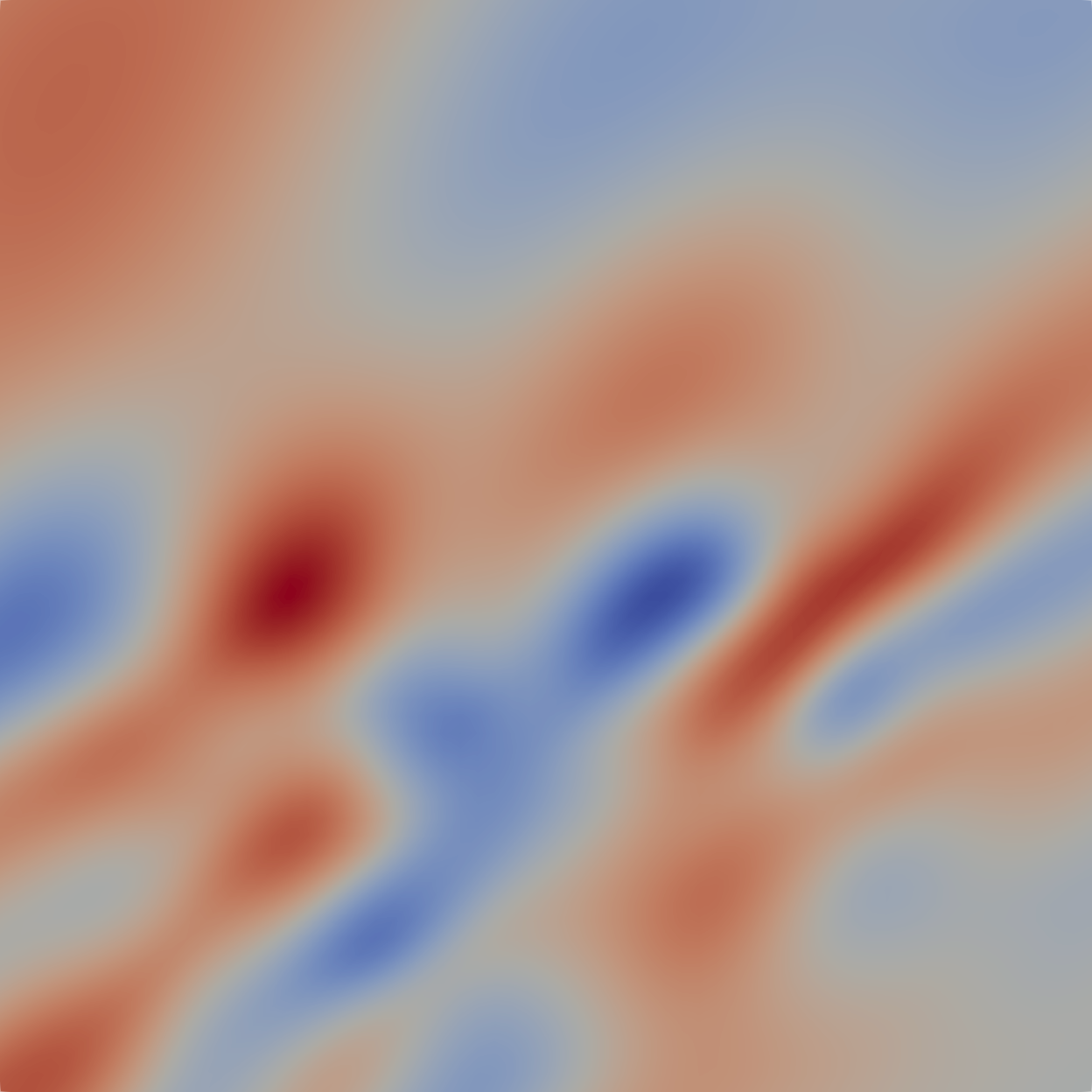}
           \end{tabular}
         }
       };
\end{tikzpicture}

%% file: tikz/spec-and-evecs-ad_diff.tex
\begin{tikzpicture}
  \begin{axis}[width=8cm, height=7cm, scale only axis, 
      xlabel = number, ylabel=eigenvalue, xmin=1.0, xmax=70,
      ymin=0.1, ymax=1e6, ymode = log, legend style={font=\small,nodes=right}, legend pos= north east]
    \addplot [color=blue, mark=*, only marks, mark size=1pt] table[x=nr,y=eig]{extraplots/ad_diff-new/eigenvalues.dat};
    \addlegendentry{$T_1 = 1$}
    \addplot [color=red, mark=*, only marks, mark size=1pt] table[x=nr,y=eig]{extraplots/ad_diff-new/eigenvalues_2to4.dat};
    \addlegendentry{$T_1 = 2$}
    \addplot [color=grassgreen, mark=*, only marks, mark size=1pt] table[x=nr,y=eig]{extraplots/ad_diff-new/eigenvalues_3to4.dat};
    \addlegendentry{$T_1 = 3$}
    \addplot [color=black,line width=1pt] table[x=nr,y=eig]{extraplots/ad_diff-new/ones.dat};
  \end{axis};
  \node[anchor=south west] at (9, -1.75)
       {
         {\renewcommand{\arraystretch}{1} 
           \begin{tabular}{c@{}c}
             \includegraphics[scale=.1]{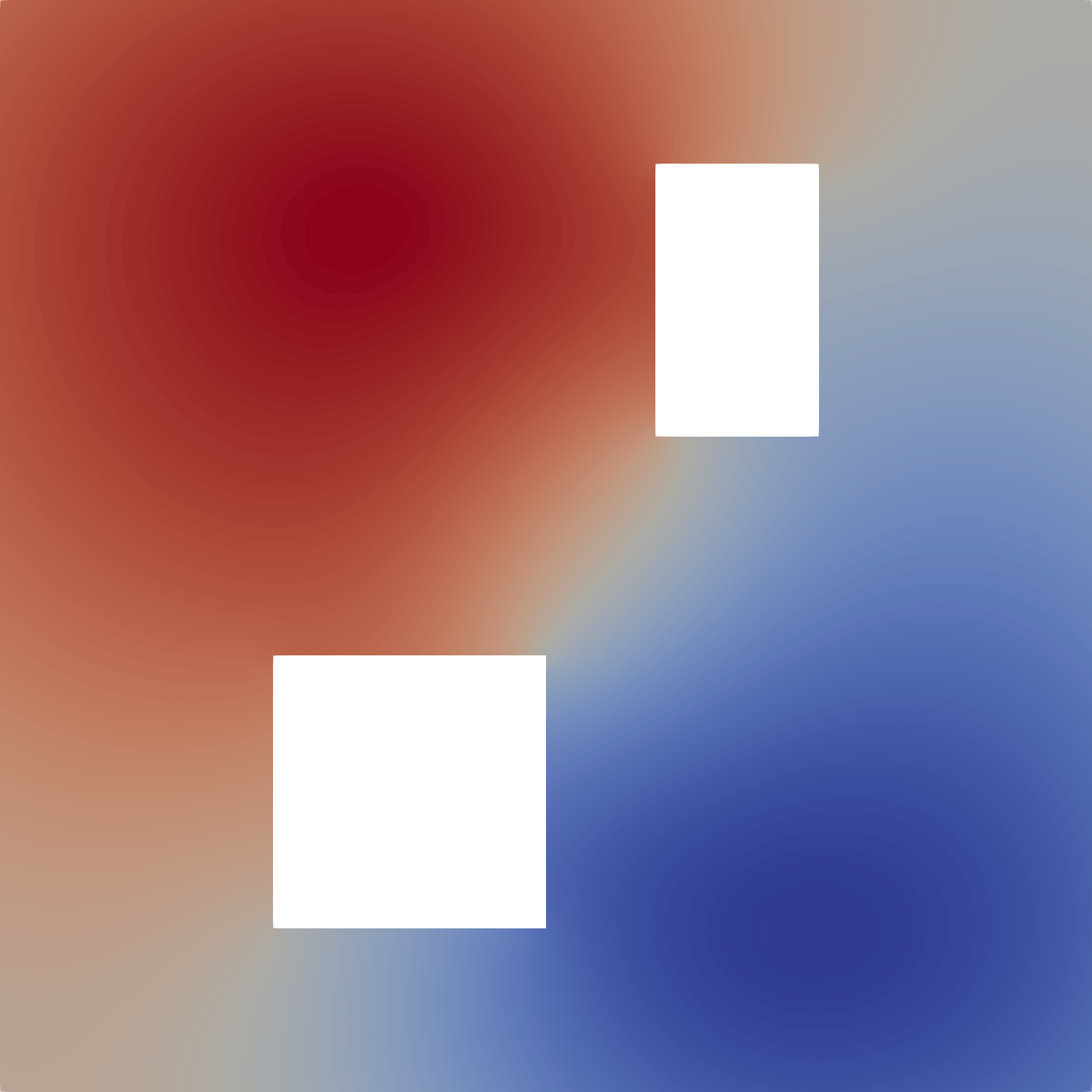} \hfill
             \includegraphics[scale=.1]{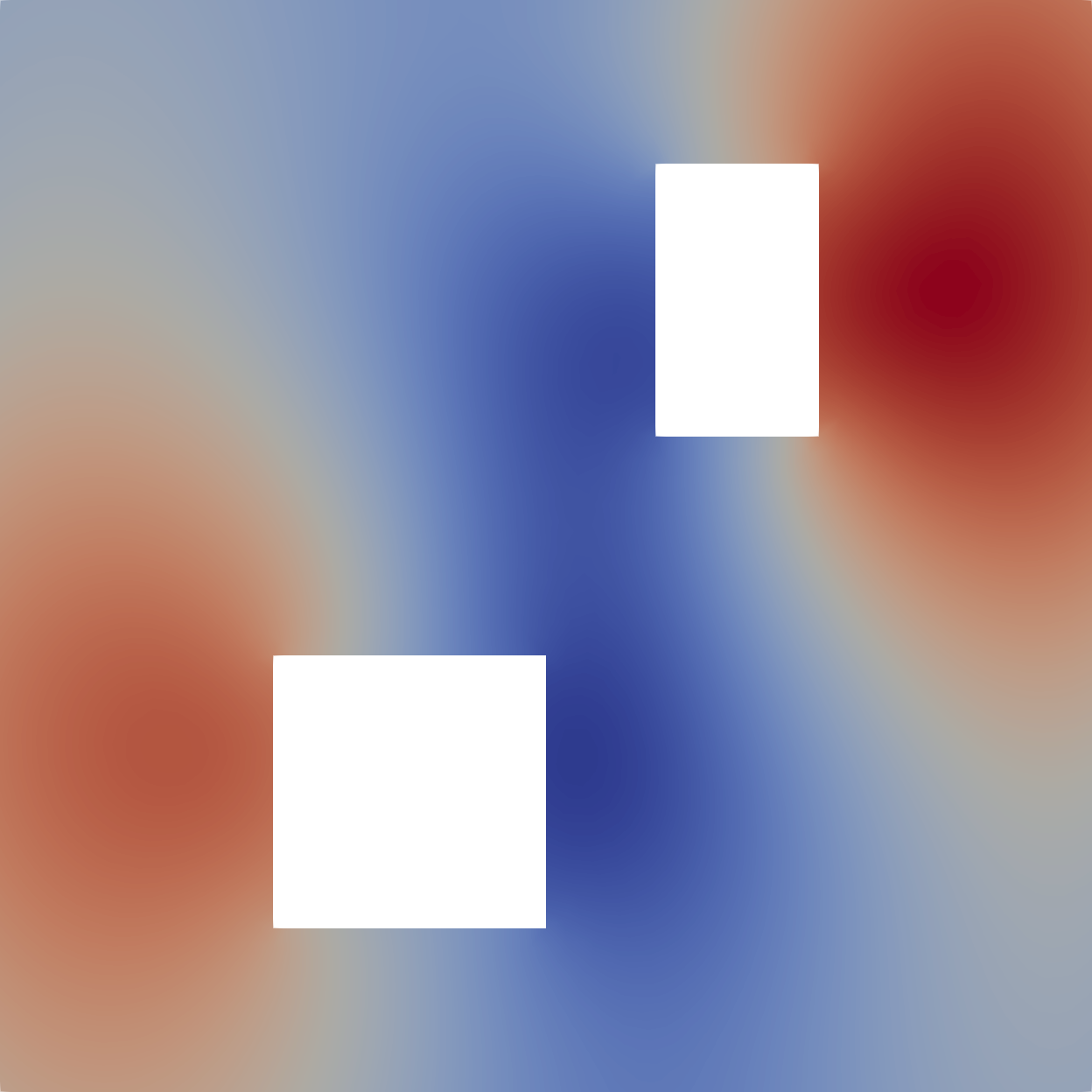} \\
             \includegraphics[scale=.1]{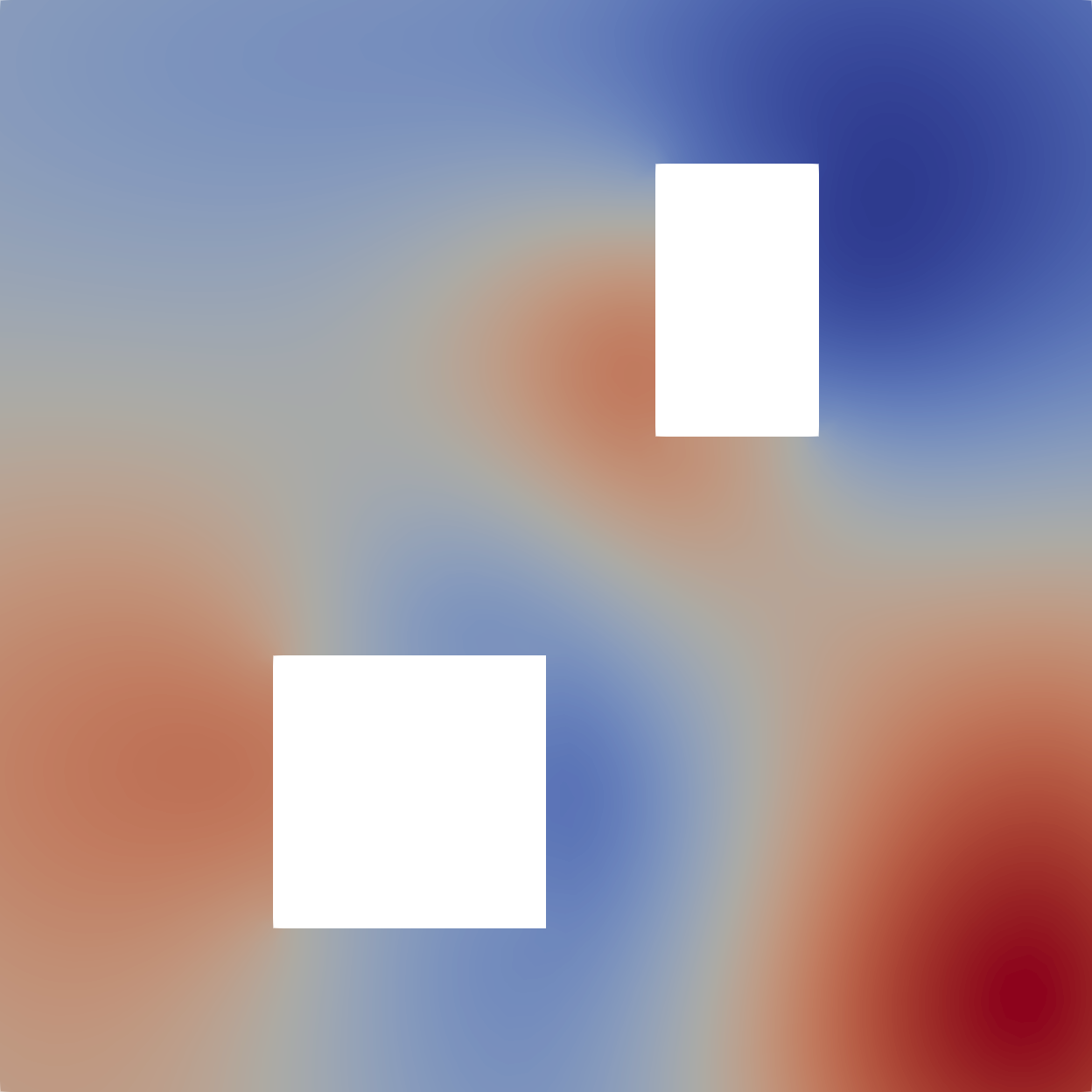} \hfill
             \includegraphics[scale=.1]{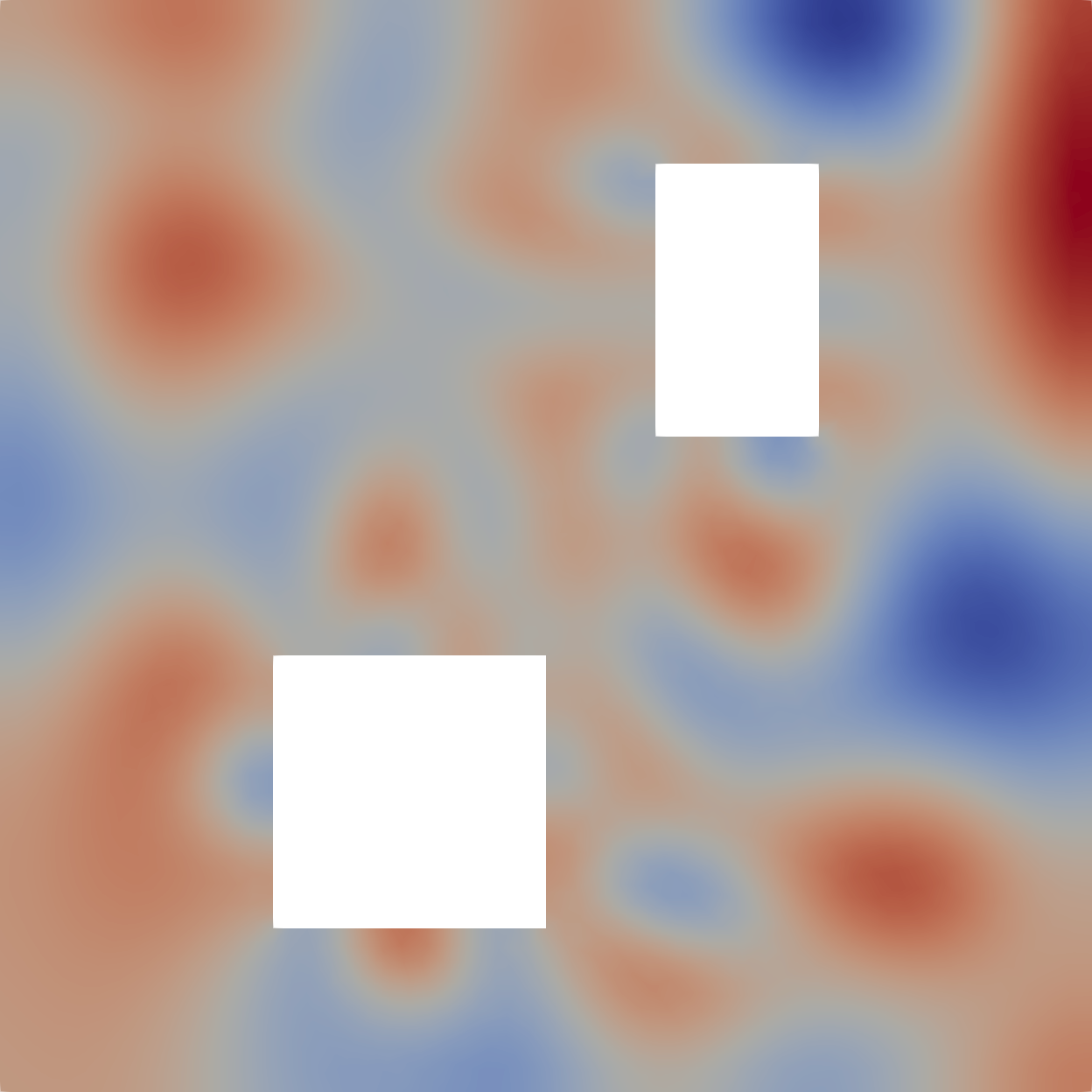}
           \end{tabular}
         }
       };
\end{tikzpicture}